\documentclass[11pt]{article}
\usepackage{amssymb}

\newtheorem{definition}{Definition}[section]
\newtheorem{theorem}[definition]{Theorem}
\newtheorem{lemma}[definition]{Lemma}
\newtheorem{corollary}[definition]{Corollary}

\newtheorem{example}[definition]{Example}
\newtheorem{note}[definition]{Note}

\def\R{{\rm I\! R}}
\def\C{{\rm I}\!\!\! {\rm C}}

\newcommand{\alg}{{\cal A}}
\newcommand{\fld}{{\cal F}}
\newcommand{\Mdf}{{\hbox{Mat}_{d+1}({\cal F})}}
\newcommand{\ls}{{\Phi = (A;
 E_0,E_1,\ldots,E_d;
A^*; E^*_0,E^*_1,\ldots,E^*_d)}}
\newcommand{\Span}[1]{{ \mbox{\rm Span}\{#1\}}} 
\newcommand{\beast}{\begin{eqnarray*}}
\newcommand{\eeast}{\end{eqnarray*}}

\begin{document}
\newenvironment{proof}{\noindent{\it Proof\/}:}{\par\noindent $\Box$\par}

\title{ \bf Two linear transformations each tridiagonal \\
with respect
to an eigenbasis of the other} 
\author{Paul Terwilliger\footnote{
 Mathematics  Department \ \ 
University of  Wisconsin  \ \ 
 480 Lincoln Drive  \ \ 
Madison, WI 53706  \ \
Email:terwilli@math.wisc.edu }}

\date{}
\maketitle
\begin{abstract}
Let $\fld $ denote a field, and let $V$ denote a  
vector space over $\fld$ with finite positive dimension.
We consider a pair of 
 linear transformations 
$A:V\rightarrow V$ and $A^*:V\rightarrow V$ satisfying
both conditions below:
\begin{enumerate}
\item There exists a basis for $V$ with respect to which
the matrix representing $A$ is diagonal, and the matrix
representing $A^*$ is irreducible tridiagonal.
\item There exists a basis for $V$ with respect to which
the matrix representing $A^*$ is diagonal, and the matrix
representing $A$ is irreducible tridiagonal.
\end{enumerate}

\medskip
\noindent 
We call such a pair a {\it Leonard pair} 
on $V$. Refining this notion a bit, we introduce the concept
of a {\it Leonard system}. 
We give a complete classification of Leonard systems.
Integral to our proof is the following result.
We show that for any Leonard pair $A,A^*$ on $V$,
there exists a sequence of scalars $\beta, \gamma, \gamma^*,
\varrho, \varrho^*$ taken from $\fld$ such that
both
\beast
0 &=&\lbrack A,A^2A^*-\beta AA^*A + 
A^*A^2 -\gamma (AA^*+A^*A)-\varrho A^*\rbrack, 
\\
0 &=& \lbrack A^*,A^{*2}A-\beta A^*AA^* + AA^{*2} -\gamma^* (A^*A+AA^*)-
\varrho^* A\rbrack,  
\eeast
where $\lbrack r,s\rbrack $ means $rs-sr$. 
The sequence is uniquely determined by the Leonard pair if the
dimension of $V$ is at least 4.
We conclude by showing how Leonard systems  
correspond to 
$q$-Racah  and related polynomials
from the Askey scheme.
\end{abstract}


\bigskip 
\def\leaderfill{\leaders\hbox to 1em{\hss.\hss}\hfill}
{\narrower {\narrower {\narrower {\narrower \small
                   \centerline{CONTENTS}

\medskip\noindent
1. Introduction                                      \leaderfill 2  \\
2. Some preliminaries                                 \leaderfill 9 \\
3. The split canonical form                          \leaderfill  10 \\
4. The primitive idempotents of a Leonard system       \leaderfill 15 \\
5. A formula for the $\varphi_i$                      \leaderfill 20 \\
6. The $D_4$ action                                    \leaderfill 20 \\
7. A result on reducibility                          \leaderfill 22 \\
8. Recurrent sequences                             \leaderfill 24 \\
9. Recurrent sequences in closed form              \leaderfill 25 \\
10. A sum                                       \leaderfill 27 \\
11. Some equations involving the split canonical form   \leaderfill 30 \\
12. Two polynomial equations for $A$ and $A^*$          \leaderfill 32 \\
13. Some vanishing products                       \leaderfill 37 \\
14. A classification of Leonard systems            \leaderfill 40 \\
15. Appendix: Leonard systems and  polynomials       \leaderfill 42 \\
References                                            \leaderfill 46 \\

\baselineskip=\normalbaselineskip \par} \par} \par} \par}


%
%
%
%
%
%
%
%
%
%
%
%
%
%
%
%
%
%
%
%


\section{Introduction}

\noindent
Throughout this paper, $\fld$ will denote
an arbitrary field. 

\medskip
\noindent We begin with the following situation
in linear
algebra.

\begin{definition}
\label{def:leonardpairtalkphilS99}
Let 
 $V$ denote a  
vector space over $\fld$ with finite positive dimension.
By a {\it Leonard pair} on $V$,
we mean an ordered pair $(A, A^*)$, where
$A:V\rightarrow V$ and $A^*:V\rightarrow V$ are linear transformations 
that 
 satisfy both (i), (ii) below. 
\begin{enumerate}
\item There exists a basis for $V$ with respect to which
the matrix representing $A^*$ is diagonal, and the matrix
representing $A$ is irreducible tridiagonal.
\item There exists a basis for $V$ with respect to which
the matrix representing $A$ is diagonal, and the matrix
representing $A^*$ is irreducible tridiagonal.

\end{enumerate}
(A tridiagonal matrix is said to be irreducible
whenever all entries immediately above and below the main
diagonal are nonzero).

\end{definition}


\begin{note} According to a common notational convention, for
a linear transformation $A$ the conjugate-transpose of $A$ is denoted
$A^*$. We emphasize we are not using this convention. In a Leonard
pair $(A,A^*)$, the linear transformations $A$ and $A^*$
are arbitrary subject
to (i),  (ii) above.
\end{note}

\noindent Our use of the name ``Leonard pair'' is motivated
by a connection to a theorem of Leonard \cite{Leodual}, \cite[p260]{BanIto} 
involving
the $q$-Racah  and related polynomials of the Askey scheme \cite{KoeSwa}.
For more information on this, we refer the reader to 
 Section 15.


\medskip
\noindent Here is an example of a Leonard pair.
Set 
$V={\fld}^4$ (column vectors), set 
\beast
A = 
\left(
\begin{array}{ c c c c }
0 & 3  &  0    & 0  \\
1 & 0  &  2   &  0    \\
0  & 2  & 0   & 1 \\
0  & 0  & 3  & 0 \\
\end{array}
\right), \qquad  
A^* = 
\left(
\begin{array}{ c c c c }
3 & 0  &  0    & 0  \\
0 & 1  &  0   &  0    \\
0  & 0  & -1   & 0 \\
0  & 0  & 0  & -3 \\
\end{array}
\right),
\eeast
and view $A$ and $A^*$  as linear transformations from $V$ to $V$.
We assume 
the characteristic of $\fld$ is not 2 or 3, to ensure
$A$ is irreducible.
Then $(A, A^*)$ is a Leonard
pair on $V$. 
Indeed, 
condition (i) in Definition
\ref{def:leonardpairtalkphilS99}
is satisfied by the basis for $V$
consisting of the columns of the 4 by 4 identity matrix.
To verify condition (ii), we display an invertible  matrix  
$P$ such that 
$P^{-1}AP$ is 
diagonal, and 
 such that 
$P^{-1}A^*P$ is
irreducible tridiagonal.
Put 
\beast
P = 
\left(
\begin{array}{ c c c c}
1 & 3  &  3    &  1 \\
1 & 1  &  -1    &  -1\\
1  & -1  & -1  & 1  \\
1  & -3  & 3  & -1 \\
\end{array}
\right).
\eeast
 By matrix multiplication $P^2=8I$, where $I$ denotes the identity,   
so $P^{-1}$ exists. Also by matrix multiplication,    
\begin{equation}
AP = PA^*.
\label{eq:apeqpasS99}
\end{equation}
Apparently
$P^{-1}AP$ equals $A^*$, and is therefore diagonal.
By (\ref{eq:apeqpasS99}), and since $P^{-1}$ is
a scalar multiple of $P$, we find
$P^{-1}A^*P$ equals $A$, and is therefore irreducible tridiagonal.  Now 
condition (ii) of  Definition 
\ref{def:leonardpairtalkphilS99}
is satisfied
by the basis for $V$ consisting of the columns of $P$.

\medskip
\noindent When working with a Leonard pair, 
it is often convenient to consider a closely related
and somewhat more abstract object, which we call
a {\it Leonard system}.
In order to define this, we first make an observation about
Leonard pairs.

%
\begin{lemma}
\label{lem:preeverythingtalkS99}
With reference to Definition
\ref{def:leonardpairtalkphilS99},
let $(A, A^*)$ denote a Leonard pair on $V$. Then
the eigenvalues of $A$ are distinct and contained in $\fld$.
Moreover, the eigenvalues of
$A^*$ are distinct and contained in $\fld$.
\end{lemma}

\begin{proof} Concerning $A$, recall by 
 Definition
\ref{def:leonardpairtalkphilS99}(ii) that there exists
a basis for $V$ consisting of eigenvectors for $A$. 
Consequently the eigenvalues of $A$ are all in $\fld$, and the minimal
polynomial of $A$ has no repeated roots. To show the eigenvalues
of $A$ are distinct, we show the minimal polynomial of $A$ has
degree equal to $\hbox{dim}\,V$. 
By Definition
\ref{def:leonardpairtalkphilS99}(i), there exists
a basis for $V$ with respect to which the matrix representing
$A$ is
 irreducible tridiagonal. Denote this matrix by $B$. On one hand,
  $A$ and $B$ have the same minimal polynomial.    
On the other hand,  
using the tridiagonal shape of $B$, we find  
$I, B, B^2, \ldots, B^d$ are linearly independent, where
$d=\hbox{dim}\,V-1$,  so the minimal  polynomial of $B$
 has  degree $d+1=\hbox{dim}\,V$.
We conclude the  
mininimal polynomial of $A$  has degree equal to 
  $\hbox{dim}\,V$,
so the eigenvalues of  $A$ are distinct.
We have now obtained our assertions about $A$,
and the case of $A^*$ is similar.
\end{proof}
\noindent To prepare for our definition of  a Leonard system,
we recall a few concepts from elementary linear algebra. 
Let $d$  denote  a  nonnegative
integer, 
and let $\Mdf$ denote the $\fld$-algebra consisting of all
$d+1$ by $d+1$ matrices with entries in $\fld$. We
index the rows and columns by $0,1,\ldots, d$.
For the rest of this paper,  $\;\alg\;$ will
denote an $\;\fld$-algebra
isomorphic to 
$\;\hbox{Mat}_{d+1}(\fld)$. 
Let $A$ denote an element of $\cal A$. By an {\it eigenvalue }
of $A$, we mean a root of the minimal polynomial of $A$.
The eigenvalues  of $A$ are contained in the algebraic closure of $\fld$.
The element $A $  will be called {\it
multiplicity-free} whenever it has $d+1$ distinct  eigenvalues,
all of which are 
in $\;\fld$.
Let $A$ denote a  multiplicity-free element of $\alg$.
Let $\theta_0, \theta_1, \ldots, \theta_d$ denote an ordering of 
the eigenvalues
of $A$, and for $0 \leq i \leq d$   put 
\beast
E_i = \prod_{{0 \leq  j \leq d}\atop
{j\not=i}} {{A-\theta_j I}\over {\theta_i-\theta_j}},
\eeast
where $I$ denotes the identity of $\cal A$.
By elementary linear algebra,
\begin{eqnarray}
&&AE_i = E_iA = \theta_iE_i \qquad \qquad  (0 \leq i \leq d),
\label{eq:primid1S99}
\\
&&
\quad E_iE_j = \delta_{ij}E_i \qquad \qquad (0 \leq i,j\leq d),
\label{eq:primid2S99}
\\
&&
\qquad \qquad \sum_{i=0}^d E_i = I.
\label{eq:primid3S99}
\end{eqnarray}
From this, 
one  finds  $E_0, E_1, \ldots, E_d$ is a  basis for the
subalgebra of $\alg$ generated by $A$.
We refer to $E_i$ as the {\it primitive idempotent} of
$A$ associated with $\theta_i$.
It is helpful to think of these primitive idempotents as follows. 
Let $V$ denote the irreducible left $\cal A$-module. Then
\begin{eqnarray}
V = E_0V + E_1V + \cdots + E_dV \qquad \qquad (\hbox{direct sum}).
\label{eq:VdecompS99}
\end{eqnarray}
For $0\leq i \leq d$, $E_iV$ is the (one dimensional) eigenspace of
$A$ in $V$ associated with the 
eigenvalue $\theta_i$, 
and $E_i$ acts  on $V$ as the projection onto this eigenspace.

\begin{definition}
\label{def:deflstalkS99}
Let $d$  denote  a  nonnegative
integer, let $\fld $ denote a field, 
and let $\alg$ 
denote an $\;\fld$-algebra isomorphic to 
$\hbox{Mat}_{d+1}(\fld)$. 
By a {\it Leonard system} in $\;\alg$, we mean a 
sequence 
\begin{equation}
\; \Phi = (A;\,E_0,\,E_1,\,\ldots,
\,E_d;\,A^*;\,E^*_0,\,E^*_1,\,\ldots,\,E^*_d)
\label{eq:ourstartingpt}
\end{equation}
 that satisfies  (i)--(v) below.
\begin{enumerate}
\item $A$,  $\;A^*\;$ are both multiplicity-free elements in $\;\alg$.
\item $E_0,\,E_1,\,\ldots,\,E_d\;$ is an ordering of the primitive 
idempotents of $\;A$.
\item $E^*_0,\,E^*_1,\,\ldots,\,E^*_d\;$ is an ordering of the primitive 
idempotents of $\;A^*$.
\item ${\displaystyle{
E_iA^*E_j = \cases{0, &if $\;\vert i-j\vert > 1$;\cr
\not=0, &if $\;\vert i-j \vert = 1$\cr}
\qquad \qquad 
(0 \leq i,j\leq d)}}$.
\item ${\displaystyle{
 E^*_iAE^*_j = \cases{0, &if $\;\vert i-j\vert > 1$;\cr
\not=0, &if $\;\vert i-j \vert = 1$\cr}
\qquad \qquad 
(0 \leq i,j\leq d).}}$
\end{enumerate}
We refer to $d$ as the {\it diameter} of $\Phi$, and say 
$\Phi$ is {\it over } $\fld$.  We sometimes write
${\cal A} = {\cal A}(\Phi)$, $\fld= \fld(\Phi)$. 
 For notational convenience, we set $E_{-1}=0$, $E_{d+1}=0$, $
E^*_{-1}=0$, $E^*_{d+1}=0$.

\end{definition}

\noindent 
To see the connection between Leonard pairs and 
Leonard systems, observe 
conditions (ii), (iv) above assert
that with respect to an appropriate  
basis consisting of eigenvectors for $A$, the matrix representing 
$A^*$ is irreducible tridiagonal.
Similarily, 
conditions (iii), (v) assert
that with respect to  an appropriate  
basis consisting of eigenvectors for $A^*$, the matrix representing 
$A$ is irreducible tridiagonal.

\medskip
\noindent 
A little later in this introduction, we will state our 
main results, which are Theorems
\ref{thm:statementinintroS99},
\ref{thm:newrepLSintroS99}, and
\ref{thm:lastchancedolangradyIntS99}.
For now, we mention some of the concepts that get used. 

\medskip
\noindent 
Let $\Phi$ 
denote the Leonard system in
(\ref{eq:ourstartingpt}), 
and let 
$\sigma :\alg \rightarrow {\cal A}'$ denote an isomorphism of
$\fld$-algebras. We write 
\begin{equation}
\Phi^{\sigma}:= 
(A^{\sigma};  E_0^{\sigma},E_1^{\sigma},\ldots, E_d^{\sigma};
A^{*\sigma}; E_0^{*\sigma}, 
 E_1^{*\sigma},\ldots,
E_d^{*\sigma}),
\label{eq:lsisoAS99}
\end{equation}
and observe 
$\Phi^{\sigma}$ 
is a Leonard  system in ${\cal A }'$.

\begin{definition}
\label{def:isolsS99o}
Let $\Phi$ 
and  
 $\Phi'$ 
denote Leonard systems over $\fld$.
 By an {\it isomorphism of Leonard  systems
 from $\Phi $ to $\Phi'$}, we mean an isomorphism of $\fld $-algebras
$\sigma :{\cal A}(\Phi) \rightarrow {\cal A}(\Phi')$ such 
that  $\Phi^\sigma = \Phi'$. 
The Leonard systems $\Phi $, $\Phi'$
are said to be {\it isomorphic} whenever there exists
an isomorphism of Leonard  systems from $\Phi $ to $\Phi'$. 
\end{definition}

\medskip
\noindent A given Leonard system  can be modified in  several
ways to get a new Leonard system. For instance, 
let $\Phi$ 
 denote the Leonard system in
(\ref{eq:ourstartingpt}), 
 and let $\alpha$, $\alpha^*$, $\beta$,
$\beta^*$ denote scalars in $\fld$ such that $\alpha \not=0$, $\alpha^*\not=0$.
Then
\beast
(\alpha A+\beta I; E_0, E_1, \ldots, E_d; 
\alpha^* A^*+\beta^* I; E^*_0, E^*_1, \ldots, E^*_d)
\eeast
is a Leonard system in $\alg$.
Also, 
\begin{eqnarray}
 \;\Phi^*&:=& (A^*; E^*_0,E^*_1,\ldots,E^*_d;A;E_0,E_1, \ldots,E_d),
\label{eq:lsdualS99}
\\
\Phi^{\downarrow}&:=& (A; E_0,E_1,\ldots,E_d;A^*;E^*_d,E^*_{d-1}, \ldots,E^*_0),
\label{eq:lsinvertS99}
\\
\Phi^{\Downarrow} 
&:=& (A; E_d,E_{d-1},\ldots,E_0;A^*;E^*_0,E^*_1, \ldots,E^*_d)
\label{eq:lsdualinvertS99}
\end{eqnarray}
are Leonard systems in $\alg$.  
 We refer to $\Phi^*$
(resp.  
 $\Phi^\downarrow$)   
(resp.  
 $\Phi^\Downarrow$) 
 as the  
{\it dual} 
(resp. {\it first inversion})
(resp.  {\it second inversion}) of  $\Phi$.
Viewing $*, \downarrow, \Downarrow$
as permutations on the set of all Leonard systems,
\begin{eqnarray}
&&\qquad \qquad \qquad  *^2 \;=\;  
\downarrow^2\;= \;
\Downarrow^2 \;=\;1,
\qquad \quad 
\label{eq:deightrelationsAS99}
\\
&&\Downarrow *\; 
=\;
* \downarrow,\qquad \qquad   
\downarrow *\; 
=\;
* \Downarrow,\qquad \qquad   
\downarrow \Downarrow \; = \;
\Downarrow \downarrow.
\qquad \quad 
\label{eq:deightrelationsBS99}
\end{eqnarray}
The group generated by symbols 
$*, \downarrow, \Downarrow $ subject to the relations
(\ref{eq:deightrelationsAS99}),
(\ref{eq:deightrelationsBS99})
is the dihedral group $D_4$.  
We recall $D_4$ is the group of symmetries of a square,
and has 8 elements.
Apparently $*, \downarrow, \Downarrow $ induce an action of 
 $D_4$ on the set of all Leonard systems.
Two Leonard systems will be called {\it relatives} whenever they
are in the same orbit of this $D_4$ action.
Assuming $d\geq 1$ to avoid
trivialities, the relatives of $\Phi$ are as follows:
\medskip

\centerline{
\begin{tabular}[t]{c|c}
        name &relative \\ \hline 
        $\Phi$ & $(A;E_0,E_1,\ldots,E_d;A^*;E^*_0,E^*_1,\ldots,E^*_d)$   \\ 
        $\Phi^\downarrow$ &
         $(A;E_0,E_1,\ldots,E_d;A^*;E^*_d,E^*_{d-1},\ldots,E^*_0)$   \\ 
        $\Phi^\Downarrow$ &
         $(A;E_d,E_{d-1},\ldots,E_0;A^*;E^*_0,E^*_1,\ldots,E^*_d)$   \\ 
        $\Phi^{\downarrow \Downarrow}$ &
         $(A;E_d,E_{d-1},\ldots,E_0;A^*;E^*_d,E^*_{d-1},\ldots,E^*_0)$   \\ 
	$\Phi^*$ &
        $(A^*;E^*_0,E^*_1,\ldots,E^*_d;A;E_0,E_1,\ldots,E_d)$   \\ 
        $\Phi^{\downarrow *}$ &
	 $(A^*;E^*_d,E^*_{d-1},\ldots,E^*_0;  
         A;E_0,E_1,\ldots,E_d)$ \\
        $\Phi^{\Downarrow *}$ &
	 $(A^*;E^*_0,E^*_1,\ldots,E^*_d;    
         A;E_d,E_{d-1},\ldots,E_0)$ \\ 
	$\Phi^{\downarrow \Downarrow *}$ &
	 $(A^*;E^*_d,E^*_{d-1},\ldots,E^*_0;    
         A;E_d,E_{d-1},\ldots,E_0)$
	\end{tabular}}
\medskip
\noindent 
We remark there may be some isomorphisms among the above Leonard
systems.

\noindent In view of our above comments, 
when we discuss a Leonard system, we are often not
interested in the orderings of the primitive idempotents, we
just care how the $A$ and $A^*$ elements interact. This brings us back
to the notion of a 
 Leonard pair.

\begin{definition}
\label{def:leonardpairS99}
Let $d$  denote  a  nonnegative
integer, let $\fld$ denote a field, 
and let $\alg$ 
denote an $\;\fld$-algebra isomorphic to 
$\hbox{Mat}_{d+1}(\fld)$. 
By a {\it Leonard pair} in $\alg$, we mean an ordered
pair $(A, A^*)$ satisfying both (i), (ii) below.
\begin{enumerate}
\item $A$ and $A^*$ are both multiplicity-free elements in $\alg$.
\item There exists an ordering 
 $ E_0,E_1,\ldots,E_d
$ of the primitive idempotents of $A$, and there exists an
ordering 
$E^*_0,E^*_1,\ldots,E^*_d$ of the primitive idempotents of $A^*$,
such that 
\beast
(A;E_0,E_1,\ldots,E_d;A^*;E^*_0,E^*_1,\ldots,E^*_d)
\eeast
 is a 
Leonard system.
\end{enumerate}
We refer to $d$ as the diameter of the pair, and say the pair is
over $\fld$.
\end{definition}

\noindent     
Let $\Phi$ denote the Leonard system in  
(\ref{eq:ourstartingpt}). Apparently
the pair $(A,A^*)$ from that line is
a Leonard pair in $\alg$, which we say is 
{\it associated } with $\Phi$. 
 Let $(A,A^*)$ denote   
a Leonard pair in $\alg$.
Then
$(A^*, A)$ is a Leonard pair in $\alg$, which we call the
{\it dual} of $(A,A^*)$.
It is routine to show two Leonard systems are relatives
if and only if their associated Leonard pairs are equal
or dual.

\medskip
\noindent In the following lemma, we make explicit the connection
between the notions of Leonard pair that appear in 
Definition \ref{def:leonardpairtalkphilS99} and
Definition \ref{def:leonardpairS99}. The proof
is routine and left to the reader.

\begin{lemma}
\label{lem:lponvandlpinaS99}
Let $V$  denote a vector space
over $\fld$ with finite positive dimension. 
Let $\hbox{End(V)}$ denote the
$\fld$-algebra consisting of all linear transformations 
from $V$ to $V$, and recall 
 $\hbox{End(V)}$ is $\fld$-algebra isomorphic to
 $\Mdf$, where $d+1=\hbox{dim}(V)$. Then for 
 all $A$ and $A^*$ in  
 $\hbox{End(V)}$, the following are equivalent.
\begin{enumerate}
\item $(A,A^*)$ is a Leonard pair on $V$, in the sense of Definition
\ref{def:leonardpairtalkphilS99}.
\item $(A,A^*)$ is a Leonard pair in 
 $\hbox{End(V)}$, in the sense of Definition
\ref{def:leonardpairS99}.
\end{enumerate}
\end{lemma}

\noindent


\medskip 
\noindent We now introduce four sequences of parameters
that we will use to describe a given Leonard system.
The first two sequences  are given 
in the following definition.

\begin{definition}
\label{def:lsdefcommentsS99}
Let $\Phi$ denote the Leonard system in 
(\ref{eq:ourstartingpt}).
For $0 \leq i \leq d$, 
we let $\theta_i $ (resp. $\theta^*_i$) denote the eigenvalue
of $A$ (resp. $A^*$) associated with $E_i$ (resp. $E^*_i$).
We refer to  $\theta_0, \theta_1, \ldots, \theta_d$ as the 
eigenvalue sequence of $\Phi$.
We refer to  $\theta^*_0, \theta^*_1, \ldots, \theta^*_d$ as the 
dual eigenvalue sequence of $\Phi$.
\end{definition}
\noindent There are two more parameter sequences of interest to us. 
Let $\Phi$ denote the Leonard system in 
(\ref{eq:ourstartingpt}). As we will show in Theorem 3.2, 
there exists an 
isomorphism of $\fld$-algebras
$\flat :\alg \rightarrow  \Mdf$ 
and there exists  scalars
$\varphi_1, \varphi_2,\ldots, \varphi_d$ in $\fld$
such that
\beast
A^\flat = 
\left(
\begin{array}{c c c c c c}
\theta_0 & & & & & {\bf 0} \\
1 & \theta_1 &  & & & \\
& 1 & \theta_2 &  & & \\
& & \cdot & \cdot &  &  \\
& & & \cdot & \cdot &  \\
{\bf 0}& & & & 1 & \theta_d
\end{array}
\right),
&&\quad 
A^{*\flat} = 
\left(
\begin{array}{c c c c c c}
\theta^*_0 &\varphi_1 & & & & {\bf 0} \\
 & \theta^*_1 & \varphi_2 & & & \\
&  & \theta^*_2 & \cdot & & \\
& &  & \cdot & \cdot &  \\
& & &  & \cdot & \varphi_d \\
{\bf 0}& & & &  & \theta^*_d
\end{array}
\right),
\eeast
where the $\theta_i, \theta^*_i$ are from
Definition \ref{def:lsdefcommentsS99}.
The sequence 
$\flat, \varphi_1, \varphi_2,\ldots, \varphi_d$ is uniquely
determined by $\Phi$.
We refer to  
$\varphi_1, \varphi_2,\ldots, \varphi_d$ as the 
{\it $\varphi$-sequence} of $\Phi$. 
We let $\phi_1, \phi_2, \ldots, \phi_d$ denote the
$\varphi$-sequence of $\Phi^\Downarrow $, and call this 
the  
{\it $\phi$-sequence} of $\Phi $.  

\medskip
\noindent The central result of this paper is 
the following classification of Leonard systems.

\begin{theorem}
\label{thm:statementinintroS99}  
\label{thm:newrepLScharagainS99} Let 
$d$ denote a nonnegative integer, let $\fld $ denote a field,
and let 
\begin{eqnarray}
&&\theta_0, \theta_1, \ldots, \theta_d; \qquad \qquad \; 
\theta^*_0, \theta^*_1, \ldots, \theta^*_d; 
\label{eq:paramlist1S99}
\\
&&\varphi_1, \varphi_2, \ldots, \varphi_d;  \qquad \qquad 
\phi_1, \phi_2, \ldots, \phi_d \qquad \quad
\label{eq:paramlist2S99}
\end{eqnarray}
denote scalars in $\fld$. 
Then there exists  a Leonard system $\Phi$ over $\fld$  with 
eigenvalue sequence $\theta_0, \theta_1, \ldots, \theta_d$, 
dual eigenvalue sequence  
$\theta^*_0, \theta^*_1, \ldots, \theta^*_d$, $\varphi$-sequence
$ \varphi_1, \varphi_2, \ldots, \varphi_d $, and $\phi$-sequence
$\phi_1, \phi_2, \ldots, \phi_d$ if and only if 
(i)--(v) hold below.
\begin{enumerate}
\item $ \varphi_i \not=0, \qquad \phi_i\not=0 \qquad \qquad \qquad\qquad (1 \leq i \leq d)$,
\item $ \theta_i\not=\theta_j,\qquad  \theta^*_i\not=\theta^*_j\qquad $
if $\;\;i\not=j,\qquad \qquad \qquad (0 \leq i,j\leq d)$,
\item $ {\displaystyle{ \varphi_i = \phi_1 \sum_{h=0}^{i-1}
{{\theta_h-\theta_{d-h}}\over{\theta_0-\theta_d}} 
\;+\;(\theta^*_i-\theta^*_0)(\theta_{i-1}-\theta_d) \qquad \;\;(1 \leq i \leq d)}}$,
\item $ {\displaystyle{ \phi_i = \varphi_1 \sum_{h=0}^{i-1}
{{\theta_h-\theta_{d-h}}\over{\theta_0-\theta_d}} 
\;+\;(\theta^*_i-\theta^*_0)(\theta_{d-i+1}-\theta_0) \qquad (1 \leq i \leq d)}}$,
\item The expressions
\begin{equation}
{{\theta_{i-2}-\theta_{i+1}}\over {\theta_{i-1}-\theta_i}},\qquad \qquad  
 {{\theta^*_{i-2}-\theta^*_{i+1}}\over {\theta^*_{i-1}-\theta^*_i}} 
 \qquad  \qquad 
\label{eq:defbetaplusoneS99int}
\end{equation} 
 are equal and independent of $i$, for $\;2\leq i \leq d-1$.  
\end{enumerate}
Moreover, if (i)--(v) hold 
above then $\Phi$ is unique up to isomorphism of Leonard systems.
\end{theorem}
\noindent The proof of Theorem 
\ref{thm:statementinintroS99}  
appears in Section 14.

\medskip
\noindent We view Theorem
\ref{thm:statementinintroS99} as a linear algebraic version
of a theorem of Leonard
 \cite{Leodual}, \cite[p260]{BanIto}. This is discussed in Section 15.

\medskip
\noindent We have found all solutions to Theorem
\ref{thm:statementinintroS99}(i)--(v) in parametric form. We will
present these in a future paper, and for now display only the
``most general'' solution. It is
\beast
\theta_i &=& \theta_0 + h(1-q^i)(1-sq^{i+1})/q^i,
\\
\theta^*_i &=& \theta^*_0 + h^*(1-q^i)(1-s^*q^{i+1})/q^i
\eeast
for $0 \leq i \leq d$, and
\beast
\varphi_i &=& hh^*q^{1-2i}(1-q^i)(1-q^{i-d-1})(1-r_1q^i)(1-r_2q^i),
\\
\phi_i &=& hh^*q^{1-2i}(1-q^i)(1-q^{i-d-1})(r_1-s^*q^i)(r_2-s^*q^i)/s^*
\eeast
for $1 \leq i \leq d$, where $q, h, h^*, r_1, r_2, s, s^*$ are scalars
in the algebraic closure of $\fld$ such that $r_1r_2 = s s^*q^{d+1}$. 
For this solution the common value of 
(\ref{eq:defbetaplusoneS99int}) equals $q+q^{-1}+1$.

\medskip
\noindent  One nice feature of the parameter sequences 
(\ref{eq:paramlist1S99}), 
(\ref{eq:paramlist2S99}) is that they are modified 
in a simple way as one passes from a given Leonard system
to a relative. To describe how this works, we use the following
notational convention.

\medskip
\begin{definition}
\label{def:d4actionlsS99} Let $\Phi$ denote a Leonard system.
For any element $g$ of the group $D_4$,
and for 
 any object 
$f$ we  associate with $\Phi$,    
we let $f^g$ denote the corresponding object
associated with the Leonard system
$\Phi^{g^{-1}}$. We have been using this convention all along;
an example is $\theta^*_i(\Phi)= \theta_i(\Phi^*)$.

\end{definition}

%

\begin{theorem}
\label{thm:newrepLSintroS99} 
\label{thm:newrepLSS99} 
Let $\Phi$ denote a Leonard system,
with eigenvalue sequence
$\theta_0, \theta_1, \ldots, \theta_d$, 
dual eigenvalue sequence  
$\theta^*_0, \theta^*_1, \ldots, \theta^*_d$, $\varphi$-sequence
$ \varphi_1, \varphi_2, \ldots, \varphi_d $, and $\phi$-sequence
$\phi_1, \phi_2, \ldots, \phi_d$.
 Then for all $g \in D_4$, 
the scalars $\theta^g_i$, 
 $\theta^{*g}_i$, 
 $\varphi^{g}_i$, 
 $\phi^{g}_i$ are as follows.

\medskip
\centerline{
\begin{tabular}[t]{|c||c|c|c|c|}\hline
        g &$\theta^g_i$&$\theta^{*g}_i$&$\varphi^g_i$&$\phi^g_i$\\ \hline \hline
        $1$ &$\theta_i$&$\theta^*_i$&$\varphi_i$&$\phi_i$\\ \hline
        $\downarrow$ &$\theta_i$&$\theta^*_{d-i}$&
	                 $\phi_{d-i+1}$&$\varphi_{d-i+1}$\\ \hline
        $\Downarrow$ &$\theta_{d-i}$&$\theta^*_i$&$\phi_i$&$\varphi_i$\\ \hline
        ${\downarrow \Downarrow}$
	 &$\theta_{d-i}$&$\theta^*_{d-i}$&$\varphi_{d-i+1}$&$\phi_{d-i+1}$\\ \hline
	$*$ &$\theta^*_i$&$\theta_i$&$\varphi_i$&$\phi_{d-i+1}$\\ \hline
        ${\downarrow *}$
            	&$\theta^*_i$&$\theta_{d-i}$&$\phi_i$&$\varphi_{d-i+1}$\\ \hline
        ${\Downarrow *}$ &$\theta^*_{d-i}$&$\theta_i$&
	                 $\phi_{d-i+1}$&$\varphi_i$\\ \hline
        ${\downarrow \Downarrow *}$
	 &$\theta^*_{d-i}$&$\theta_{d-i}$&$\varphi_{d-i+1}$&$\phi_i$\\ \hline
	\end{tabular}}
\medskip
\end{theorem}

\noindent The proof of 
Theorem 
\ref{thm:newrepLSintroS99}  appears in Section 6.

\medskip
\noindent  We show the elements of a Leonard pair
satisfy the following relations.

\begin{theorem}
\label{thm:lastchancedolangradyIntS99}
\label{eq:lastchancedolangradyS99}
Let $\fld$ denote
a field,
and let $(A,A^*)$ denote a Leonard pair over $\fld$.
 Then
there  exists a sequence of scalars $\beta, \gamma, \gamma^*,
\varrho, \varrho^*$ taken from $\fld$ such that
both
\begin{eqnarray}
0 &=&\lbrack A,A^2A^*-\beta AA^*A + 
A^*A^2 -\gamma (AA^*+A^*A)-\varrho A^*\rbrack, 
\label{eq:qdolangrady199nint}
\label{eq:qdolangrady199n}
\\
0 &=& \lbrack A^*,A^{*2}A-\beta A^*AA^* + AA^{*2} -\gamma^* (A^*A+AA^*)-
\varrho^* A\rbrack,  \qquad  \quad 
\label{eq:qdolangrady2S99nint}
\label{eq:qdolangrady2S99n}
\end{eqnarray}
where $\lbrack r,s \rbrack $ means $\,rs-sr$.
The sequence is uniquely determined by the Leonard pair if the
diameter is at least 3.
\end{theorem}
\noindent The proof of 
Theorem \ref{thm:lastchancedolangradyIntS99} appears at the
end of Section 12.

\medskip
\noindent The relations 
(\ref{eq:qdolangrady199nint}), 
(\ref{eq:qdolangrady2S99nint})
previously appeared in \cite{TersubIII}.
In that paper, the author
considers a combinatorial  object called  
a $P$- and $Q$-polynomial association scheme
\cite{BanIto}, \cite{bcn},
\cite{Leopandq},
\cite{Tercharpq},
\cite{Ternew}.
He shows 
 that for these schemes the adjacency matrix $A$ 
 and a certain diagonal matrix  $A^*$
satisfy
(\ref{eq:qdolangrady199nint}), 
(\ref{eq:qdolangrady2S99nint}). 
In this context the algebra 
generated by $A$ and $A^*$ is known as the subconstituent
algebra or the Terwilliger algebra
\cite{Cau},
\cite{Col},
\cite{Curbip12},
\cite{Curbip2},
\cite{Curthin},
\cite{CurNom},
\cite{Go},
\cite{HobIto},
\cite{Tan},
\cite{TersubI},
\cite{TersubII}.

\medskip
\noindent A special case of 
(\ref{eq:qdolangrady199nint}),
(\ref{eq:qdolangrady2S99nint}) occurs in the context of quantum
groups. Setting $\beta = q^2 + q^{-2}$, 
$\gamma=0 $, $\gamma^*=0$, $\varrho=0$, $\varrho^*=0$   
in 
 (\ref{eq:qdolangrady199nint}),
(\ref{eq:qdolangrady2S99nint}), one obtains
\begin{eqnarray}
0 &=& A^3 A^* - \lbrack 3 \rbrack_q A^2A^*A +  
 \lbrack 3 \rbrack_q AA^*A^2 - A^*A^3,  
\label{eq:sSerre1}
\\
0 &=& A^{*3} A - \lbrack 3 \rbrack_q A^{*2}AA^* +  
 \lbrack 3 \rbrack_q A^*AA^{*2} - AA^{*3},
\label{eq:sSerre2}
\end{eqnarray}
where 
\beast
\lbrack 3 \rbrack_q := {{q^3 - q^{-3}}\over {q - q^{-1}}}.
\eeast
The equations 
(\ref{eq:sSerre1}), 
(\ref{eq:sSerre2}) are known as the $q$-{\it Serre relations}, and 
are among the defining relations for the quantum affine algebra
$U_q({\widehat {sl}}_2)$   
\cite{CP1}, \cite{CP2}. 

\medskip
\noindent
A special case of 
(\ref{eq:qdolangrady199nint}),
(\ref{eq:qdolangrady2S99nint}) has come up in the context of
exactly solvable models in  statistical mechanics.
Setting $\beta=2,\gamma=0, \gamma^*=0,\varrho=16, \varrho^*=16$
in 
 (\ref{eq:qdolangrady199nint}),
(\ref{eq:qdolangrady2S99nint}), one obtains
\begin{eqnarray}
\lbrack A, \lbrack A, \lbrack A,A^*\rbrack \rbrack \rbrack &=& 16
\lbrack A,A^*\rbrack,
\label{eq:DG1S99}
\\
\lbrack A^*, \lbrack A^*, \lbrack A^*,A\rbrack \rbrack \rbrack &=& 16
\lbrack A^*,A\rbrack.
\label{eq:DG2S99}
\end{eqnarray}
The equations 
(\ref{eq:DG1S99}), 
(\ref{eq:DG2S99}) 
are known as the Dolan--Grady
relations
\cite{CKOns},
\cite{Davsuper},
\cite{DolGra}, 
\cite{Ugl}.  
We remark the Lie algebra over $\C$ generated by two symbols $A, A^*$
subject to 
(\ref{eq:DG1S99}), 
(\ref{eq:DG2S99}) (where we interpret 
 $\lbrack \,,\,\rbrack $ as the Lie bracket)
 is infinite dimensional and is known as the 
Onsager algebra 
\cite{Dav}, \cite{Per}. 
The author would like to thank Anatol N. Kirillov for pointing out
the connection to statistical mechanics.


\medskip
\noindent 
We  mention the relations
(\ref{eq:qdolangrady199nint}),
(\ref{eq:qdolangrady2S99nint})  are satisfied by the
generators of 
both the classical
and quantum ``Quadratic  Askey-Wilson algebra''
introduced by   Granovskii, Lutzenko, and 
Zhedanov \cite{GYLZmut}.
See 
\cite{GYZnature},
\cite{GYZTwisted},
\cite{GYZlinear},
\cite{GYZspherical},
\cite{Zhidd}, 
\cite{ZheCart},
\cite{Zhidden}
for more information on this algebra.

\medskip
\noindent Given a field $\fld$, and given scalars 
$\beta, \gamma, \gamma^*,
\varrho, \varrho^*$ taken from $\fld$,
it is natural in light of our above comments
to consider the associative $\fld$-algebra generated
by two symbols $A$, $A^*$ subject to the relations
(\ref{eq:qdolangrady199nint}),
(\ref{eq:qdolangrady2S99nint}). It appears to be an open problem
to give a basis for this algebra, and to describe  all 
its irreducible representations.

\section{Some preliminaries }


\noindent We now turn to the business of proving 
the results that we displayed in the Introduction.
This will take most of the paper,
up through the end of Section 14.
We begin with some simple observations. 

\begin{definition}
\label{def:verygeneralS99}
In this section, $d$ will denote a nonnegative integer, $\fld $ will 
denote a field, and $\alg $ will denote an $\fld$-algebra isomorphic
to $\Mdf$. We let $V$ denote the irreducible left $\alg$-module. 
We let $A$ and $A^*$ denote  multiplicity-free elements of
$\alg$.  We let  
$E_0, E_1, \ldots E_d$  denote an  ordering
of the primitive idempotents of $A$, and we let
$E^*_0, E^*_1, \ldots E^*_d$ denote an  ordering
of the primitive idempotents of $A^*$.
For $0 \leq i \leq d$, we let 
 $\theta_i$ (resp. $\theta^*_i$) denote
the eigenvalue of $A$ (resp. $A^*$) associated
with $E_i$ (resp. $E^*_i$).
\end{definition}

\begin{definition}
\label{def:modmeaningS99}
With reference to Definition
\ref{def:verygeneralS99}, 
by an $(A, A^*)$-module, we mean  a subspace $W$ of $V$  
such that 
 $AW\subseteq W$  and 
 $A^*W\subseteq W$.
Let $W$ denote an $(A,A^*)$-module.
We say $W$ is irreducible
whenever $W\not=0$, and $W$ contains  
no $(A, A^*)$-modules other than $0$ and $W$. 
\end{definition}

\noindent With reference to Definition
\ref{def:verygeneralS99},
let $W$ denote an  $(A,A^*)$-module.
Since $AW\subseteq W$, 
we find
$W =\sum_{i\in S}E_iV$,
where $S$ is an appropriate subset of 
$\lbrace 0,1,\ldots, d\rbrace $. 
We now consider which subsets $S$ can occur.

\begin{lemma}
\label{lem:modcharS99}
With reference to 
Definition
\ref{def:verygeneralS99},
let $S$ denote a  subset of $\lbrace 0,1,\ldots, d \rbrace $,
and put 
$W = 
\sum_{i\in S} E_iV$.
Then the following are equivalent.
\begin{enumerate}
\item 
$W$ 
is an $(A,A^*)$-module.
\item 
$E_iA^*E_j=0 $ for all  
$j \in S$ and for all  $i \in 
 \lbrace 0,1,\ldots, d \rbrace \backslash S$.
\end{enumerate}
\end{lemma}

\begin{proof} $(i)\rightarrow (ii)$
Let the integers $i,j$ be given. 
Observe
$E_jV\subseteq W$ since $j\in S$, and $A^*W\subseteq W$, so 
$A^*E_jV\subseteq W$. Observe  $E_iW=0$ since $i \notin S$,
so 
$E_iA^*E_jV=0$. It follows  
$E_iA^*E_j=0$.   

\noindent 
$(ii)\rightarrow (i)$ 
Recall $E_iV$ is an eigenspace of $A$ for 
$0 \leq i \leq d$, so 
$AW\subseteq W$. We now show $A^*W\subseteq W$.
For notational convenience  set $J=\sum_{i\in S}E_i$.
Let ${\overline S} $ denote the complement of
$S$ in $\lbrace 0,1,\ldots, d\rbrace $, and set
 $K=\sum_{i\in {\overline S}}E_i$.
Then $J+K=I$ and $KA^*J=0$. Combining these,
we find $A^*J = JA^*J$. Applying this to $V$, and
observing $JV=W$, we routinely find $A^*W\subseteq W$.
\end{proof}

\noindent Interchanging the roles of $A$ and $A^*$ in
the previous lemma, we immediately  obtain the following result.

\begin{lemma}
\label{lem:modcharstarS99}
With reference to 
Definition
\ref{def:verygeneralS99},
let $S^*$ denote a  subset of $\lbrace 0,1,\ldots, d \rbrace $,
and put 
$W = 
\sum_{i\in S^*} E^*_iV.
$
Then the following are equivalent.
\begin{enumerate}
\item 
$W$ 
is an $(A,A^*)$-module.
\item
$E^*_iAE^*_j=0 $ for all  
$j \in S^*$ and for all  $i \in 
 \lbrace 0,1,\ldots, d \rbrace \backslash S^*$.
\end{enumerate}
\end{lemma}


\noindent The following  
 constants will be of use to us.

\begin{definition}
\label{def:aidefS99}
With reference to Definition 
\ref{def:verygeneralS99},
we define
\begin{eqnarray}
a_i = \hbox{tr}\, AE^*_i, \qquad  \quad
a^*_i = \hbox{tr}\, A^*E_i, \qquad \qquad  
 (0 \leq i \leq d), \qquad
\label{eq:defofaiS99}
\end{eqnarray}
where $tr$ means trace.
\end{definition}

\begin{lemma}
\label{lem:sumeigequalssumaiS99}
With reference to Definition 
\ref{def:verygeneralS99}
and Definition
\ref{def:aidefS99},
\begin{eqnarray}
\theta_0+\theta_1+\cdots + \theta_d &=& a_0+a_1+\cdots+ a_d,
\label{eq:sumeigequalssumaiAS99}
\\
\theta^*_0+\theta^*_1+\cdots + \theta^*_d &=& a^*_0+a^*_1+\cdots+ a^*_d.
\label{eq:sumeigequalssumaistarAS99}
\end{eqnarray}

\end{lemma}
\begin{proof} 
To get 
(\ref{eq:sumeigequalssumaiAS99}),  
take the trace of both sides in the equation
$A = A \sum_{i=0}^d E^*_i$,
and evaluate the result using   
the left equation in 
(\ref{eq:defofaiS99}).
Line (\ref{eq:sumeigequalssumaistarAS99}) is similarily obtained.

\end{proof}

\section{The split canonical form }

\begin{definition}
\label{def:splitconsetupS99}
In this section, $d$ will denote a nonnegative integer, $\fld $ will 
denote a field, and $\alg $ will denote an $\fld$-algebra isomorphic
to $\Mdf$. We let $V$ denote the irreducible left $\alg$-module.
We let 
\beast
\ls
\eeast
denote a Leonard system in $\alg$, with eigenvalue  sequence
$\theta_0, \theta_1, \ldots, \theta_d$ and dual eigenvalue 
sequence
 $\theta^*_0, \theta^*_1, \ldots, \theta^*_d$.
\end{definition}

\noindent 
With reference to Definition
\ref{def:splitconsetupS99}, when studying 
$\Phi$, it is tempting to represent
one of  $A$ and $A^*$ by a  
 diagonal matrix,  and the other by an irreducible
tridiagonal matrix.
This approach has some merit, but  we are going to
do something else.
Our goal in this section is to prove the following theorem.

\begin{theorem}
\label{thm:splitformexistS99}
With reference to Definition
\ref{def:splitconsetupS99}, there  
 exists  
an isomorphism of $\fld$-algebras 
$\flat :\alg \rightarrow  \Mdf$
and there exists 
scalars $\varphi_1, \varphi_2,\ldots, \varphi_d$ in $\fld$
such that
\begin{eqnarray}
A^\flat = 
\left(
\begin{array}{c c c c c c}
\theta_0 & & & & & {\bf 0} \\
1 & \theta_1 &  & & & \\
& 1 & \theta_2 &  & & \\
& & \cdot & \cdot &  &  \\
& & & \cdot & \cdot &  \\
{\bf 0}& & & & 1 & \theta_d
\end{array}
\right),
&&\quad 
A^{*\flat} = 
\left(
\begin{array}{c c c c c c}
\theta^*_0 &\varphi_1 & & & & {\bf 0} \\
 & \theta^*_1 & \varphi_2 & & & \\
&  & \theta^*_2 & \cdot & & \\
& &  & \cdot & \cdot &  \\
& & &  & \cdot & \varphi_d \\
{\bf 0}& & & &  & \theta^*_d
\end{array}
\right).
\label{eq:ahrtandastarhrtS99}
\end{eqnarray}
The sequence
$\flat, \varphi_1, \varphi_2,\ldots, \varphi_d$ is uniquely
determined by $\Phi$. Moreover 
 $\varphi_i\not=0$ for $1\leq i \leq d$.

\end{theorem}

\noindent 
 We begin with an irreducibility result.

\begin{lemma}
\label{lem:VirredaastarmoduleS99}
With reference to 
Definition
\ref{def:splitconsetupS99}, 
The module $V$ is irreducible as an $(A,A^*)$-module. 
\end{lemma}

\begin{proof}
Let $W$ denote a nonzero  $(A,A^*)$-module in $V$. We show
$W=V$.
 Since  $AW\subseteq W$, 
there exists a  subset $S$ of $\lbrace 0,1,\ldots, d\rbrace
$ such that $W=\sum_{h\in S}E_hV$.
From Lemma \ref{lem:modcharS99}(ii)
and the definition of a Leonard system,
we find
\begin{equation}
j \in S  \quad \hbox{and}\quad
|i-j|=1 \quad \rightarrow \quad i \in S,
\label{eq:implicosplitS99}
\end{equation}
for $0 \leq i,j\leq d$. Observe $S\not=\emptyset $ since $W\not=0$.
Combining this with 
(\ref{eq:implicosplitS99}), we find 
$S=\lbrace 0,1,\ldots, d\rbrace$, 
 so $W=V$.

\end{proof}

%

\begin{definition}
\label{def:VijS99}
With reference to 
Definition
\ref{def:splitconsetupS99}, 
we set 
\begin{equation}
V_{ij} =  \Biggl(\sum_{h=0}^i E^*_hV\Biggr)\cap
 \Biggl(\sum_{k=j}^d E_kV\Biggr)
\label{eq:defofvijS99}
\end{equation}
for all integers $i,j$. We interpret the sum on the left in
(\ref{eq:defofvijS99}) to be 0 (resp. $V$) if $i<0$  (resp. $i>d$).
Similarily, we interpret the sum on the right in
(\ref{eq:defofvijS99}) to be V (resp. $0$) if $j<0$  (resp. $j>d$).

\end{definition}

\begin{lemma}
\label{lem:thevijbasicfactsS99}
With reference to 
Definition
\ref{def:splitconsetupS99} and 
Definition \ref{def:VijS99}, we have
\begin{enumerate}
\item
$V_{i0} = 
 E^*_0V+E^*_1V+\cdots+E^*_iV \qquad \qquad  (0 \leq i \leq d),$
\item
$V_{dj} = 
 E_jV+E_{j+1}V+\cdots+E_dV \qquad \qquad (0 \leq j\leq d).$
\end{enumerate}

\end{lemma}

\begin{proof} 
To get (i), set $j=0$ in 
(\ref{eq:defofvijS99}), and apply 
(\ref{eq:VdecompS99}). Line  (ii) is similarily obtained.
\end{proof}

\begin{lemma}
\label{lem:howaactsonvijS99}
With reference to 
Definition
\ref{def:splitconsetupS99} 
and Definition \ref{def:VijS99}, the following (i)--(iv) hold
for $0 \leq i,j\leq d$.
\begin{enumerate}
\item
$(A-\theta_jI)V_{ij} \subseteq V_{i+1,j+1}$,
\item
$AV_{ij} \subseteq V_{ij} + V_{i+1,j+1}$, 
\item
$(A^*-\theta^*_iI)V_{ij} \subseteq V_{i-1,j-1}$,
\item
$A^*V_{ij} \subseteq V_{ij} + V_{i-1,j-1}$.
\end{enumerate}
\end{lemma}

\begin{proof} (i) Recall $E^*_rAE^*_h = 0 $ for  $r-h>1$,
so 
\begin{equation}
(A-\theta_jI)\sum_{h=0}^i E^*_hV \subseteq 
\sum_{h=0}^{i+1} E^*_hV.
\label{eq:aminusthetajAS99}
\end{equation}
Using
(\ref{eq:primid1S99}), we obtain
\begin{equation}
(A-\theta_jI)\sum_{h=j}^d E_hV \subseteq 
\sum_{h=j+1}^d E_hV.
\label{eq:aminusthetajBS99}
\end{equation}
Evaluating $(A-\theta_jI)V_{ij}$ using 
(\ref{eq:defofvijS99}),
(\ref{eq:aminusthetajAS99}), 
(\ref{eq:aminusthetajBS99}), 
 we
routinely find it is contained in $V_{i+1,j+1}$.

\noindent (ii) Immediate from  (i) above.

\noindent (iii),(iv) Apply (i), (ii) above 
to $\Phi^{\downarrow \Downarrow *}$ 

\end{proof}

\begin{lemma}
\label{lem:manyvijzeroS99}
With reference to 
Definition
\ref{def:splitconsetupS99} and  
Definition \ref{def:VijS99}, we have 
\begin{equation}
V_{ij}= 0 \quad \hbox{if}\quad i<j\qquad  \qquad (0 \leq i,j\leq d).
\label{eq:manyvijzeroS99}
\end{equation}
\end{lemma}

\begin{proof} We show the sum
\begin{equation}
V_{0r}+V_{1,r+1}+\cdots +V_{d-r,d}
\label{eq:manyvijzeroAS99}
\end{equation}
equals 0 for $0 <r \leq d$. Let $r$ be given, and let $W$ denote the sum
in 
(\ref{eq:manyvijzeroAS99}). Applying Lemma 
\ref{lem:howaactsonvijS99}(ii),(iv), we find
$AW\subseteq W$ and $A^*W\subseteq W$, so $W$ is
an $(A,A^*)$-module.
Applying Lemma  
\ref{lem:VirredaastarmoduleS99}, we find
$W=0$ or $W=V$.
 By Definition
\ref{def:VijS99}, each term in 
(\ref{eq:manyvijzeroAS99}) is contained in 
\begin{equation}
E_rV+E_{r+1}V+\cdots +E_dV,
\label{eq:manyvijzeroBS99}
\end{equation}
so $W$ is contained in 
(\ref{eq:manyvijzeroBS99}). Apparently $W\not=V$, so
$W=0$. We have now shown  
(\ref{eq:manyvijzeroAS99}) is zero for $0<r\leq d$, and 
(\ref{eq:manyvijzeroS99}) follows.

\end{proof}

\begin{lemma}
\label{lem:vijgivesdirsumS99}
With reference to 
Definition
\ref{def:splitconsetupS99} and  
Definition \ref{def:VijS99}, we have
\begin{equation}
\hbox{dim}\,V_{ii}=1 \qquad \qquad  (0\leq i\leq d)
\label{eq:vijgivesdirsumBS99}
\end{equation}
and 
\begin{equation}
V= V_{00}+V_{11}+\cdots+V_{dd} \qquad \qquad  (\hbox{direct sum}).
\label{eq:vijgivesdirsumAS99}
\end{equation}

\end{lemma}

\begin{proof} We  first show
(\ref{eq:vijgivesdirsumAS99}).
Let $W$ denote the sum on the right in 
(\ref{eq:vijgivesdirsumAS99}).
Observe
$AW\subseteq W$ by 
Lemma 
\ref{lem:howaactsonvijS99}(ii), and 
$A^*W\subseteq W$ by 
Lemma 
\ref{lem:howaactsonvijS99}(iv), so $W$ is
an $(A,A^*)$-module.
 Applying Lemma
\ref{lem:VirredaastarmoduleS99}, we find
$W=0$ or $W=V$.
Observe $W$ contains 
 $V_{00}$, and $V_{00}=E^*_0V$ is nonzero,
 so
$W\not=0$.
It follows $W=V$, and in other words 
\begin{equation}
V= V_{00}+V_{11}+\cdots+V_{dd}.
\label{eq:vijgivesdirsumCS99}
\end{equation}
We show the sum 
(\ref{eq:vijgivesdirsumCS99}) is direct.  To do this,
we show
\beast
 (V_{00}+V_{11}+\cdots +V_{i-1,i-1})\cap V_{ii} =0
\eeast
for $1 \leq i \leq d$.
Let the integer $i$ be given. By
(\ref{eq:defofvijS99}) we find 
\beast
V_{jj} \subseteq E^*_0V+E^*_1V+\cdots +E^*_{i-1}V
\eeast
for $0 \leq j \leq i-1$, and
\beast
V_{ii}\subseteq E_iV+E_{i+1}V+\cdots +E_dV.
\eeast
It follows
\beast
 &&(V_{00}+V_{11}+\cdots +V_{i-1,i-1})\cap V_{ii}
\\
 && \qquad  \subseteq
 (E^*_0V+E^*_1V+\cdots +E^*_{i-1}V)\cap
(E_iV+E_{i+1}V+\cdots +E_dV) \qquad \qquad 
\\
&& \qquad  = V_{i-1,i}
\\
&& \qquad  = 0
\eeast
in view of Lemma 
\ref{lem:manyvijzeroS99}.
We have now shown
the sum (\ref{eq:vijgivesdirsumCS99}) is direct,
so 
(\ref{eq:vijgivesdirsumAS99}) holds.
It remains to show
(\ref{eq:vijgivesdirsumBS99}). In view of
(\ref{eq:vijgivesdirsumAS99}), it suffices to show
$V_{ii}\not=0$ for $0 \leq i \leq d$. Suppose there
exists an integer $i$ $(0 \leq i \leq d)$ such that
$V_{ii}=0$. We observe $i\not=0$, since $V_{00}=E^*_0V$
is nonzero, and $i\not=d$, since $V_{dd}=E_dV$ is nonzero.
Set 
\beast
U = V_{00}+V_{11}+\cdots+V_{i-1,i-1},
\eeast
and observe
 $U\not=0$ and $U\not=V$ by our remarks above.
By 
Lemma 
\ref{lem:howaactsonvijS99}(ii) and since $V_{ii}=0$,
we find 
$AU\subseteq U$.
By 
Lemma 
\ref{lem:howaactsonvijS99}(iv) we find 
$A^*U\subseteq U$. Now $U$ is an $(A,A^*)$-module.
Applying Lemma \ref{lem:VirredaastarmoduleS99},
 we find $U=0$ or $U=V$,
contradicting our comments above.
We conclude 
 $V_{ii}\not=0$ for $0 \leq i \leq d$, and
(\ref{eq:vijgivesdirsumBS99})
follows.

\end{proof}

\begin{lemma}
\label{lem:animprovementonhowactS99}
With reference to 
Definition
\ref{def:splitconsetupS99}  and 
Definition \ref{def:VijS99}, the following 
 (i)--(iv)   hold.
\begin{enumerate}
\item
$(A-\theta_iI)V_{ii} = V_{i+1,i+1} \qquad \qquad (0\leq i \leq d-1)$,
\item
$(A-\theta_dI)V_{dd} = 0$,
\item
$(A^*-\theta^*_iI)V_{ii} = V_{i-1,i-1} \qquad \qquad (1\leq i \leq d)$,
\item
$(A^*-\theta^*_0I)V_{00} = 0$.
\end{enumerate}

\end{lemma} 

\begin{proof} (i) Let the integer $i$ be given.
Recall  
$(A-\theta_iI)V_{ii}$ is contained in $V_{i+1,i+1}$
by Lemma
\ref{lem:howaactsonvijS99}(i), and 
$V_{i+1,i+1}$ has dimension 1 by
(\ref{eq:vijgivesdirsumBS99}), so it suffices to show
\begin{equation}
(A-\theta_iI)V_{ii}\not=0.
\label{eq:aviinotzeroS99}
\end{equation}
Assume 
$(A-\theta_iI)V_{ii}=0$, and 
set
\beast
W = V_{00}+V_{11}+\cdots +V_{ii}.
\eeast
By 
Lemma \ref{lem:vijgivesdirsumS99}, and since
$0 \leq i<d$, we find 
$W\not=0$ and $W\not=V$.
Observe 
$AV_{ii}  \subseteq V_{ii}$ by our above assumption;
combining this with 
 Lemma   
\ref{lem:howaactsonvijS99}(ii),
we find $AW \subseteq W$.
By 
Lemma \ref{lem:howaactsonvijS99}(iv), we find
$A^*W\subseteq W$. Now $W$ is an $(A, A^*)$-module.
Applying Lemma 
\ref{lem:VirredaastarmoduleS99}, we find 
$W=0$ or $W=V$, contradicting our above remarks.
We conclude 
(\ref{eq:aviinotzeroS99}) holds, and the result follows. 

\noindent (ii) Recall $V_{dd} = E_dV$  by Lemma 
\ref{lem:thevijbasicfactsS99}(ii).

\noindent (iii), (iv)
Apply (i),(ii) above to $\Phi^{\downarrow \Downarrow *}$.

\end{proof}
\noindent{\it Proof of Theorem
\ref{thm:splitformexistS99}:
}
We first show 
 there exists  
an isomorphism of $\fld$-algebras 
$\flat :\alg \rightarrow  \Mdf$
and there exists 
nonzero scalars $\varphi_1, \varphi_2,\ldots, \varphi_d$ in $\fld$
such that
(\ref{eq:ahrtandastarhrtS99}) holds.
Let $V$ denote the irreducible left $\alg$-module, and
let the 
subspaces $V_{ij}$ of $V$ be as in 
Definition \ref{def:VijS99}. Let $v_0$ denote a nonzero
vector in $V_{00}=E^*_0V$, and let $v_1, v_2, \ldots, v_d$
denote the vectors in $V$ satisfying 
\begin{equation}
(A-\theta_iI)v_i = v_{i+1} \qquad \qquad (0 \leq i \leq d-1).
\label{eq:recdefofvss99}
\end{equation}
Combining  Lemma 
\ref{lem:animprovementonhowactS99}(i) and 
(\ref{eq:vijgivesdirsumBS99}), we find 
$v_i$ is a basis for $V_{ii}$, for $0 \leq i \leq  d$.
By this  and 
(\ref{eq:vijgivesdirsumAS99}), 
we find 
\begin{equation}
v_0, v_1, \ldots, v_d
\label{eq:vibasisforVS99}
\end{equation}
is a basis for $V$.
For all $X \in \alg$, let $X^\flat$ denote the matrix
in $\Mdf$ that represents the $X$ action on $V$ with respect
to the basis 
(\ref{eq:vibasisforVS99}). By elementary linear algebra, we 
find the map $\flat: X \rightarrow  X^\flat $
is an isomorphism of $\fld$-algebras from $\alg$ to $\Mdf$.
We check $A^\flat$ and  $A^{*\flat}$ have
the required  form
(\ref{eq:ahrtandastarhrtS99}). 
To get $A^\flat$,
we consider the action of $A$ on the basis
(\ref{eq:vibasisforVS99}).
Most of this action is given in 
(\ref{eq:recdefofvss99}), but we still need $Av_d$.
Recall $v_d$ is contained in  
$V_{dd}$, 
so by Lemma
\ref{lem:animprovementonhowactS99}(ii),
\begin{equation}
Av_d = \theta_dv_d.
\label{eq:AlasttermS99}
\end{equation}
Combining   
(\ref{eq:recdefofvss99}), 
(\ref{eq:AlasttermS99}), 
we find 
the equation on the left in 
(\ref{eq:ahrtandastarhrtS99})
holds.
%
To get $A^{*\flat}$, 
we consider the action of 
$A^*$ on the basis  
(\ref{eq:vibasisforVS99}).
By Lemma 
\ref{lem:animprovementonhowactS99}(iii), there
exist nonzero scalars $\varphi_1, \varphi_2,\ldots, \varphi_d$
in $\fld $ such that 
\begin{equation}
(A^*-\theta^*_iI)v_i = \varphi_iv_{i-1} \qquad \qquad (1\leq i \leq d).
\label{eq:wegodownS99}
\end{equation}
By Lemma 
\ref{lem:animprovementonhowactS99}(iv), 
\begin{equation}
(A^*-\theta^*_0I)v_0 = 0.
\label{eq:wegodownlastS99}
\end{equation}
Combining 
(\ref{eq:wegodownS99}), 
(\ref{eq:wegodownlastS99}),
we obtain 
the equation on the right in 
(\ref{eq:ahrtandastarhrtS99}).
%

\noindent 
Concerning uniqueness, let $\flat':\alg \rightarrow \Mdf$ denote an
isomorphism of $\fld$-algebras,
and let $\varphi'_1, \varphi'_2,\ldots, \varphi'_d$ denote 
scalars in $\fld $ such that 
(\ref{eq:ahrtandastarhrtS99}) hold.
We show
$\flat = \flat'$. To do this, we show the composition
$\sigma = \flat^{-1}\flat' $ equals 1.
Observe $\sigma $ is an isomorphism of $\fld$-algebras from
$\Mdf$ to itself, so
by elementary linear algebra, 
there exists an invertible matrix $X$ in $\Mdf$ such that
\begin{equation}
Y^\sigma = X^{-1}YX \qquad \qquad (\forall Y \in \Mdf).
\label{eq:splitcanonformuniqueAS99}
\end{equation}
Apparently 
\begin{equation}
XY^\sigma = YX \qquad \qquad (\forall Y \in \Mdf).
\label{eq:splitcanonformuniqueBS99}
\end{equation}
If $Y=A^{*\flat}$  then 
 $Y^\sigma =A^{*\flat'}$; using these values in
(\ref{eq:splitcanonformuniqueBS99}), we get 
\begin{equation}
XA^{*\flat'}= 
 A^{*\flat}X.
\label{eq:splitcanonformuniqueCS99}
\end{equation}
Both $A^{*\flat}$ and 
 $A^{*\flat'}$ are
 upper triangular,  and
each
has $ii$ entry $\theta^*_i$ for $0 \leq i \leq d$;
multiplying out each side of 
(\ref{eq:splitcanonformuniqueCS99}) using this, 
and recalling $\theta^*_0, \theta^*_1, \ldots, \theta^*_d$ are
distinct, 
we find $X$ is upper triangular.
If $Y=A^\flat$
then 
 $Y^\sigma =A^{\flat'}$; using these values in
(\ref{eq:splitcanonformuniqueBS99}), 
we get
\begin{equation}
XA^{\flat'}= 
 A^{\flat}X.
\label{eq:splitcanonformuniqueDS99}
\end{equation}
The matrices 
$A^\flat$  and  
 $A^{\flat'}$ are equal, and are given by the equation
 on the left in 
(\ref{eq:ahrtandastarhrtS99});
multiplying out each side of
(\ref{eq:splitcanonformuniqueDS99}) using  this,
and recalling 
 $\theta_0, \theta_1, \ldots, \theta_d$ are distinct, 
we find $X$ is 
a  scalar multiple of the identity.
It follows $\sigma=1$, so $\flat=\flat'$.
The scalars $\varphi_i'$ are determined by 
$\Phi$ and $\flat'$, so 
 $\varphi'_i=\varphi_i$ for $1 \leq i \leq d$.
This completes the proof.
\par \noindent $\Box $ \par 

\begin{definition}
\label{def:varphidefS99}
With reference to Definition
\ref{def:splitconsetupS99}, 
by the split canonical form of $\Phi$, we mean
the Leonard system $\Phi^\flat$, where
$\flat=\flat(\Phi)$  
is from 
Theorem
\ref{thm:splitformexistS99}.
By the $\varphi$-sequence of $\Phi$,  we mean the
seqence $\varphi_1, \varphi_2,\ldots, \varphi_d$
from Theorem
\ref{thm:splitformexistS99}.
For notational convenience, we define
$\varphi_0=0$, $\varphi_{d+1}=0$.
\end{definition}

\begin{lemma}
\label{lem:paramsdetisoS99}
Let $\Phi $ and $\Phi'$  denote a Leonard systems over 
$\fld$.
Then the following are equivalent.
\begin{enumerate}
\item $\Phi$ and $\Phi'$ are isomorphic.
\item $\Phi$ and $\Phi'$ share  the same eigenvalue sequence,
dual eigenvalue sequence, and $\varphi$-sequence. 
\end{enumerate}
\end{lemma}

\begin{proof} $(i)\rightarrow (ii)$ Routine.

\noindent $(ii)\rightarrow (i)$ Observe $\Phi$ and $\Phi'$
are both isomorphic to a common split canonical form, so they are  
isomorphic.

\end{proof}

\noindent  
It is helpful to consider the following parameters.

\begin{definition}
\label{def:phidefS99}
With reference to Definition
\ref{def:splitconsetupS99}, 
let $\phi_1, \phi_2, \ldots, \phi_d$ denote the
$\varphi$-sequence of $\Phi^\Downarrow $. 
Let the map $\flat = \flat(\Phi)$ be as in 
Theorem
\ref{thm:splitformexistS99}, and let 
 $\diamondsuit= \flat^\Downarrow$ denote the corresponding map
 for $\Phi^\Downarrow $.
Observe
\beast
A^\diamondsuit = 
\left(
\begin{array}{c c c c c c}
\theta_d & & & & & {\bf 0} \\
1 & \theta_{d-1} &  & & & \\
& 1 & \theta_{d-2} &  & & \\
& & \cdot & \cdot &  &  \\
& & & \cdot & \cdot &  \\
{\bf 0}& & & & 1 & \theta_0
\end{array}
\right),
&&\quad 
A^{*\diamondsuit} = 
\left(
\begin{array}{c c c c c c}
\theta^*_0 &\phi_1 & & & & {\bf 0} \\
 & \theta^*_1 & \phi_2 & & & \\
&  & \theta^*_2 & \cdot & & \\
& &  & \cdot & \cdot &  \\
& & &  & \cdot & \phi_d \\
{\bf 0}& & & &  & \theta^*_d
\end{array}
\right).
\eeast
We call 
 $\phi_1, \phi_2, \ldots, \phi_d$  the
$\phi$-sequence of $\Phi $.  For notational convenience,
we define $\phi_0=0$, $\phi_{d+1}=0$. 
\end{definition}

\begin{lemma}
\label{lem:newparamsdetisoS99}
Let $\Phi $ and $\Phi'$  denote a Leonard systems over 
$\fld$.
Then the following are equivalent.
\begin{enumerate}
\item $\Phi$ and $\Phi'$ are isomorphic.
\item $\Phi$ and $\Phi'$ share  the same eigenvalue sequence,
dual eigenvalue sequence, and $\phi$-sequence. 
\end{enumerate}
\end{lemma}

\begin{proof}  Apply Lemma
\ref{lem:paramsdetisoS99} to 
the second inversions of $\Phi$ and $\Phi'$.

\end{proof}

\section{The primitive idempotents of a Leonard system}

\noindent 
In this section,  we 
consider a Leonard system in 
 split canonical form, and compute the entries
 of the primitive idempotents.
We begin by considering a more general situation.

\begin{definition}
\label{def:setupforsectionS99} Let 
$d$ denote a nonnegative integer and let $\fld $ denote a field.
In this section, we let 
$A$ and $A^*$ denote any matrices in $\Mdf$ of the form 
\beast
A = 
\left(
\begin{array}{c c c c c c}
\theta_0 & & & & & {\bf 0} \\
1 & \theta_1 &  & & & \\
& 1 & \theta_2 &  & & \\
& & \cdot & \cdot &  &  \\
& & & \cdot & \cdot &  \\
{\bf 0}& & & & 1 & \theta_d
\end{array}
\right),
&&\quad 
A^* = 
\left(
\begin{array}{c c c c c c}
\theta^*_0 &\varphi_1 & & & & {\bf 0} \\
 & \theta^*_1 & \varphi_2 & & & \\
&  & \theta^*_2 & \cdot & & \\
& &  & \cdot & \cdot &  \\
& & &  & \cdot & \varphi_d \\
{\bf 0}& & & &  & \theta^*_d
\end{array}
\right),
\eeast
such that
\begin{eqnarray}
&&\theta_i \not= \theta_j, \qquad   
\theta^*_i \not= \theta^*_j \qquad \hbox{if} \quad i\not=j \qquad  
(0 \leq i,j\leq d), \qquad \qquad 
\\
&& 
 \qquad  \varphi_i \not=0 \qquad (1 \leq i \leq d).
\end{eqnarray}
We observe $A$ (resp. $A^*$) is multiplicity-free,
with eigenvalues
$\theta_0,\theta_1,\ldots, \theta_d$ (resp.  
$\theta^*_0,\theta^*_1,\ldots, \theta^*_d$).
For $0 \leq i \leq d$, 
we let $E_i$ (resp. $E^*_i$) denote the primitive idempotent for 
$A$ (resp. $A^*$) associated with $\theta_i$ (resp. $\theta^*_i$). 
For notational convenience, we  define $\varphi_0=0$, 
$\varphi_{d+1}=0$.
\end{definition}

\noindent    
With reference to 
Definition \ref{def:setupforsectionS99},
we do not assume $A, A^*$ come from a Leonard system,
so we do not expect any $D_4$ action; however we do
have the following result.

\begin{lemma}
\label{lem:insteadofd4actionS99}
With reference to Definition
\ref{def:setupforsectionS99},  let $G$ denote the 
diagonal matrix in 
$\Mdf$ with  diagonal entries
\begin{equation}
G_{ii} = 
\varphi_1\varphi_2\cdots \varphi_i \qquad \qquad (0 \leq i \leq d).
\label{eq:gremindfreeS99}
\end{equation}
Let $Z$ denote the matrix in 
$\Mdf$ with $ij^{\hbox{th}}$ entry 1 if $i+j=d$, and 0 if
 $i+j\not=d$, for $0 \leq i,j\leq d$. 
Then (i)--(iii) hold below.
\begin{enumerate}
\item
The matrices
$G^{-1}A^{*t}G$ and   
$G^{-1}A^tG$ equal 
\beast
\left(
\begin{array}{c c c c c c}
\theta^*_0 & & & & & {\bf 0} \\
1 & \theta^*_1 &  & & & \\
& 1 & \theta^*_2 &  & & \\
& & \cdot & \cdot &  &  \\
& & & \cdot & \cdot &  \\
{\bf 0}& & & & 1 & \theta^*_d
\end{array}
\right),
\qquad 
\left(
\begin{array}{c c c c c c}
\theta_0 &\varphi_1 & & & & {\bf 0} \\
 & \theta_1 & \varphi_2 & & & \\
&  & \theta_2 & \cdot & & \\
& &  & \cdot & \cdot &  \\
& & &  & \cdot & \varphi_d \\
{\bf 0}& & & &  & \theta_d
\end{array}
\right),
\eeast
respectively.
\item The matrices  
$ZA^{t}Z$ and  
$ZA^{*t}Z$ equal
\beast
\left(
\begin{array}{c c c c c c}
\theta_d & & & & & {\bf 0} \\
1 & \theta_{d-1} &  & & & \\
& 1 & \theta_{d-2} &  & & \\
& & \cdot & \cdot &  &  \\
& & & \cdot & \cdot &  \\
{\bf 0}& & & & 1 & \theta_0
\end{array}
\right),
\qquad 
\left(
\begin{array}{c c c c c c}
\theta^*_d &\varphi_d & & & & {\bf 0} \\
 & \theta^*_{d-1} & \varphi_{d-1} & & & \\
&  & \theta^*_{d-2} & \cdot & & \\
& &  & \cdot & \cdot &  \\
& & &  & \cdot & \varphi_1 \\
{\bf 0}& & & &  & \theta^*_0
\end{array}
\right),
\eeast
respectively. 
\item 
 The matrices 
$ZGA^{*}G^{-1}Z$ and
$ZGAG^{-1}Z$ equal 
\beast
\left(
\begin{array}{c c c c c c}
\theta^*_d & & & & & {\bf 0} \\
1 & \theta^*_{d-1} &  & & & \\
& 1 & \theta^*_{d-2} &  & & \\
& & \cdot & \cdot &  &  \\
& & & \cdot & \cdot &  \\
{\bf 0}& & & & 1 & \theta^*_0
\end{array}
\right),
\qquad 
\left(
\begin{array}{c c c c c c}
\theta_d &\varphi_d & & & & {\bf 0} \\
 & \theta_{d-1} & \varphi_{d-1} & & & \\
&  & \theta_{d-2} & \cdot & & \\
& &  & \cdot & \cdot &  \\
& & &  & \cdot & \varphi_1 \\
{\bf 0}& & & &  & \theta_0
\end{array}
\right),
\eeast
respectively.
\end{enumerate}
We remark $Z=Z^{-1}$.
\end{lemma}

\begin{proof} Routine matrix multiplication.

\end{proof}

\noindent We remark that in the above lemma only (i) and (iii) will
be used later in the paper; we include (ii) for the sake of completeness.

\noindent
It is convenient to use the following notation.

\begin{definition}
\label{def:thetauS99} Suppose we are given a 
nonnegative integer $d$ and a sequence
of scalars 
\beast
\theta_0,\theta_1,\ldots, \theta_d;\qquad \qquad   
\theta^*_0,\theta^*_1,\ldots, \theta^*_d
\eeast
taken from a field $\fld$. 
Then for  $0 \leq i \leq d+1$, 
we let   
$\tau_i$, 
$\tau^*_i$,
$\eta_i$,
$\eta^*_i$ denote the following polynomials in 
$\fld\lbrack \lambda \rbrack $.  
\begin{eqnarray}
\tau_i &:=& 
(\lambda-\theta_0)(\lambda - \theta_1)\cdots (\lambda-\theta_{i-1}),  
\label{eq:taudeffreeS99}
\\
\tau^*_i &:=& 
(\lambda-\theta^*_0)(\lambda - \theta^*_1)\cdots (\lambda-\theta^*_{i-1}),  
\label{eq:taudeffreedualS99}
\\
\eta_i &:=& 
(\lambda-\theta_d)(\lambda - \theta_{d-1})\cdots (\lambda-\theta_{d-i+1}),  
\label{eq:etafreeS99}
\\
\eta^*_i &:=& 
(\lambda-\theta^*_d)(\lambda - \theta^*_{d-1})\cdots 
(\lambda-\theta^*_{d-i+1}).  
\label{eq:etafreedualS99}
\end{eqnarray}
We observe each of $\tau_i$, 
 $\tau^*_i$, 
$\eta_i$,
 $\eta^*_i$
is  monic of degree $i$.  
\end{definition}

\begin{lemma}
\label{lem:theEisS99}  
With reference to 
Definition \ref{def:setupforsectionS99},
Pick any integer $r$ $(0 \leq r\leq d)$.
Then  the primitive idempotent 
$E_r$  of $A$
  has 
$ij^{\hbox{th}}$ entry 
\begin{equation}
{{\tau_j(\theta_r)\eta_{d-i}(\theta_r)}\over
{\tau_r(\theta_r)\eta_{d-r}(\theta_r)}}
\label{eq:entriesofefreeS99}
\end{equation}
for $0 \leq i,j\leq d$. 
We are using the notation
(\ref{eq:taudeffreeS99}), 
(\ref{eq:etafreeS99}).

\end{lemma}

\begin{proof} Let the integers $i,j$ be given.
Computing the $ij$ entry of  $AE_r=\theta_rE_r$ using
matrix multiplication, and taking into account the form
of $A$ in Definition
 \ref{def:setupforsectionS99},
we find
\beast
(E_r)_{i-1,j}
= (\theta_r-\theta_i)(E_r)_{ij} 
\eeast
if $i\geq 1$.
Replacing $i$ by $i+1$ in the above line, we find
\begin{equation}
(E_r)_{ij}
= (\theta_r-\theta_{i+1})(E_r)_{i+1,j} 
\label{eq:firstrecS99}
\end{equation}
if $i\leq d-1$.
Using the recursion 
(\ref{eq:firstrecS99}), 
 we routinely find
\begin{eqnarray}
(E_r)_{ij} &=& (\theta_r-\theta_{i+1})
(\theta_r-\theta_{i+2}) \cdots
(\theta_r-\theta_d)(E_r)_{dj}
\nonumber
\\
&=& \eta_{d-i}(\theta_r)(E_r)_{dj}.
\label{eq:firstpartrecS99}
\end{eqnarray}
Computing the $dj$ entry of 
$E_rA=\theta_rE_r$ using
matrix multiplication, and taking into account the form
of $A$, we find
\beast
(E_r)_{d,j+1} 
= (\theta_r-\theta_j)(E_r)_{dj} 
\eeast
if $j\leq d-1$.  Replacing $j$ by $j-1$ in the above line we find 
\begin{equation}
(E_r)_{dj} 
= (\theta_r-\theta_{j-1})(E_r)_{d,j-1} 
\label{eq:secondrecS99}
\end{equation}
if $j\geq 1$.
Using the recursion 
(\ref{eq:secondrecS99}),
 we routinely find
\begin{eqnarray}
(E_r)_{dj} &=& (\theta_r-\theta_{j-1})
(\theta_r-\theta_{j-2}) \cdots
(\theta_r-\theta_0)(E_r)_{d0}
\nonumber
\\
&=& \tau_{j}(\theta_r)(E_r)_{d0}.
\label{eq:secondpartrecS99}
\end{eqnarray}
Combining 
(\ref{eq:firstpartrecS99}),
(\ref{eq:secondpartrecS99}), we find
\begin{equation}
(E_r)_{ij} =
\tau_j(\theta_r)\eta_{d-i}(\theta_r)c,
\label{eq:togetherS99}
\end{equation}
where we abbreviate  $c= (E_r)_{d0}$. We now  find $c$.
Since $A$ is lower triangular, and since 
$E_r$ is a polynomial in 
$A$, we see 
$E_r$ is lower triangular.
Recall $E_r^2=E_r$, so the diagonal entry of $(E_r)_{rr}$
equals 0 or 1. We show 
$(E_r)_{rr}=1$.
Setting $i=r$, $j=r$ in (\ref{eq:togetherS99}),
\begin{equation}
(E_r)_{rr}=
\tau_r(\theta_r)\eta_{d-r}(\theta_r)c.
\label{eq:rrentryoferS99}
\end{equation}
Observe 
$\tau_r(\theta_r)\not=0$ and $\eta_{d-r}(\theta_r)\not=0$ by
Definition
\ref{def:thetauS99},  
and since  the eigenvalues 
are distinct.
Observe $c\not=0$; otherwise $E_r=0$ in view of 
(\ref{eq:togetherS99}). Apparently the right side of
(\ref{eq:rrentryoferS99}) is not 0, so 
$(E_r)_{rr}\not=0$, and we conclude  
$(E_r)_{rr}=1$. Setting
$(E_r)_{rr}=1$ in 
(\ref{eq:rrentryoferS99}),  solving for $c$,
 and evaluating  
(\ref{eq:togetherS99}) using the result, we find
the 
$ij^{\hbox{th}}$ entry 
of $E_r$ is given by 
(\ref{eq:entriesofefreeS99}).

\end{proof}

\begin{example}
\label{ex:primidsdeq2freeS99}
With reference to  
Definition \ref{def:setupforsectionS99},
 for $d=2$, the primitive idempotents of
$A$ are as follows.  
\beast
E_0 = 
\left(
\begin{array}{c c c }
1 & 0  &  0
\\
{{1}\over {\theta_0-\theta_1}} & 0 & 0
\\
{{1}\over
{(\theta_0 - \theta_1)
(\theta_0 - \theta_2)}}
&
0
&
0
\end{array}
\right),
\qquad \quad E_1 = 
\left(
\begin{array}{c c c }
0 & 0  &  0
\\
{{1}\over {\theta_1-\theta_0}} & 1 & 0 \\
{{1}\over
{(\theta_1 - \theta_0)
(\theta_1 - \theta_2)}}
&
{{1}\over {\theta_1-\theta_2}}
&
0
\end{array}
\right),
\eeast
\beast
E_2 = 
\left(
\begin{array}{c c c }
0 & 0  &  0 \\
0 & 0 & 0 \\
{{1}\over
{(\theta_2 - \theta_1)
(\theta_2 - \theta_0)}}
&
{{1}\over {\theta_2-\theta_1}}
&
1
\end{array}
\right).
\eeast
\end{example}

\begin{lemma}
\label{lem:theEistarsS99}  
With reference to 
Definition \ref{def:setupforsectionS99},
Pick any integer $r$ $(0 \leq r \leq d)$. Then  
the primitive idempotent $E^*_r$  of $A^*$ 
  has 
$ij^{\hbox{th}}$ entry 
\begin{equation}
{{\varphi_1 \varphi_2\cdots \varphi_j}
\over{\varphi_1\varphi_2\cdots \varphi_i}}
\,{{\tau^*_i(\theta^*_r)\eta^*_{d-j}(\theta^*_r)}\over{
\tau^*_r(\theta^*_r)\eta^*_{d-r}(\theta^*_r)}}
\label{eq:theestarentriesS99}
\end{equation}
for $0 \leq i,j\leq d$. We are using the notation
(\ref{eq:taudeffreedualS99}),
(\ref{eq:etafreedualS99}).
\end{lemma}

\begin{proof} Let the matrix $G$ be as in  Lemma
\ref{lem:insteadofd4actionS99},
and  set
$A':=G^{-1}A^{*t}G$; this matrix is
 given on  
the left in 
Lemma 
\ref{lem:insteadofd4actionS99}(i).
Let $E'_r$ denote the primitive idempotent of 
$A'$ associated with the eigenvalue $\theta^*_r$.
We find $E'_r$ in two ways. 
On one hand, applying 
 Lemma  
\ref{lem:theEisS99} to $A'$, we find  
$E'_r$ 
  has 
$ij^{\hbox{th}}$ entry 
\begin{equation}
{{\tau^*_j(\theta^*_r)\eta^*_{d-i}(\theta^*_r)}\over
{\tau^*_r(\theta^*_r)\eta^*_{d-r}(\theta^*_r)}}
\label{eq:entriesofestarAfreeS99}
\end{equation}
for $0 \leq i,j\leq d$.  
On the other hand, by elementary linear algebra
\beast
E'_r = G^{-1}E^{*t}_rG,
\eeast
so  
$E'_r$ 
  has 
$ij^{\hbox{th}}$ entry 
\begin{eqnarray}
G_{ii}^{-1} (E^*_r)_{ji}G_{jj}
=  
{{\varphi_1 \varphi_2\cdots \varphi_j}
\over{\varphi_1\varphi_2\cdots \varphi_i}} 
\,(E^*_r)_{ji}
\label{eq:theestarentriesyzS99}
\end{eqnarray}
for $0 \leq i,j\leq d$.
Equating (\ref{eq:entriesofestarAfreeS99}) and the
right side of 
(\ref{eq:theestarentriesyzS99}),
and solving for
  $(E^*_r)_{ji}$,
we routinely obtain the result.

\end{proof}

\begin{example}
\label{ex:primidsdeq2dualfreeS99}
With reference to  
Definition \ref{def:setupforsectionS99},
 for $d=2$, the primitive idempotents of $A^*$
 are  given by 
\beast
E^*_0 = 
\left(
\begin{array}{c c c }
1 
&
{{\varphi_1}\over{\theta^*_0-\theta^*_1}}
& 
{{\varphi_1\varphi_2}\over{(\theta^*_0-\theta^*_1)(\theta^*_0-\theta^*_2)}}
\\
0 & 0 & 0
\\
0
&
0
&
0
\end{array}
\right),
\qquad \quad 
E^*_1 = 
\left(
\begin{array}{c c c }
0
& 
{{\varphi_1}\over{\theta^*_1-\theta^*_0}}
&
{{\varphi_1\varphi_2}\over{(\theta^*_1-\theta^*_0)(\theta^*_1-\theta^*_2)}}
\\
0
&
1
&
{{\varphi_2}\over{\theta^*_1-\theta^*_2}}
\\
0
&
0
&
0
\end{array}
\right),
\eeast
\beast
E^*_2 = 
\left(
\begin{array}{c c c }
0 & 0  &
{{\varphi_1\varphi_2}\over{(\theta^*_2-\theta^*_1)(\theta^*_2-\theta^*_0)}}
 \\
0 & 0 & 
{{\varphi_2}\over{\theta^*_2-\theta^*_1}}
 \\
0
&
0
&
1
\end{array}
\right).
\eeast

\end{example}

\noindent  We now  
restate Lemma
\ref{lem:theEisS99}   and Lemma 
\ref{lem:theEistarsS99}  
in the context of of Leonard systems.

\begin{theorem}
\label{thm:formofprimidS99}
Let
$\ls $
denote
a Leonard system, with eigenvalue sequence 
$\theta_0,\theta_1,\ldots, \theta_d$,  
dual eigenvalue sequence $\theta^*_0,\theta^*_1,\ldots, \theta^*_d$,
and $\varphi$-sequence 
$\varphi_1,\varphi_2,\ldots \varphi_d$.  
Let the map  $\flat=\flat(\Phi)$ be as in 
Theorem
\ref{thm:splitformexistS99}.
Pick any integer $r$ $(0\leq r \leq d)$. Then
$E_r^\flat$ is the matrix
in $\Mdf$ with 
$ij^{\hbox{th}}$ entry 
\begin{equation}
{{\tau_j(\theta_r)\eta_{d-i}(\theta_r)}\over
{\tau_r(\theta_r)\eta_{d-r}(\theta_r)}}
\label{eq:entriesofefreefinS99}
\end{equation}
for $0 \leq i,j\leq d$.  
Moreover,
$E^{*\flat}_r$  is the matrix in 
$\Mdf$ with 
$ij^{\hbox{th}}$ entry 
\begin{equation}
{{\varphi_1 \varphi_2\cdots \varphi_j}
\over{\varphi_1\varphi_2\cdots \varphi_i}}
\,{{\tau^*_i(\theta^*_r)\eta^*_{d-j}(\theta^*_r)}\over{
\tau^*_r(\theta^*_r)\eta^*_{d-r}(\theta^*_r)}}
\label{eq:theestarentriesfinS99}
\end{equation}
for $0 \leq i,j\leq d$. We are using the notation
(\ref{eq:taudeffreeS99})--(\ref{eq:etafreedualS99}).
\end{theorem}

\begin{proof} By Theorem 
\ref{thm:splitformexistS99}, the matrices $A^\flat$ and
$A^{*\flat}$ are of the form given in Definition
\ref{def:setupforsectionS99}.   
To get $E^\flat_r$,
apply Lemma
\ref{lem:theEisS99} 
to $A^\flat$, and observe
$E^\flat_r$ is the primitive idempotent
$A^\flat$ associated with $\theta_r$.
To get $E^{*\flat}_r$,
apply Lemma
\ref{lem:theEistarsS99}  
to $A^{*\flat}$, and observe
$E^{*\flat}_r$ is the primitive idempotent
$A^{*\flat}$ associated with $\theta^*_r$.

\end{proof}
\medskip
\noindent We finish this section with a few observations.

\begin{lemma}
\label{lem:commentoneis1S99}
With reference to Definition 
\ref{def:setupforsectionS99},
\begin{enumerate}
\item 
${\displaystyle{
E_iA^*E_j = \left\{ \begin{array}{ll}
                   0  & \mbox{if $\;j-i > 1; $ } \\
				  \not= 0 & \mbox{if $\;j-i=1 $ }
				   \end{array}
				\right. \qquad \qquad  (0 \leq i,j\leq d),
}}$
\item 
${\displaystyle{
E^*_iAE^*_j = \left\{ \begin{array}{ll}
                   0  & \mbox{if $\;i-j > 1; $ } \\
				  \not= 0 & \mbox{if $\;i-j=1 $ }
				   \end{array}
				\right. \qquad \qquad  (0 \leq i,j\leq d).
}}$
\end{enumerate}

\end{lemma}
\begin{proof}
(i)  For $0 \leq r \leq d$,  consider the pattern of zero entries in 
$E_r$.
By Lemma
\ref{lem:theEisS99},    
we find $E_r$ has all entries 0  in rows
$0,1,\ldots, r-1$ and columns $r+1,r+2,\ldots, d$.
Moreover, the $rr$ entry of $E_r$ is 1. 
Multiplying out $E_iA^*E_j$ using this and
the shape of $A^*$,
we find that if 
$ j-i>1$, then
all entries are 0, and 
if  $j-i=1$, then
the $ij$ entry
 is
 $\varphi_j$ and hence nonzero.
The result follows.

\noindent (ii) Very similar to the proof of (i) above.

\end{proof}

\begin{lemma}
\label{lem:def41andlsS99}
With reference to 
Definition
\ref{def:setupforsectionS99},   the following are equivalent.
\begin{enumerate}
\item 
$(A,A^*)$ is a Leonard pair in
$\Mdf$. 
\item
$(A;E_0,E_1,\ldots, E_d;A^*;E^*_0,E^*_1,\ldots,E^*_d)$
is a Leonard system in $\Mdf$.
\end{enumerate}
\noindent Suppose (i), (ii) hold. Then
the Leonard system  in (ii) 
has eigenvalue sequence
$\theta_0,\theta_1,\ldots, \theta_d$,
 dual eigenvalue
sequence 
$\theta^*_0,\theta^*_1,\ldots, \theta^*_d$,
and 
$\varphi$-sequence 
 $\varphi_1, \varphi_2,\ldots, \varphi_d$.

\end{lemma}

\begin{proof} Let $\Phi$ denote the sequence in part (ii).

\noindent 
$(i)\rightarrow (ii)$
We verify $\Phi$ 
satisfies the conditions (i)--(v)
of Definition \ref{def:deflstalkS99}, where we take
${\cal A} = \Mdf$. 
Conditions (i)--(iii) are immediate from the construction, so 
consider conditions (iv), (v). We assume $(A,A^*)$ is a Leonard
pair, so it is associated with some Leonard system.
It follows
\begin{eqnarray}
E_iA^*E_j=0 \quad &\hbox{iff}& \quad  E_jA^*E_i=0 \qquad \qquad 
(0 \leq i,j\leq d),
\label{eq:symmetryAS99}
\\
E^*_iAE^*_j=0 \quad &\hbox{iff}& \quad  E^*_jAE^*_i=0 \qquad \qquad 
(0 \leq i,j\leq d).
\label{eq:symmetryBS99}
\end{eqnarray}
Combining 
(\ref{eq:symmetryAS99}),
(\ref{eq:symmetryBS99})
with Lemma
\ref{lem:commentoneis1S99}, 
we find 
$\Phi$ satisfies the conditions (iv), (v)
of Definition \ref{def:deflstalkS99}.
We have now shown
$\Phi$ is a Leonard system in $\Mdf$.

\noindent 
$(ii)\rightarrow (i)$ Immediate from Definition
\ref{def:leonardpairS99}.

\noindent Now suppose (i), (ii) hold. 
From the construction 
$\Phi$ has eigenvalue sequence
$\theta_0,\theta_1,\ldots, \theta_d$ and 
 dual eigenvalue
sequence 
$\theta^*_0,\theta^*_1,\ldots, \theta^*_d$.
From the form of the matrices $A$ and $A^*$ in
Definition
\ref{def:setupforsectionS99},  
we find the map $\flat=\flat(\Phi)$ 
from Theorem
\ref{thm:splitformexistS99}
is the identity. Now $\Phi$ has $\varphi$-sequence
 $\varphi_1, \varphi_2,\ldots, \varphi_d$ in view
 of Definition
\ref{def:varphidefS99}. 
\end{proof}

\section{A formula for the $\varphi_i$}


\noindent In this section, we continue to consider the
situation of
Definition
\ref{def:setupforsectionS99}.
We obtain
the scalars $\varphi_i$ from that definition
in terms of 
the scalars $a_i$,  $a^*_i$ introduced in  
Definition
\ref{def:aidefS99}. To do this, we first obtain the $a_i$ and 
$a^*_i$
in terms of the $\varphi_i$.

\begin{lemma}
\label{thm:aiintermsofphivarphiS99}
With reference to Definition \ref{def:setupforsectionS99}
and Definition
\ref{def:aidefS99}, 
\begin{eqnarray}
a_i = \theta_i + {{\varphi_i}\over {\theta^*_i-\theta^*_{i-1}}}
+ {{\varphi_{i+1}}\over {\theta^*_i-\theta^*_{i+1}}}
\qquad \qquad (0 \leq i \leq d), 
\label{eq:aiintermsofvarphiS99}
\end{eqnarray}
where we recall
$\varphi_0=0$,
 $\varphi_{d+1}=0$,
and where  $\theta^*_{-1}$, $\theta^*_{d+1}$ denote 
 indeterminants.
Moreover, 
\begin{equation}
a^*_i= \theta^*_i + {{\varphi_i}\over {\theta_i-\theta_{i-1}}}
+ {{\varphi_{i+1}}\over {\theta_i-\theta_{i+1}}}
\qquad \qquad (0 \leq i \leq d),
\label{eq:aisintermsofvarphiS99} 
\end{equation}
where
 $\theta_{-1}$, $\theta_{d+1}$ denote 
 indeterminants.
\end{lemma}

\begin{proof} Concerning 
(\ref{eq:aiintermsofvarphiS99}),
recall $a_i$ equals the trace of $AE^*_i$, 
and this is the sum of the diagonal entries in
$AE^*_i$. Computing these entries by matrix multiplication,
and taking into account the form of $A$ in
Definition
\ref{def:setupforsectionS99}, we find that 
for $0 \leq j \leq d$, the $jj$ entry
\begin{eqnarray}
(AE^*_i)_{jj} &=& 
  \theta_j(E^*_i)_{jj}
 +(E^*_i)_{j-1,j}
\label{eq:multoutaestar}
\end{eqnarray}
where we interpret the term on the right in
(\ref{eq:multoutaestar}) to be zero if $j=0$.
By Lemma 
\ref{lem:theEistarsS99},
the diagonal entry $(E^*_i)_{jj}$ equals 1 if $j=i$, and
0 if $j\not=i$. Moreover, the entry
$(E^*_i)_{j-1,j}$ equals $\varphi_i(\theta^*_i-\theta^*_{i-1})^{-1}$
if $j=i$,    
$\varphi_{i+1}(\theta^*_i-\theta^*_{i+1})^{-1}$ if $j=i+1$, 
and  0 if $j\notin \lbrace i,i+1 \rbrace $. 
Evaluating 
(\ref{eq:multoutaestar}) using the above information, we readily
obtain (\ref{eq:aiintermsofvarphiS99}).
The proof of 
(\ref{eq:aisintermsofvarphiS99}) is very similar, and omitted.

\end{proof}


\begin{lemma}
\label{lem:sumaieiequalsvarphiS99} 
With reference to Definition \ref{def:setupforsectionS99}
and Definition
\ref{def:aidefS99}, 
pick any integer $i$ $(1\leq i \leq d)$. Then 
the scalar $\varphi_i$ equals each of the following four expressions.
\begin{eqnarray}
(\theta^*_i-\theta^*_{i-1})\sum_{h=0}^{i-1} (\theta_h-a_h),
&&\qquad \quad    
(\theta^*_{i-1}-\theta^*_i)\sum_{h=i}^{d} (\theta_h-a_h),
\qquad \qquad 
\label{eq:firsttwosumsgivingvarphiS99}
\\
(\theta_i-\theta_{i-1})\sum_{h=0}^{i-1} (\theta^*_h-a^*_h),
&&\qquad \quad  
(\theta_{i-1}-\theta_i)\sum_{h=i}^{d} (\theta^*_h-a^*_h).
\label{eq:secondtwosumsgivingvarphiS99}
\end{eqnarray}
\end{lemma}

\begin{proof} 
To see 
$\varphi_i$ equals
the expression
on the left in 
(\ref{eq:firsttwosumsgivingvarphiS99}),
eliminate
each of $a_0, a_1, \ldots, a_{i-1}$ in that expression 
using 
(\ref{eq:aiintermsofvarphiS99}), and simplify.
The 
two expressions  
in (\ref{eq:firsttwosumsgivingvarphiS99}) 
 are equal  by 
(\ref{eq:sumeigequalssumaiAS99}).
To see
 $\varphi_i$ equals
the expression
on the left in 
(\ref{eq:secondtwosumsgivingvarphiS99}),
eliminate
each of $a^*_0, a^*_1, \ldots, a^*_{i-1}$ in that expression
using 
(\ref{eq:aisintermsofvarphiS99}), and simplify.
The  two expressions in 
(\ref{eq:secondtwosumsgivingvarphiS99})
are equal by
(\ref{eq:sumeigequalssumaistarAS99}).

\end{proof}

\section{The $D_4$ action}

\begin{definition}
\label{def:secsetd4S99}
In this section,  $d$ will denote
a nonnegative integer, $\fld$ will denote a field, and
$\alg$ will denote an $\fld$-algebra isomorphic to
$\Mdf$. We let  
\begin{equation}
\ls
\label{eq:presentlsS99}
\end{equation}
denote a Leonard system in $\alg$,
with eigenvalue sequence $\theta_0, \theta_1, \ldots, \theta_d$, 
dual eigenvalue sequence 
$\theta^*_0, \theta^*_1, \ldots, \theta^*_d$, 
$\varphi$-sequence $\varphi_1, \varphi_2, \ldots, \varphi_d$,
and $\phi$-sequence 
 $\phi_1, \phi_2, \ldots, \phi_d$.
\end{definition}

\noindent With reference to 
Definition
\ref{def:secsetd4S99},
we now consider what happens to
the $\theta_i, \theta^*_i, \varphi_i, \phi_i$  
when   
$\Phi$ is replaced 
by a relative. We begin with a simple observation.

\begin{lemma}
\label{def:theaiaisnewdefS99} 
With reference to
Definition \ref{def:secsetd4S99} and 
Definition 
\ref{def:aidefS99}, 
\begin{eqnarray}
&&\theta^\downarrow_i = \theta_i, \qquad \qquad 
\theta^\Downarrow_i = \theta_{d-i},
\label{eq:d4onthS99}
\\
&&a^\downarrow_i = a_{d-i}, \qquad \qquad 
a^\Downarrow_i = a_i,
\label{eq:d4onaiS99}
\end{eqnarray}
for $0 \leq i \leq d$. (We are using the notation of Definition
\ref{def:d4actionlsS99}).
\end{lemma}

\begin{proof} Recall $\theta_i$ is the eigenvalue of $A$ associated
with $E_i$, and  $a_i$ is the trace of $AE^*_i$.
By 
Definition
\ref{def:d4actionlsS99}, we find that for 
all $g \in D_4$,  $\theta^g_i$ is the eigenvalue
of $A^g$ associated with $E^g_i$, and $a^g_i$ is the trace of 
$A^gE^{*g}_i$. By 
(\ref{eq:lsinvertS99}) we have $A^\downarrow =A$,
$E^\downarrow_i=E_i$, and $E^{*\downarrow}_i = E^*_{d-i}$.
By 
(\ref{eq:lsdualinvertS99}) we have $A^\Downarrow =A$,
$E^\Downarrow_i=E_{d-i}$, and $E^{*\Downarrow}_i = E^*_i$.
The result follows.

\end{proof}

\begin{lemma}
\label{lem:sumaieiequalsvarphiLSS99} 
With reference to
Definition \ref{def:secsetd4S99}
and Definition
\ref{def:aidefS99},  
 for $1\leq i \leq d$, 
the scalar $\varphi_i$ equals each of the following four expressions. 
\begin{eqnarray}
(\theta^*_i-\theta^*_{i-1})\sum_{h=0}^{i-1} (\theta_h-a_h),
&&\qquad \quad    
(\theta^*_{i-1}-\theta^*_i)\sum_{h=i}^{d} (\theta_h-a_h),
\qquad \qquad 
\label{eq:firsttwosumsgivingvarphiLSS99}
\\
(\theta_i-\theta_{i-1})\sum_{h=0}^{i-1} (\theta^*_h-a^*_h),
&&\qquad \quad  
(\theta_{i-1}-\theta_i)\sum_{h=i}^{d} (\theta^*_h-a^*_h).
\label{eq:secondtwosumsgivingvarphiLSS99}
\end{eqnarray}
\end{lemma}

\begin{proof} Apply Lemma 
\ref{lem:sumaieiequalsvarphiS99}
to the split canonical form of $\Phi$. 

\end{proof}

\begin{lemma}
\label{lem:sumaieiequalsphiS99} 
With reference to
Definition \ref{def:secsetd4S99}
and 
Definition \ref{def:aidefS99},  
 for $1\leq i \leq d$, 
the scalar $\phi_i$ equals each of the following four expressions.
\begin{eqnarray}
(\theta^*_{i}-\theta^*_{i-1}) \sum_{h=0}^{i-1}(\theta_{d-h}-a_h),
&&\quad  
(\theta^*_{i-1}-\theta^*_{i})
\sum_{h=i}^d (\theta_{d-h}-a_h), 
\label{eq:phifirstsumS99}
\\
(\theta_{d-i}-\theta_{d-i+1}) \sum_{h=0}^{i-1}(\theta^*_{h}-a^*_{d-h}),
&&\;\; 
(\theta_{d-i+1}-\theta_{d-i})
\sum_{h=i}^{d} (\theta^*_{h}-a^*_{d-h}). \qquad \qquad 
\label{eq:phisecondsumS99}
\end{eqnarray}

\end{lemma}

\begin{proof} By definition $\phi_i=\varphi^\Downarrow_i$. To find
$\varphi^\Downarrow_i$,
apply
Lemma 
 \ref{lem:sumaieiequalsvarphiLSS99} 
to $\Phi^\Downarrow$, and  evaluate the result
using
(\ref{eq:deightrelationsAS99}), 
(\ref{eq:deightrelationsBS99}), and
Lemma 
\ref{def:theaiaisnewdefS99}. 

\end{proof}

\noindent We are now ready to prove 
Theorem
\ref{thm:newrepLSintroS99} from the Introduction. 

%
%

\medskip
\noindent {\it Proof of Theorem
\ref{thm:newrepLSintroS99}}: Referring to the table in the theorem statement, 
the $\theta^g_i$ and $\theta^{*g}_i$
are readily obtained  using  
(\ref{eq:deightrelationsAS99}), 
(\ref{eq:deightrelationsBS99}), and
(\ref{eq:d4onthS99}).
Comparing the expressions on the left in
(\ref{eq:firsttwosumsgivingvarphiLSS99}), 
(\ref{eq:secondtwosumsgivingvarphiLSS99}), we 
 see $\varphi^*_i =\varphi_i$.
By definition
$\varphi^\Downarrow_i =\phi_i$, and recall $\Downarrow$
is an involution, so 
$\phi^\Downarrow_i =\varphi_i$.
Applying $*$ to the expression on the left in 
(\ref{eq:phifirstsumS99}), and 
replacing $i$ by $d-i+1$ in the result, we 
get the expression on the right in
(\ref{eq:phisecondsumS99});
it follows
 $\phi^*_i =\phi_{d-i+1}$.
Applying $\downarrow $ to the expression on
the left in 
(\ref{eq:firsttwosumsgivingvarphiLSS99}),
and 
replacing $i$ by $d-i+1$ in the result, we 
get the expression
on the right in 
(\ref{eq:phifirstsumS99}); it follows
$\varphi_i^\downarrow =  
\phi_{d-i+1}$.  
By this and since $\downarrow $ is an involution,
we find 
 $\phi_i^\downarrow = 
\varphi_{d-i+1}$.
The remaining entries of the table are routinely obtained 
using (\ref{eq:deightrelationsAS99}), 
(\ref{eq:deightrelationsBS99}).
\par\noindent $\Box$\par




\medskip
\noindent 
We interpret the data in   
Theorem 
\ref{thm:newrepLSS99} as follows. 
Let $\Phi$ be as in 
Definition \ref{def:secsetd4S99}.
By the {\it parameter square} of $\Phi$, we mean
the diagram below.

\medskip
\centerline{
         \begin{tabular}[t]{|c c c c c c c c c c c|}\hline
	$\theta^*_0$ & & & & &$ \phi_1$ & & & & & $\theta_d$ \\ 
	& $\theta^*_1$  & & & &$ \phi_2$ & & & &  $\theta_{d-1}$ & \\ 
	&  & $\cdot$ & & & $ \cdot $ & &  &  $\cdot$ &  &\\ 
	&  & & $\cdot$ & & $ \cdot $ &   &  $\cdot$ & &  &\\ 
	&  & &  & &  &   &   & &  &\\ 
	&  & &  & &  &   &   & &  &\\ 
	$\varphi_1$& $\varphi_2$  & $\cdot$ & $\cdot$  & &  & 
	&  $\cdot $ &$\cdot$ & $\varphi_{d-1}$ & $\varphi_d$\\ 
	&  & &  & &  &   &   & &  &\\ 
	&  & &  & &  &   &   & &  &\\ 
	&  & & $\cdot$ & & $ \cdot $ &   &  $\cdot$ & &  &\\ 
	&  & $\cdot$ & & & $ \cdot $ & &  &  $\cdot$ &  &\\ 
	& $\theta_1$  & & & &$ \phi_{d-1}$ & & & &  $\theta^*_{d-1}$ & \\ 
	$\theta_0$ & & & & &$ \phi_{d}$ & & & & & $\theta^*_d$ \\ \hline
	\end{tabular}
}

\medskip

\medskip

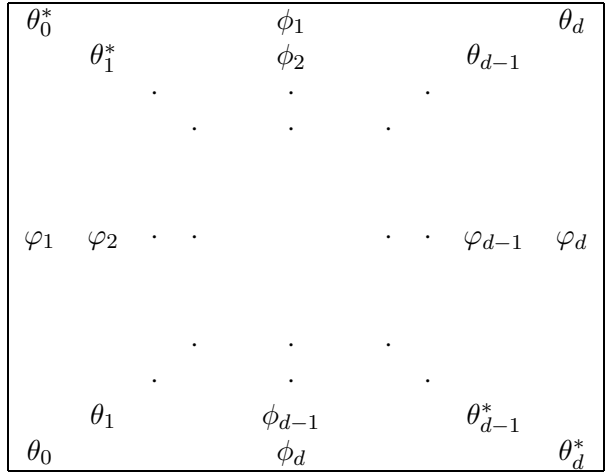
\centerline{Fig. 1.  The parameter square of a Leonard system}  

\bigskip
\noindent To get the parameter square for $\Phi^*$,
reflect the parameter square of $\Phi$
about the horizontal line running through the center.
To get the parameter square of $\Phi^\downarrow$,  reflect 
 the parameter square of $\Phi$
about the  diagonal running from the bottom  left corner
to the top right corner.
To get the parameter square of $\Phi^\Downarrow$,  reflect 
 the parameter square of $\Phi$
about the  diagonal running from the bottom right corner
to the top left corner. 

\medskip
\noindent We finish this section with a comment.

%
%
%

\begin{lemma}
\label{lem:varphiminusphiS99}
With reference to Definition \ref{def:secsetd4S99},
for $1 \leq i \leq d$, 
\begin{enumerate}
\item
${\displaystyle{
\varphi_i-\phi_i \;=\;(\theta^*_i-\theta^*_{i-1})\sum_{h=0}^{i-1}(\theta_h-\theta_{d-h}),  
}}$
\item
${\displaystyle{
\varphi_i-\phi_{d-i+1} \;=\;(\theta_i-\theta_{i-1})
\sum_{h=0}^{i-1}(\theta^*_h-\theta^*_{d-h}).  
}}$
\end{enumerate}
\end{lemma}

\begin{proof}(i) Subtract the expression on the left in 
(\ref{eq:phifirstsumS99})
from the 
expression on the left in 
(\ref{eq:firsttwosumsgivingvarphiLSS99}).

\noindent (ii) Apply (i) above to $\Phi^*$.

\end{proof}

\section{A result on reducibility}


\noindent In this section, we return to the situation of 
Definition
\ref{def:setupforsectionS99}.
We extend the domain of definition of the  $\phi_i$ scalars
to the level of 
Definition
\ref{def:setupforsectionS99}
using 
Lemma 
\ref{lem:varphiminusphiS99}(i). We 
 use the resulting constants
to get a  reducibility result.
\begin{definition}
\label{def:extenddefphiS99}  
With reference to Definition
\ref{def:setupforsectionS99}, we define
\begin{equation}
\phi_i=\varphi_i-(\theta^*_i-\theta^*_{i-1})\sum_{h=0}^{i-1}(\theta_h-\theta_{d-h})
\qquad \qquad (1 \leq i \leq d).
\label{eq:extenddefphiS99}
\end{equation}
For notational convenience, we set $\phi_0=0$, $\phi_{d+1}=0$.
\end{definition}

\begin{lemma}
\label{lem:whatifsubmodS99}  
With reference to Definition
\ref{def:setupforsectionS99},
let $V=\fld^{d+1}$ denote the irreducible  left module for $\Mdf$,
and let 
 $W$ denote a nonzero $(A, A^*)$-module in $V$.
  Then there exists an integer 
$r$ $(0 \leq r \leq d)$ such that both
\begin{eqnarray}
W = \sum_{h=r}^d E^*_hV,\qquad \qquad 
W = \sum_{h=0}^{d-r} E_hV. 
\label{eq:modulefirstdecS99}
\label{eq:moduleseconddecS99}
\end{eqnarray}
Moreover, the scalar $\phi_r$ from  
Definition
\ref{def:extenddefphiS99}  
is zero.
\end{lemma}

\begin{proof} 
Since $W$ is nonzero and $A^*W\subseteq W$, there exists a
nonempty subset
$S^*$ of $\lbrace 0,1,\ldots, d\rbrace $ such that 
$W = \sum_{i\in S^*} E^*_iV$. 
Recall by Lemma 
\ref{lem:commentoneis1S99}(ii) that $E^*_{i+1}AE^*_i\not=0$
for $0 \leq i \leq d-1$.
Combining this 
with Lemma 
\ref{lem:modcharstarS99}(ii), 
we 
find $i \in S^* $ implies $i+1 \in  S^*$ for
$0 \leq i \leq d-1$. It follows     
$S^*=
\lbrace r,r+1,\ldots, d\rbrace $  for some 
integer $r$ $(0 \leq r \leq d)$.
Since $W$ is nonzero and $AW\subseteq W$,
 there exists a nonempty subset
$S$ of $\lbrace 0,1,\ldots, d\rbrace $ such that 
$W = \sum_{i\in S} E_iV$. 
Recall by Lemma 
\ref{lem:commentoneis1S99}(i)
that 
$E_{i-1}A^*E_{i}\not=0$
for $1 \leq i \leq d$.
Combining this with
 Lemma 
\ref{lem:modcharS99}(ii), 
 we 
find $i \in S $ implies $i-1 \in  S$ for
$1 \leq i \leq d$. It follows     
$S=
\lbrace 0,1,\ldots, s\rbrace $  for some 
integer $s$ $(0 \leq s \leq d)$. 
Considering the dimension of $W$ we find $|S| =|S^*|$,
so $s=d-r$, and
(\ref{eq:modulefirstdecS99}) follows.
It remains to show $\phi_r=0$. This holds by definition if
$r=0$, so assume $r\geq 1$. To get $\phi_r=0$ in this case,
 we first show 
\begin{equation}
a_r+a_{r+1}+\cdots + a_d = \theta_0+\theta_1+\cdots + \theta_{d-r}.
\label{eq:steponetoshowphiS99}
\end{equation}
For convenience, we abbreviate 
$E=\sum_{h=0}^{d-r}E_h$ and 
$E^*=\sum_{h=r}^{d}E^*_h$. 
We  show $AE$ and $AE^*$ have the same trace.
To do this, we  
put $X=A(E-E^*)$, 
and show $X$ has trace 0.
In fact $X^2=0$. To see this, 
we show $XV\subseteq W$ and $XW=0$. 
Each of $EV$, $E^*V$ equals $W$
by (\ref{eq:modulefirstdecS99}),
so $(E-E^*)V\subseteq W$.
Recall $AW\subseteq W$, so $XV\subseteq W$.
Observe each of $E$, $E^*$
acts as the identity on $W$, so $(E-E^*)W=0$, and it follows
$XW=0$.  We have now shown $X^2=0$, so  $X$ has trace 0,
and $AE$, $AE^*$ have the same trace.
We now compute these traces. By
(\ref{eq:primid1S99}) and since each $E_h$ has trace 1, we find
 $AE$ has trace  $\sum_{h=0}^{d-r}\theta_h$.
Using Definition
\ref{def:aidefS99}, we routinely find
$AE^*$ has trace 
 $\sum_{h=r}^{d}a_h$. 
We now have
(\ref{eq:steponetoshowphiS99}).
Eliminating  the left side of 
(\ref{eq:steponetoshowphiS99})
using
the equation on the right in
(\ref{eq:firsttwosumsgivingvarphiS99}), we find $\phi_r=0$.

\end{proof}

\begin{theorem}
\label{thm:submodcondS99}  
With reference to Definition
\ref{def:setupforsectionS99},
let $V=\fld^{d+1}$ denote the irreducible left module for $\Mdf$, 
and suppose the scalars $\phi_1, \phi_2,\ldots, \phi_d$ from
 (\ref{eq:extenddefphiS99})
are all nonzero. Then 
 $V$ is irreducible as an 
$(A,A^*)$-module.

\end{theorem}

\begin{proof} Let $W$ denote a nonzero $(A,A^*)$-module 
in $V$. We show $W=V$. Let $r$ denote the  integer
associated with $W$ from Lemma
\ref{lem:whatifsubmodS99}.
From that lemma and our present assumption we find
$r$ is not one of $1,2,\ldots, d$, so $r=0$.
 Setting $r=0$ in
(\ref{eq:modulefirstdecS99}), we find $W=V$.

\end{proof}

\section{Recurrent sequences } 

\noindent  It is going to turn out that the eigenvalue
sequence and dual eigenvalue sequence of a Leonard system
each satisfy a certain recurrence.
In this section, we set the stage by considering
this recurrence from several points of view.

\begin{definition}
\label{def:recseqsettupS99}
In this section, $\fld $ will denote a field,
$d$ will denote a nonnegative integer,  and 
 $\; \theta_0, \theta_1, \ldots,  \theta_d\; $  will denote a sequence
 of scalars taken from $\fld$.

\end{definition}

\begin{definition}
\label{lem:beginthreetermS99}
With reference to Definition
\ref{def:recseqsettupS99}, let $\beta, \gamma, \varrho$ denote scalars
in $\fld$.
\begin{enumerate} 
\item The sequence 
  $\; \theta_0, \theta_1, \ldots,  \theta_d \;$ 
is said to be recurrent whenever $\theta_{i-1}\not=\theta_i$ for
$2 \leq i \leq d-1$, and 
\begin{equation}
{{\theta_{i-2}-\theta_{i+1}}\over {\theta_{i-1}-\theta_i}} 
\label{eq:thethingwhichisbetaS99}
\end{equation}
is independent of
$i$, for $2 \leq i \leq  d-1$.
\item The sequence 
  $\; \theta_0, \theta_1, \ldots,  \theta_d \;$ 
is said to be $\beta$-recurrent whenever 
\begin{equation}
\theta_{i-2}\,-\,(\beta+1)\theta_{i-1}\,+\,(\beta +1)\theta_i \,-\,\theta_{i+1}
\label{eq:betarecS99}
\end{equation}
is zero for 
$2 \leq i \leq d-1$.
\item The sequence 
  $\; \theta_0, \theta_1, \ldots,  \theta_d \;$ 
is said to be $(\beta,\gamma)$-recurrent whenever 
\begin{equation}
\theta_{i-1}\,-\,\beta \theta_i\,+\,\theta_{i+1}=\gamma 
\label{eq:gammathreetermS99}
\end{equation}
 for 
$1 \leq i \leq d-1$.
\item The sequence 
  $\; \theta_0, \theta_1, \ldots,  \theta_d \;$ 
is said to be $(\beta,\gamma,\varrho)$-recurrent whenever 
\begin{equation}
\theta^2_{i-1}-\beta \theta_{i-1}\theta_i+\theta^2_i 
-\gamma (\theta_{i-1} +\theta_i)=\varrho
\label{eq:varrhothreetermS99}
\end{equation}
 for 
$1 \leq i \leq d$.
\end{enumerate}
\end{definition}

\begin{lemma} 
\label{lem:recvsbrecS99}
With reference to Definition
\ref{def:recseqsettupS99}, the following are equivalent.
\begin{enumerate}
\item The sequence 
  $\; \theta_0, \theta_1, \ldots,  \theta_d \;$ 
is  recurrent.
\item There exists $\beta \in \fld$ such that
  $\; \theta_0, \theta_1, \ldots,  \theta_d \;$  is $\beta$-recurrent,
and $\theta_{i-1}\not=\theta_i$ for
$2 \leq i \leq d-1$.
\end{enumerate}
\noindent Suppose (i), (ii), and that $d\geq 3$. Then
the common value of
(\ref{eq:thethingwhichisbetaS99}) 
equals $\beta +1$.

\end{lemma}
\begin{proof} Routine.

\end{proof}

\begin{lemma}
\label{lem:brecvsbgrecS99}
With reference to Definition
\ref{def:recseqsettupS99}, the
following are equivalent for all $\beta \in \fld$.
\begin{enumerate}
\item The sequence 
  $\; \theta_0, \theta_1, \ldots,  \theta_d \;$ 
is  $\beta$-recurrent.
\item There exists $\gamma \in \fld$ such that
  $\; \theta_0, \theta_1, \ldots,  \theta_d \;$  
  is $(\beta,\gamma)$-recurrent.
\end{enumerate}
\end{lemma}
\begin{proof} 
$(i)\rightarrow (ii) $ 
For $2\leq i \leq d-1$, the expression
(\ref{eq:betarecS99}) is zero by assumption,
so
\beast
\theta_{i-2}\,-\,\beta \theta_{i-1}\,+\,\theta_i \;= \;
\theta_{i-1}\,-\,\beta \theta_i\,+\,\theta_{i+1}.
\eeast
Apparently the left side of 
(\ref{eq:gammathreetermS99}) is independent of $i$, and
the result follows.

\noindent 
$(ii)\rightarrow (i) $  Subtracting the equation 
(\ref{eq:gammathreetermS99}) at $i$ from the corresponding equation
obtained by replacing $i$ by $i-1$, we find
(\ref{eq:betarecS99}) is zero 
for $2\leq i \leq d-1$.

\end{proof}

\begin{lemma}
\label{lem:bgrecvsbgdrecS99}
With reference to Definition
\ref{def:recseqsettupS99}, the following (i),(ii) hold
for all $\beta, \gamma \in \fld$.
\begin{enumerate}
\item  Suppose 
  $\; \theta_0, \theta_1, \ldots,  \theta_d \;$ 
is  $(\beta,\gamma)$-recurrent. Then 
there exists $\varrho \in \fld$ such that
  $\; \theta_0, \theta_1, \ldots,  \theta_d \;$  is $(\beta,\gamma,\varrho)$-recurrent.
\item  Suppose 
  $\; \theta_0, \theta_1, \ldots,  \theta_d \;$ 
is  $(\beta,\gamma,\varrho)$-recurrent, and that $\theta_{i-1}\not=\theta_{i+1}$
for $1 \leq i\leq d-1$. Then
  $\; \theta_0, \theta_1, \ldots,  \theta_d \;$  is $(\beta,\gamma)$-recurrent.
\end{enumerate}
\end{lemma}

\begin{proof} 
Let  $p_i$ denote the expression on the left in
(\ref{eq:varrhothreetermS99}),
and observe
 \beast
p_i-p_{i+1} &=& 
(\theta_{i-1}-\theta_{i+1})(\theta_{i-1}-\beta \theta_i +\theta_{i+1} - \gamma)
\eeast
for $1 \leq i \leq d-1$. 
Assertions (i), (ii) are both routine consequences of this.

\end{proof}

\section{Recurrent sequences in closed form } 

\noindent In this section,  we obtain some formula
involving recurrent sequences.

\begin{definition}
\label{def:brecseqsettupS99}
In this section, $\fld $ will denote a field,
$d$ will denote a nonnegative integer,  and 
 $\; \beta, \theta_0, \theta_1, \ldots,  \theta_d\; $  will denote 
 scalars in $\fld$ such that 
 $\; \theta_0, \theta_1, \ldots,  \theta_d\; $ is $\beta$-recurrent. 
We   let ${\cal F}^{cl}$ denote the algebraic closure of
$\fld$. For all $q \in 
{\cal F}^{cl}$, we let 
$\fld\lbrack q \rbrack $  denote the field extention of
$\fld$ generated by $q$.

\end{definition}

\begin{lemma}
\label{lem:closedformthreetermS99}
With reference to Definition \ref{def:brecseqsettupS99}, 
the following  (i)--(iv) hold.
\begin{enumerate}
\item  Suppose $\beta \not=2$, $\beta \not=-2$, and pick
$q \in 
{\cal F}^{cl}$ such that 
 $q+q^{-1}=\beta $. Then there exists  scalars 
 $\alpha_1, \alpha_2, \alpha_3$  in  
$\fld\lbrack q \rbrack $ such that
\begin{equation}
\theta_i = \alpha_1 + \alpha_2 q^i + \alpha_3 q^{-i}
\qquad \qquad (0 \leq i \leq d). 
\label{eq:closedformthreetermIS99}
\end{equation}
\item Suppose $\beta = 2$ and $\hbox{char}(\fld) \not=2$. Then there
exists 
 $\alpha_1, \alpha_2, \alpha_3 $  in $\fld $ such that
\begin{equation}
\theta_i = \alpha_1 + \alpha_2 i + \alpha_3 i^2 
\qquad \qquad (0 \leq i \leq d).  
\label{eq:closedformthreetermIIS99}
\end{equation}
\item Suppose $\beta = -2$ and  $\hbox{char}(\fld) \not=2$. Then there
exists 
 $\alpha_1, \alpha_2, \alpha_3 $  in $\fld $ such that
\begin{equation}
\theta_i = \alpha_1 + \alpha_2 (-1)^i + \alpha_3 i(-1)^i 
\qquad \qquad (0 \leq i \leq d).  
\label{eq:closedformthreetermIIIS99}
\end{equation}
\item Suppose $\beta = 0$ and  $\hbox{char}(\fld) =2$. Then there
exists 
 $\alpha_1, \alpha_2, \alpha_3 $  in $\fld $ such that
\begin{equation}
\theta_i = \alpha_1 + \alpha_2 i + \alpha_3 \Biggl( {{i }\atop {2}}\Biggr) 
\qquad \qquad (0 \leq i \leq d),  
\label{eq:closedformthreetermIVS99}
\end{equation}
where we interpret the binomial coefficient as follows:  
\beast
 \Biggl( {{i }\atop {2}}\Biggr) 
 = \left\{ \begin{array}{ll}
                   0  & \mbox{if $\;i=0\;$ or $\;i=1\;$ (mod $4$), } \\
				  1 & \mbox{if $\;i=2 \;$ or $\;i=3\;$ (mod $4$). }
				   \end{array}
				\right. 
\eeast
\end{enumerate}
\end{lemma}

\begin{proof} (i).
We assume $d\geq 2$; otherwise the result is trivial.
Let $q$ be given, and consider the equations 
(\ref{eq:closedformthreetermIS99}) for $i=0,1,2$. These
equations are linear in $\alpha_0, \alpha_1, \alpha_2$.
We routinely find the coefficient matrix is nonsingular,
so there exist
$\alpha_0, \alpha_1, \alpha_2$ in $\fld \lbrack q \rbrack $
such that
(\ref{eq:closedformthreetermIS99})  holds for $i=0,1,2$.
Using these scalars, let $\varepsilon_i $ denote the
left  side of 
(\ref{eq:closedformthreetermIS99}) minus the 
right  side of 
(\ref{eq:closedformthreetermIS99}), for $0 \leq i \leq d$.
On one hand, 
$\varepsilon_0$,
$\varepsilon_1$,
$\varepsilon_2$ are zero from the construction.
On the other hand,
one readily checks
\beast
\varepsilon_{i-2}\,-\,(\beta+1)\varepsilon_{i-1}\,+\,(\beta +1)\varepsilon_i \,-\,\varepsilon_{i+1}=0 
\eeast
for 
$2 \leq i \leq d-1$.
Combining these facts, we find
$\varepsilon_i=0$ for $0 \leq i \leq d$, and the result follows.

\noindent (ii)--(iv) Similar to the proof of (i) above.

\end{proof}

\begin{lemma}
\label{lem:closedformcommentthreetermS99}
With reference to Definition \ref{def:brecseqsettupS99}, 
assume 
$\;\theta_0, \theta_1, \ldots, \theta_d \;$ are distinct.
Then (i)--(iv) hold below.
\begin{enumerate}
\item  Suppose $\beta \not=2$, $\beta \not=-2$,
and pick
$q \in 
{\cal F}^{cl}$ such that 
 $q+q^{-1}=\beta $. Then  
$q^i \not=1 $ for 
$1 \leq i \leq d$.   
\item Suppose $\beta = 2$ and $\hbox{char}(\fld) \not=2$. Then 
 $\hbox{char}(\fld)=0$ or 
 $\hbox{char}(\fld)>d$. 
\item Suppose $\beta = -2$ and  $\hbox{char}(\fld) \not=2$. Then 
 $\hbox{char}(\fld)=0$ or 
 $\hbox{char}(\fld)>d/2$. 
\item Suppose $\beta = 0$ and  $\hbox{char}(\fld) =2$. Then $d\leq 3$. 

\end{enumerate}
\end{lemma}
\begin{proof} (i) Using
(\ref{eq:closedformthreetermIS99}), we find 
$q^i=1 $ implies $\theta_i=\theta_0$ for $1 \leq i \leq d$.

\noindent (ii) Using
(\ref{eq:closedformthreetermIIS99}), we find
that for
  $1 \leq i \leq d$,
if $i$ is congruent to $0$ modulo the characteristic of $\fld$,
then  
 $\theta_i=\theta_0$, a contradiction. The result follows.

\noindent (iii) Using  
(\ref{eq:closedformthreetermIIIS99}), we find that for 
any even integer $i$,  
 $(1 \leq i \leq d)$,
if $i$ is congruent to $0$ modulo the characteristic of $\fld$,
 then $\theta_i=\theta_0$, a contradiction. The result follows. 

\noindent (iv) Suppose $d\geq 4$. Applying  
(\ref{eq:closedformthreetermIVS99}), we find $\theta_0=\theta_4$, 
a contradiction.

\end{proof}

\begin{lemma}
\label{lem:symeigformulaS99}
With reference to Definition \ref{def:brecseqsettupS99}, 
assume $\; \theta_0, \theta_1, \ldots, \theta_d \;$ are  distinct. 
Pick any integers $i,j,r,s$ $(0 \leq i,j,r,s \leq d)$
and assume $\;i+j=r+s\;$, $\;r\not=s$.
Then  (i)--(v) hold below.
\begin{enumerate}
\item Suppose $\beta \not=2$, $\beta\not=-2$. Then
\begin{equation}
{{\theta_i-\theta_{j}}\over {\theta_r-\theta_s}}
\;= \;{{q^{i}-q^j}\over {q^r-q^s}},
\label{eq:symeiggencasenewthreetermS99}
\end{equation}
where  $q+q^{-1}=\beta $.
\item Suppose $\beta = 2$ and
 $\hbox{char}(\fld) \not=2$. 
 Then
\begin{equation}
{{\theta_i-\theta_{j}}\over {\theta_r-\theta_s}}
\;=\; {{i-j}\over {r-s}}.
\label{eq:symeigbeta2newthreetermS99}
\end{equation}
\item Suppose $\beta = -2$ and  $\hbox{char}(\fld)\not=2$.
  Then
\begin{equation}
{{\theta_i-\theta_{j}}\over {\theta_r-\theta_s}}
\;=\;
 \left\{ \begin{array}{ll}
            (-1)^{i+r} {{i-j}\over {r-s}}  & \mbox{if $\;i+j\;$ is even}, \\
	(-1)^{i+r}  & \mbox{if $\;i+j\;$ is odd.}
				   \end{array}
				\right. 
\end{equation}
\item Suppose $\beta = 0$ and
 $\hbox{char}(\fld) =2$. 
 Then
\begin{equation}
{{\theta_i-\theta_{j}}\over {\theta_r-\theta_s}}
\;=\; 
 \left\{ \begin{array}{ll}
            0  & \mbox{if $\;i=j$}, \\
	1  & \mbox{if $\;i\not=j$.}
				   \end{array}
				\right. 
\end{equation}

\end{enumerate}
\end{lemma}

\begin{proof} To get (i), evaluate the left side in 
(\ref{eq:symeiggencasenewthreetermS99}) using
(\ref{eq:closedformthreetermIS99}), and simplify the
result. The cases (ii)--(iv) are very similar.
%

\end{proof}
\noindent We finish this section with an observation.

\begin{lemma}
\label{lem:dualitypreservedS99} 
With the notation and 
assumptions of 
Lemma \ref{lem:symeigformulaS99}, the scalar
\beast
{{\theta_i-\theta_{j}}\over {\theta_r-\theta_s}}
\eeast
depends only on $i,j,r,s$ and $\beta$, and not
on $\theta_0, \theta_1, \ldots, \theta_d$.
\end{lemma}

\begin{proof}
This is immediate from 
the data in Lemma \ref{lem:symeigformulaS99}.

\end{proof}

\section{A sum}

\begin{definition}
\label{def:taupolyfreeagainS99}
Throughout this section,   
 $\fld $ will 
denote a field,
$d$ will denote an  integer at least 1, and
$\theta_0, \theta_1,\ldots, \theta_d$ will
denote a sequence of 
distinct scalars
in $\fld$. We let $\beta $ denote any scalar in $\fld$.
\end{definition}

\noindent  
With reference to 
Definition \ref{def:taupolyfreeagainS99},
 we now consider the sums
\begin{equation}
 \sum_{h=0}^{i-1} {{\theta_h-\theta_{d-h}}\over {\theta_0-\theta_d}},
\label{eq:omegaabovediagS99}
\end{equation}
where $0 \leq i \leq d+1$.
Denoting the sum in 
(\ref{eq:omegaabovediagS99}) by $\vartheta_i$, we remark 
\begin{equation}
 \vartheta_0=0,\qquad \vartheta_1=1,\qquad  \vartheta_d=1,\qquad
\vartheta_{d+1}=0.
\label{eq:varthetaprelimsS99}           
\end{equation}
Moreover 
\begin{equation}
\vartheta_i= \vartheta_{d-i+1} \qquad \qquad (0 \leq i \leq d+1), 
\label{eq:varthetaprelims2S99}
\end{equation}
and 
\begin{equation}
\vartheta_{i+1}-\vartheta_i={{\theta_{i}-\theta_{d-i}}
\over {\theta_0-\theta_{d}}}
\qquad \qquad (0 \leq i \leq d). 
\label{eq:easyrecAS99}
\end{equation}

\noindent 
It turns out the sums
(\ref{eq:omegaabovediagS99})
play an important role a bit later, so we will 
examine them carefully.
We  begin by giving explicit formulae for
the sums 
(\ref{eq:omegaabovediagS99}) under the assumption
the sequence 
$\theta_0, \theta_1,\ldots, \theta_d$ is recurrent. 
To avoid trivialities,  we assume $d\geq 3$. 
\begin{lemma} 
\label{lem:symeigvalsformulaS99}
 With reference 
to Definition
\ref{def:taupolyfreeagainS99},
assume $d\geq 3$,  
and assume 
$\theta_0, \theta_1,\ldots, \theta_d$ is $\beta$-recurrent. 
Then for all integers $i$ $(0 \leq i \leq d+1)$, we have the following.
\begin{enumerate}
\item Suppose $\beta \not=2$, $\beta\not=-2$. Then
\begin{equation}
 \sum_{h=0}^{i-1} 
{{\theta_h-\theta_{d-h}}\over {\theta_0-\theta_d}}
= {{(q^i-1)(q^{d-i+1}-1)}\over {(q-1)(q^d-1)}},
\label{eq:alphagencaseS99}
\end{equation}
where  $q+q^{-1}=\beta $.
\item Suppose $\beta = 2$ and 
$\hbox{char}(\fld) \not=2$. Then
\begin{equation}
 \sum_{h=0}^{i-1} 
{{\theta_h-\theta_{d-h}}\over {\theta_0-\theta_d}}
= {{i(d-i+1)}\over {d}}.
\label{eq:alphacasebeta2S99}
\end{equation}
\item Suppose $\beta = -2$,  
$\hbox{char}(\fld) \not=2$, and 
$\;d\;$ odd. Then
\begin{equation}
 \sum_{h=0}^{i-1} 
{{\theta_h-\theta_{d-h}}\over {\theta_0-\theta_d}}
=  \left\{ \begin{array}{ll}
                   0  & \mbox{if $\;i\; $ is even } \\
				  1 & \mbox{if $\;i\;$ is odd.}
				   \end{array}
				\right. 
\label{eq:alphacasebetamin2S99}
\end{equation}                     
\item Suppose $\beta = -2$, 
$\hbox{char}(\fld) \not=2$, and 
  $\;d\;$ even. Then 
\begin{equation}
 \sum_{h=0}^{i-1} 
{{\theta_h-\theta_{d-h}}\over {\theta_0-\theta_d}}
 = 
  \left\{ \begin{array}{ll}
                   i/d  & \mbox{if $\;i\; $ is even } \\
	(d-i+1)/d \quad 	   & \mbox{if $\;i\;$ is odd. }
				   \end{array}
				\right.  
\label{eq:alphacasebetamin2deveS99}
\end{equation}
\item Suppose $\beta = 0$, 
$\hbox{char}(\fld) =2$, and 
  $\;d=3$. Then 
\begin{equation}
 \sum_{h=0}^{i-1} 
{{\theta_h-\theta_{d-h}}\over {\theta_0-\theta_d}}
 = 
  \left\{ \begin{array}{ll}
                   0  & \mbox{if $\;i\; $ is even } \\
	1 \quad 	   & \mbox{if $\;i\;$ is odd. }
				   \end{array}
				\right.  
\label{eq:alphacasebeta0char2S99}
\end{equation}

\end{enumerate}
\end{lemma}

\begin{proof} 
The  above sums
can  be  computed directly from 
Lemma \ref{lem:symeigformulaS99}.

\end{proof}

\noindent 
We mention some recursions satisfied by the 
 sums
(\ref{eq:omegaabovediagS99}).

\begin{lemma}
\label{lem:shortcharofsumsAS99}
With reference to Definition
\ref{def:taupolyfreeagainS99},
 assume 
$\theta_0, \theta_1,\ldots, \theta_d$ is recurrent, 
 and 
put 
\begin{equation}
\vartheta_i = \sum_{h=0}^{i-1} 
{{\theta_h-\theta_{d-h}}\over {\theta_0-\theta_d}}
\qquad \qquad (0 \leq i \leq d+1).
\label{eq:defofvarthetareminderS99}
\end{equation}
Then (i),(ii) hold below.
\begin{enumerate}
\item 
${\displaystyle{
\vartheta_{i+1} = \vartheta_i\,
{{\theta_{i}-\theta_{d-1}}\over {\theta_{i-1}-\theta_d}}
+ 1 
\qquad \qquad (1 \leq i \leq d), 
}}$
\item 
${\displaystyle{
\vartheta_{i} = \vartheta_{i+1}\,
{{\theta_{i}-\theta_{1}}\over {\theta_{i+1}-\theta_0}}
+ 1 
\qquad \qquad (0 \leq i \leq d-1). 
}}$
\end{enumerate}
\end{lemma}

\begin{proof} (i) These equations are readily verified
case by case, using 
Lemma 
\ref{lem:symeigvalsformulaS99}.

\noindent (ii) Apply (i) above to the sequence $\theta_d,\theta_{d-1},
\ldots, \theta_0$, and use 
(\ref{eq:varthetaprelims2S99}).

\end{proof}

\begin{lemma}
\label{lem:charofsumsS99}
With reference to 
Definition
\ref{def:taupolyfreeagainS99},
assume
$\theta_0, \theta_1,\ldots, \theta_d$ is recurrent.
Let $r$ denote any integer in the range $1\leq r \leq d+1$, 
and suppose we are given scalars
 $\vartheta_1, \vartheta_2, \ldots,   
\vartheta_r $  in $\fld$ such that
\begin{equation}
\vartheta_{i+1} = \vartheta_i\,
{{\theta_{i}-\theta_{d-1}}\over {\theta_{i-1}-\theta_d}}
+ \vartheta_1 
\qquad \qquad (1 \leq i \leq r-1). 
\label{eq:startrecursionvarthS99}
\end{equation}
Then
\begin{equation}
\vartheta_i = \vartheta_1\,\sum_{h=0}^{i-1} 
{{\theta_h-\theta_{d-h}}\over {\theta_0-\theta_d}}
\qquad \qquad (1 \leq i \leq r).
\label{eq:wegetvarthetabackS99}
\end{equation}
\end{lemma}

\begin{proof} Define
\begin{equation}
\vartheta'_i = \vartheta_i- \vartheta_1\,
\sum_{h=0}^{i-1} 
{{\theta_h-\theta_{d-h}}\over {\theta_0-\theta_d}}
\label{eq:defvarthprimeS99}
\end{equation}
for $1 \leq i \leq r$, and observe
$\vartheta'_1=0$.
Combining   
Lemma 
\ref{lem:shortcharofsumsAS99}(i)
and 
(\ref{eq:startrecursionvarthS99}), we routinely find
\begin{equation}
\vartheta'_{i+1} = \vartheta'_i\,
{{\theta_{i}-\theta_{d-1}}\over {\theta_{i-1}-\theta_d}}
\qquad \qquad (1 \leq i \leq r-1).
\end{equation}
Apparently
$\vartheta'_i = 0 $ for $1 \leq i \leq r$, and the result
follows.

\end{proof}

\noindent We mention an identity that will be  useful later.
\begin{lemma}
\label{lem:factabtsumS99}
With reference to 
Definition
\ref{def:taupolyfreeagainS99},
assume $\theta_0, \theta_1,\ldots, \theta_d$ is recurrent.
Then
\begin{equation}
{{\theta_0-\theta_1+\theta_{i-1}-\theta_{i}}\over {\theta_0-\theta_{i}}}
\,\sum_{h=0}^{i-1} 
{{\theta_h-\theta_{d-h}}\over {\theta_0-\theta_d}}
=
{{\theta_0+\theta_{i-1}-\theta_{d-i+1}-\theta_d}\over {\theta_0-\theta_d}},
\qquad 
\label{lem:factabtsumAS99}
\end{equation}
for $ 1 \leq i \leq d$.   (Caution: the numerator on the far  
left in 
(\ref{lem:factabtsumAS99}) 
might be zero).
\end{lemma}

\begin{proof}
Add 
(\ref{eq:easyrecAS99}) and 
Lemma 
\ref{lem:shortcharofsumsAS99}(ii), solve the resulting equation
for $\vartheta_{i+1}$, and replace $i$ by $i-1$ in the result.

\end{proof}

\noindent Here is another recursion.

\begin{lemma}
\label{lem:varthetacharacS99}
With reference to 
Definition
\ref{def:taupolyfreeagainS99},
assume $\theta_0, \theta_1,\ldots, \theta_d$ is $\beta$-recurrent,
and put 
\begin{equation}
\vartheta_i = \sum_{h=0}^{i-1} 
{{\theta_h-\theta_{d-h}}\over {\theta_0-\theta_d}}
\qquad \qquad (0 \leq i \leq d+1).
\label{eq:remembervarthS99}
\end{equation}
Then the sequence
 $\vartheta_0, \vartheta_1, \ldots,   
\vartheta_{d+1} $   is $\beta$-recurrent.
%
\end{lemma}

\begin{proof} 
We show 
\begin{equation}
\vartheta_{i-2}\;-\;(\beta+1)\vartheta_{i-1} \;+\; 
(\beta+1)\vartheta_{i}\;-\;\vartheta_{i+1}
\label{vartheta3termS99}
\end{equation}
is zero for $2\leq i \leq d$.
First observe by (\ref{eq:gammathreetermS99})
that
\begin{equation}
\theta_{j-1}-\beta \theta_j+\theta_{j+1} \;=\;
\theta_{d-j-1}-\beta \theta_{d-j}+\theta_{d-j+1} \qquad \qquad
(1 \leq j \leq d-1).
\label{eq:wewillgetthereS99}
\end{equation}
Eliminating $\vartheta_{i-2}$, 
$\vartheta_{i-1}$, 
$\vartheta_{i}$, 
$\vartheta_{i+1}$ in 
(\ref{vartheta3termS99}) 
using 
(\ref{eq:remembervarthS99}), 
then cancelling terms where possible, and then simplifying the result 
using 
(\ref{eq:wewillgetthereS99}), we get zero.

\end{proof}
\noindent

\noindent For completness sake, we include a lemma concerning the 
converse to Lemma 
\ref{lem:varthetacharacS99}. We do not use the result,
so we will not dwell on the proof.

\begin{lemma}
\label{note:varthetacharacconvS99}
With reference to 
Definition
\ref{def:taupolyfreeagainS99},
assume $\theta_0, \theta_1,\ldots, \theta_d$ is $\beta$-recurrent.
Let  $\vartheta_0, 
 \vartheta_1, \ldots,$ $   
\vartheta_{d+1}$  denote a $\beta$-recurrent sequence of
scalars taken from  $\fld$, such that  
 $\vartheta_0=0$, 
$\vartheta_{d+1}=0$, and 
$\vartheta_1=\vartheta_d$.
Then
\beast
\vartheta_i = \vartheta_1\sum_{h=0}^{i-1} 
{{\theta_h-\theta_{d-h}}\over {\theta_0-\theta_d}}
\qquad \qquad (0 \leq i \leq d+1).
\eeast

\end{lemma}

\begin{proof} Routine calculation using
Lemma \ref{lem:closedformthreetermS99},
Lemma \ref{lem:closedformcommentthreetermS99},
and 
Lemma 
\ref{lem:symeigvalsformulaS99}.

\end{proof}

\section{Some equations involving the split canonical form}

\noindent In this section, we  return to the situation of
Definition
\ref{def:setupforsectionS99}, and determine 
when the products $E_dA^*E_i$ vanish
for $0 \leq i \leq d-2$.
We begin with a definition.

\begin{definition}
\label{def:varthetadefS99n}
With reference to 
Definition \ref{def:setupforsectionS99}, 
we define
\begin{equation}
\vartheta_i = \varphi_i-(\theta^*_i-\theta^*_0)(\theta_{i-1}-\theta_d)
\qquad  \qquad (1 \leq i \leq d),
\label{eq:varthetadefS99n}
\label{eq:defofvarphiprimeS99}
\end{equation}
and  
 $\vartheta_0=0$, $\vartheta_{d+1}=0$. We observe $\vartheta_1$ equals
 the scalar $\phi_1$ from
 Definition \ref{def:extenddefphiS99}.
\end{definition}
\noindent Our goal
in this section is to  
 prove the following 
theorem.

\begin{theorem}
\label{thm:whenproductsvanishS99}
With reference to 
Definition 
\ref{def:setupforsectionS99},
assume $d\geq 2$. Then the following are equivalent.
\begin{enumerate}
\item $E_dA^*E_i = 0 \qquad \qquad (0 \leq i \leq d-2)$.
\item  
${\displaystyle{
\vartheta_{i+1}= \vartheta_i\,{{\theta_i-\theta_{d-1}}\over {\theta_{i-1}-\theta_d}} + \vartheta_1
\qquad \qquad (1 \leq i \leq d-1),
}}$
\end{enumerate}
where $\vartheta_1, \vartheta_2, \ldots, \vartheta_d$ are
from 
(\ref{eq:defofvarphiprimeS99}).
\end{theorem}

\noindent 
To prove the above theorem, it 
is 
 advantageous 
to consider the linear combination
\beast
\sum_{i=0}^{d-2} E_dA^*E_i(\theta_i-\theta_{d-1}).
\eeast

\begin{lemma}
\label{lem:whyconsiderlincombS99}
With
reference to Definition \ref{def:setupforsectionS99}, assume $d\geq 2$.
Then 
\begin{eqnarray}
\sum_{i=0}^{d-2} E_dA^*E_i(\theta_i-\theta_{d-1}) = E_d(A^*-a^*_d I)(A-\theta_{d-1}I),
\label{eq:whyconsiderlincombAS99}
\end{eqnarray}
where the scalar $a^*_d $ is from
Definition \ref{def:aidefS99}.

\end{lemma}

\begin{proof} Since $E_d$ is a rank one idempotent,
we find 
$E_dA^*E_d$ is a scalar multiple of $E_d$. Taking the trace,
we find
\begin{equation}
E_dA^*E_d = a^*_d E_d.
\label{eq:alphameaningS99}
\end{equation}
We may now argue
\beast
&& E_d(A^*-a^*_dI)(A-\theta_{d-1}I)
\\
&& \qquad = \;
 E_d(A^*-a^*_d I)(A-\theta_{d-1}I)\sum_{i=0}^d E_i
 \qquad \qquad
 (\hbox{by }\; 
(\ref{eq:primid3S99}))
 \\
&& \qquad = \;
 \sum_{i=0}^d E_d(A^*-a^*_d I)E_i(\theta_i-\theta_{d-1})
 \qquad \qquad \quad
 (\hbox{by }\; 
(\ref{eq:primid1S99}))
\\
&& \qquad = \;
 \sum_{i=0}^{d-2} E_d(A^*-a^*_d I)E_i(\theta_i-\theta_{d-1}) 
 \qquad \qquad  \quad
 (\hbox{by }\; 
(\ref{eq:alphameaningS99}))
\\
&& \qquad = \;
 \sum_{i=0}^{d-2} E_dA^*E_i(\theta_i-\theta_{d-1})
 \qquad \qquad \qquad \qquad 
 (\hbox{by }\; 
(\ref{eq:primid2S99})),
\eeast
as desired.

\end{proof}

\begin{corollary}
\label{cor:onewaytoshowvanishS99}
With
reference to Definition \ref{def:setupforsectionS99}, assume $d\geq 2$,
and let the scalar $a^*_d $ be as in 
Definition \ref{def:aidefS99}.
Then
the following are equivalent.
\begin{enumerate}
\item $E_dA^*E_i=0 \qquad \qquad (0 \leq i \leq d-2)$.
\item 
$E_d(A^*-a^*_d I)(A-\theta_{d-1}I) =0$.
\end{enumerate}
\end{corollary}

\begin{proof}
$(i)\rightarrow (ii)$ Immediate from  
Lemma \ref{lem:whyconsiderlincombS99}.

\noindent
$(ii)\rightarrow (i)$  Multiply both sides
of 
(\ref{eq:whyconsiderlincombAS99})
on the right by each $E_0, E_1, \ldots, E_{d-2}$, and
simplify using 
(\ref{eq:primid2S99}).

\end{proof}

\begin{lemma}
\label{lem:entriesofaproductS99}
With 
reference to Definition \ref{def:setupforsectionS99},
assume $d\geq 2$, and  consider 
the matrix 
\begin{equation}
E_d(A^*-a^*_d I)(A-\theta_{d-1}I),
\label{eq:eastaraprodS99}
\end{equation}
where $a^*_d $ is from
Definition \ref{def:aidefS99}.
For the matrix 
(\ref{eq:eastaraprodS99}), 
 all entries in rows  $0,1,\ldots, d-1$ are zero.
The entries in the $d^{\hbox{th}}$ row of
(\ref{eq:eastaraprodS99}) 
are as follows.
For $0 \leq i \leq d$, the $di^{\hbox{th}}$ entry
of 
(\ref{eq:eastaraprodS99}) 
is 
$\tau_i(\theta_d)\tau_d(\theta_d)^{-1}$ times
\begin{equation}
\vartheta_{i+1} \,-\, \vartheta_i\,
{{\theta_i-\theta_{d-1}}\over {\theta_{i-1}-\theta_d}} 
\,-\,\vartheta_d,
\label{eq:dientryofprodS99}
\end{equation}
where the $\vartheta_j$ are from
Definition \ref{def:varthetadefS99n},
  and where
  $\theta_{-1}$ is an 
 indeterminant. 
\end{lemma}

\begin{proof} 
To obtain our first assertion,
observe by 
(\ref{eq:entriesofefreeS99}) that 
 in the matrix $E_d$, 
all entries in rows $0,1,\ldots, d-1$ are zero.
$E_d$ is the factor on the left in 
(\ref{eq:eastaraprodS99}).
By matrix 
multiplication we see that in 
(\ref{eq:eastaraprodS99}),
all entries in rows $0,1,\ldots, d-1$ are zero.
We now consider the $d^{\hbox{th}}$ row of 
(\ref{eq:eastaraprodS99}). To compute it, we recall
the 
$d^{\hbox{th}}$ row of  $E_d$.
By
(\ref{eq:entriesofefreeS99}),
we find the entries
\beast
(E_d)_{di} &=& {{\tau_i(\theta_d) }\over {\tau_d(\theta_d)}}
\\
&=&  {{1}\over{(\theta_d-\theta_i)(\theta_d-\theta_{i+1})\cdots (\theta_d-\theta_{d-1})}}
\eeast
for $0 \leq i \leq d$.
Multiplying out 
(\ref{eq:eastaraprodS99}) using this,
we routinely find
its $di^{\hbox{th}}$ entry  is
$\tau_i(\theta_d)\tau_d(\theta_d)^{-1}$ times 
\beast
(\theta^*_i-a^*_d)(\theta_i -\theta_{d-1}) \;+ \;
(\theta^*_{i+1}-a^*_d)(\theta_d -\theta_{i}) \;+\; \varphi_{i+1}
\;-\; \varphi_i \,{{\theta_i-\theta_{d-1}}\over {\theta_{i-1}-\theta_d}}
\eeast
for $0 \leq i \leq d$. Eliminating $a^*_d$ in 
the above line  
using 
(\ref{eq:aisintermsofvarphiS99}),
and 
eliminating $\varphi_i$, $\varphi_{i+1}, \varphi_d$ in the result using 
(\ref{eq:defofvarphiprimeS99}), we obtain 
(\ref{eq:dientryofprodS99}).
\end{proof}

%
%

\noindent {\it Proof of Theorem 
\ref{thm:whenproductsvanishS99}}:
$\;(i)\rightarrow (ii)$ 
Apparently Corollary
\ref{cor:onewaytoshowvanishS99}(i) holds, so 
Corollary
\ref{cor:onewaytoshowvanishS99}(ii) holds, 
and  
the matrix
(\ref{eq:eastaraprodS99}) is zero.
Applying
Lemma 
\ref{lem:entriesofaproductS99}, 
we find
 the expressions 
(\ref{eq:dientryofprodS99})
are zero for $0 \leq i \leq d$.
Setting $i=0$, $\vartheta_0=0$ in 
(\ref{eq:dientryofprodS99}), we find
$\vartheta_1=\vartheta_d$. 
Using this to eliminate 
 $\vartheta_d$ in 
(\ref{eq:dientryofprodS99}), we obtain the equations
in
Theorem
\ref{thm:whenproductsvanishS99}(ii).

\noindent 
$(ii)\rightarrow (i)$ Setting $i=d-1$ in 
the given equations,
we find $\vartheta_d = \vartheta_1$. Using this to eliminate
$\vartheta_1$ on the right in the remaining given equations,
we 
find the expressions 
(\ref{eq:dientryofprodS99}) are zero for 
 $1 \leq i \leq d-1$.
We routinely find the expression 
(\ref{eq:dientryofprodS99}) is zero for $i=0$ and $i=d$,
so  (\ref{eq:dientryofprodS99}) is zero for $0 \leq i \leq d$.
Applying Lemma
\ref{lem:entriesofaproductS99}, we find the matrix
(\ref{eq:eastaraprodS99}) is zero.
Applying Corollary
\ref{cor:onewaytoshowvanishS99}, we find
$E_dA^*E_i$ vanishes for $0 \leq i \leq d-2$.
\par\noindent $\Box$\par

\section{Two polynomial equations for $A$ and $A^*$}

\noindent In this section, we show the elements $A$ and $A^*$ in
a Leonard pair satisfy two cubic polynomial equations.
We begin with a comment on the situation of Definition
\ref{def:verygeneralS99}.

\begin{lemma}
\label{lem:antisymS99n}
With reference 
to 
Definition \ref{def:verygeneralS99}, 
suppose 
\begin{equation}
E_iA^*E_j = 0 \quad \hbox{if} \quad |i-j|>1, \qquad \qquad (0 \leq i,j\leq d),
\label{eq:makingthingsworkS99n}
\end{equation}
and let 
${\cal D}$ denote the subalgebra of $\alg$ generated by $A$.
Then
\beast
\Span{XA^*Y-YA^*X \;|\;X, Y \in {\cal D}} \,=\,\lbrace XA^*-A^*X\;|\;X \in 
{\cal D}\rbrace. \qquad \quad 
\eeast

\end{lemma}

\begin{proof} For notational convenience set $E_{-1}=0$, $E_{d+1}=0$.
We claim that 
for $0 \leq i \leq d$, 
\begin{equation}
E_iA^*E_{i+1} - E_{i+1}A^*E_i \;=\;L_i A^* - A^* L_i,
\label{eq:liconversionS99n}
\end{equation}
where  $L_i := E_0+E_1+\cdots +E_i $.
To  see 
(\ref{eq:liconversionS99n}), observe 
by
(\ref{eq:primid3S99})
and 
(\ref{eq:makingthingsworkS99n})
that for  $0 \leq j \leq d$,  both
\begin{eqnarray}
E_jA^*&=&E_{j}A^*E_{j-1}+E_jA^*E_j + E_{j}A^*E_{j+1},
\label{eq:Aexpand1S99n}
\\
A^*E_j&=& E_{j-1}A^*E_j+E_jA^*E_j + E_{j+1}A^*E_j.
\label{eq:Aexpand2S99n}
\end{eqnarray}
Summing  
(\ref{eq:Aexpand1S99n})
over 
 $j=0,1,\ldots,i$, summing 
(\ref{eq:Aexpand2S99n}) over 
 $j=0,1,\ldots,i$, and taking the difference between the
 two sums, 
we obtain 
(\ref{eq:liconversionS99n}).
Observe 
$\cal D$ is spanned by both $E_0, E_1, \ldots, E_d$
and $L_0, L_1, \ldots,L_d$, so 
\beast
&&\Span{XA^*Y-YA^*X \;|\;X, Y \in {\cal D}}
\\
&& \qquad \qquad =\; \Span{E_iA^*E_j-E_jA^*E_i\;|\;
0 \leq i,j\leq d}\\
 && \qquad \qquad   
=\; \Span{E_iA^*E_{i+1}-E_{i+1}A^*E_i\;|\;
0 \leq i\leq d}\\
 &&\qquad \qquad =\; \Span{L_iA^*-A^*L_i\;|\;0 \leq i \leq d} \\
&&\qquad \qquad =\;\lbrace XA^*-A^*X\;|\;X \in 
{\cal D}\rbrace,
\eeast
and we are done.
\end{proof}

\noindent We now assume the situation of Definition
\ref{def:setupforsectionS99}, and consider
the implications of 
Lemma \ref{lem:antisymS99n}.

\begin{lemma}
\label{lem:dolangradyS99n}
With reference to 
Definition \ref{def:setupforsectionS99}, 
assume 
\begin{equation}
E_iA^*E_j= 0 \quad \hbox{if}\quad  i-j>1, \qquad \qquad (0 \leq i,j\leq d).
\label{eq:dolangradypreS99n}
\end{equation}
Then 
there exists scalars 
 $\beta, 
\gamma, \varrho $ in $\fld$ such that 
\begin{eqnarray}
0&=&\lbrack A, A^2A^*-\beta AA^*A+A^*A^2 -\gamma(AA^*+A^*A)
-\varrho A^*\rbrack, \qquad 
\label{eq:dolangrady1textS99n}
\end{eqnarray}
where $[r,s]$ means $rs-sr$.

\end{lemma}

\begin{proof}  
First assume $d\geq 3$.
Combining 
(\ref{eq:dolangradypreS99n}) and Lemma
\ref{lem:commentoneis1S99}(i), we obtain 
(\ref{eq:makingthingsworkS99n});
applying
Lemma \ref{lem:antisymS99n},
we find 
there exists scalars $\alpha_1, \alpha_2,\ldots,\alpha_d$ in $\fld$ such that
\begin{equation}
A^2A^*A - A A^*A^2 \;= \; \sum_{i=1}^d \alpha_i(A^iA^*-A^*A^i).
\label{eq:qaskeywilsonCS99n}
\end{equation}
We show $\alpha_i=0\;$ for $4 \leq i \leq d$. 
Suppose not, and set
\beast
t:= \hbox{max}\lbrace i \;|\;4 \leq i \leq d,\;\;\alpha_i \not=0\rbrace.
\eeast
Computing the $t0$ entry of each term in  
 (\ref{eq:qaskeywilsonCS99n}),
we readily find 
\beast
0=
\alpha_t(\theta^*_0-\theta^*_t),
\eeast
an impossibility. We  now have $\,\alpha_i=0\, $ for $4 \leq i \leq d$, 
 so (\ref{eq:qaskeywilsonCS99n}) becomes
\begin{eqnarray}
\displaystyle{
{{A^2A^*A-AA^*A^2
= \alpha_1(AA^*-A^*A)\;+\;\alpha_2(A^2A^*-A^*A^2)
}\atop
{
\qquad \qquad \qquad \qquad \qquad +\;
\alpha_3(A^3A^*-A^*A^3).}}}
\label{eq:qaskeywilsonDS99n}
\end{eqnarray}
We show $\alpha_3\not=0$. Suppose $\alpha_3=0$.
Computing the $30$ entry
in 
(\ref{eq:qaskeywilsonDS99n}),
we readily find
$\theta^*_1-\theta^*_2=0$,
an impossibility.
We have now shown $\alpha_3\not=0$.
Set 
\begin{eqnarray}
\displaystyle{
{{C:= \alpha_1 A^*\,+\,\alpha_2(AA^*+A^*A)\,+\,\alpha_3(A^2A^*+A^*A^2)
\qquad \qquad \qquad 
}\atop {\qquad  \qquad \qquad \qquad \qquad \qquad \qquad  \,+\;
(\alpha_3-1)AA^*A.}}}
\label{eq:qaskeywilsonFS99n}
\end{eqnarray}
Observe $AC-CA$ equals  
\beast
&&\alpha_1(AA^*-A^*A)\;+\;\alpha_2(A^2A^*-A^*A^2)\;+\;
\alpha_3(A^3A^*-A^*A^3)
\\
&& \qquad \qquad \qquad \qquad \qquad \qquad \;+\;
AA^*A^2-A^2A^*A,
\eeast
and this equals 0  
in view of (\ref{eq:qaskeywilsonDS99n}).
Hence $A$ and $C$ commute.
Dividing 
$C$ by
 $\alpha_3$ and using
(\ref{eq:qaskeywilsonFS99n}), 
 we find $A$ commutes with
\beast
A^2A^*-\beta AA^*A + A^*A^2 -\gamma (AA^*+A^*A)-\varrho A^*,
\eeast
where
\beast
&&\beta := \alpha_3^{-1} - 1,\qquad \gamma := -\alpha_2 \alpha_3^{-1},
\qquad \varrho := -\alpha_1\alpha_3^{-1}.
\eeast
We now have
(\ref{eq:dolangrady1textS99n})
for the case $d\geq 3$. For $d\leq 2$, we adjust our argument a bit. 
Let $\alpha_3 $ denote any nonzero scalar in $\fld$.
By our initial comments, 
and since $A^3$ is a linear
combination of $I, A,A^2$, we find there exists
scalars $\alpha_1, \alpha_2$ in $\fld$ such that
(\ref{eq:qaskeywilsonDS99n}) holds.
Proceeding as before, we obtain 
(\ref{eq:dolangrady1textS99n}).

\end{proof}

\noindent Concerning the converse to
 Lemma 
\ref{lem:dolangradyS99n}, we have the following.

\begin{lemma}
\label{lem:conclusionabttripleprodS99n}
With reference to 
Definition \ref{def:setupforsectionS99}, 
suppose there exists scalars
 $\beta, 
\gamma, \varrho $ in $\fld$ such that 
\begin{eqnarray}
0&=&\lbrack A, A^2A^*-\beta AA^*A+A^*A^2 -\gamma(AA^*+A^*A)
-\varrho A^*\rbrack. \qquad 
\label{eq:dolangrady1textS99nt}
\end{eqnarray}
Then 
\begin{eqnarray}
E_iA^*E_j &=& 0  
\quad \hbox{if} \quad 1<i-j <d, \qquad \quad 
(0 \leq i,j\leq d).
\qquad \quad 
\label{eq:mostprodzeroAS99n}
\end{eqnarray}

\end{lemma}

\begin{proof} 
We define a two variable polynomial $p \in \fld\lbrack \lambda, \mu \rbrack $
by
\beast
p(\lambda,\mu) = \lambda^2-\beta \lambda \mu +
\mu^2-\gamma (\lambda +\mu)-\varrho.
\eeast
We claim
\begin{equation}
0=E_iA^*E_jp(\theta_i,\theta_j)\quad  \hbox{if}\quad i\not=j,
\qquad \qquad (0 \leq i,j\leq d).
\label{eq:cl1convS99}
\end{equation}
To prove this, set
\beast
C &=& 
A^2A^*-\beta AA^*A+A^*A^2 -\gamma(AA^*+A^*A) -\varrho A^*,
\eeast
so that $AC=CA$.
For $0 \leq i,j\leq d$,
\beast
0 &=& E_i(AC-CA)E_j
\\
&=& (\theta_i-\theta_j)E_iCE_j,
\eeast
so if $i\not=j$,
\beast
0 &=& E_iCE_j
\\
&=&
E_iA^*E_j
(\theta_i^2-\beta\theta_i\theta_j +\theta_j^2 -
\gamma(\theta_i+\theta_j)-\varrho)
\\
&=&
E_iA^*E_jp(\theta_i,\theta_j),
\eeast
and we have 
(\ref{eq:cl1convS99}).
We next claim 
\begin{equation}
p(\theta_i,\theta_j) \not=0 \quad \hbox{if}\quad 1 < i-j< d,
\qquad \qquad (0 \leq i,j\leq d).
\label{eq:whenpisnotzeroS99}
\end{equation}
Recall $E_{i-1}A^*E_i\not=0$ for $1 \leq i \leq d$  by Lemma 
\ref{lem:commentoneis1S99}(i).
By this and 
(\ref{eq:cl1convS99}), we find $p(\theta_{i-1},\theta_i)=0$
for $1 \leq i \leq d$. Since $p$ is symmetric in its arguments,
we find $\theta_{i-1}$ and $\theta_{i+1}$ are the roots of 
$p(\lambda,\theta_i)$ for $1 \leq i \leq d-1$. By this
and since $\theta_0,\theta_1,\ldots, \theta_d$ are distinct,
we obtain
(\ref{eq:whenpisnotzeroS99}).
Combining 
(\ref{eq:cl1convS99}) and 
(\ref{eq:whenpisnotzeroS99}) we obtain
(\ref{eq:mostprodzeroAS99n}).


\end{proof}

\begin{lemma}
\label{lem:theentriesofdolangradyS99n}
With reference to 
Definition \ref{def:setupforsectionS99}, 
let $\beta, \gamma,  \varrho $ denote any
scalars in $\fld$, and consider the commutator 
\begin{equation}
\lbrack A,A^2A^*-\beta AA^*A + A^*A^2 -\gamma (AA^*+A^*A)-\varrho A^*\rbrack. 
\label{eq:entriesdolangradyS99n}
\end{equation}
Then the entries of  
(\ref{eq:entriesdolangradyS99n})
are as  follows.
\begin{enumerate}
\item The $i+1,i-2$ entry is 
\beast
\theta^*_{i-2}\;-\;(\beta+1) \theta^*_{i-1} 
\;+\;(\beta+1) \theta^*_i \;-\;\theta^*_{i+1},
\eeast
for $2 \leq i \leq d-1$.
\item The $i,i-2$ entry is 
\beast
&&
\vartheta_{i-2}\;-\;(\beta+1) \vartheta_{i-1} 
\;+\;(\beta+1) \vartheta_i \;-\;\vartheta_{i+1}
\\
&&+\;(\theta^*_{i-2}-\theta^*_0)
(\theta_{i-3}\;-\;(\beta+1) \theta_{i-2} 
\;+\;(\beta+1) \theta_{i-1} \;-\;\theta_i)
\\
&&+\;(\theta_i-\theta_d)
(\theta^*_{i-2}\;-\;(\beta+1) \theta^*_{i-1} 
\;+\;(\beta+1) \theta^*_i \;-\;\theta^*_{i+1})
\\
&&+\;(\theta^*_{i-2}-\theta^*_i)
(\theta_{i-2}\;-\;\beta \theta_{i-1} 
\;+\;\theta_i\;-\;\gamma)
\eeast
for $2 \leq i \leq d$,  where  
$\vartheta_0, \vartheta_1,\ldots, \vartheta_{d+1}$  are from
Definition \ref{def:varthetadefS99n}.
\item The $i,i-1$ entry is
\beast
&&\varphi_{i-1}(\theta_{i-2}-\beta \theta_{i-1}+\theta_i-\gamma)\;-\;
\varphi_{i+1}(\theta_{i-1}-\beta \theta_i+\theta_{i+1}-\gamma)
\\
&&+\;(\theta^*_{i-1}-\theta^*_i)
(\theta^2_{i-1}-\beta \theta_{i-1}\theta_i + \theta_i^2
-\gamma (\theta_{i-1} +\theta_i) -\varrho),
\eeast
for $1 \leq i \leq d$.
\item The $ii$ entry is 
\beast
&&\varphi_i(\theta^2_{i-1}-\beta \theta_{i-1}\theta_i + \theta_i^2
-\gamma (\theta_{i-1} +\theta_i) -\varrho)
\\
&&-\;\varphi_{i+1}(\theta^2_i-\beta \theta_i\theta_{i+1} + \theta_{i+1}^2
-\gamma (\theta_{i} +\theta_{i+1}) -\varrho),
\eeast
for $0 \leq i \leq d$.
\item The $i-1,i$ entry is 
\beast
&&\varphi_i(\theta_{i-1}-\theta_i)(\theta^2_{i-1}-\beta \theta_{i-1}\theta_i + \theta_i^2
-\gamma (\theta_{i-1} +\theta_i) -\varrho),
\eeast
for $1 \leq i \leq d$.
\end{enumerate}
All remaining entries  in 
(\ref{eq:entriesdolangradyS99n})
are zero. In the above 
formulae, we assume $\varphi_0=0$, $\varphi_{d+1}=0$,
and that $\theta_{-1}$, $\theta_{d+1}$, 
 $\theta^*_{d+1}$ 
are indeterminants.

\end{lemma}

\begin{proof} Routine matrix multiplication.

\end{proof}

\begin{theorem}
\label{lem:forsplitshapedolangpreS99n}
With reference to 
Definition \ref{def:setupforsectionS99}, 
let $\beta, \gamma, \varrho $ denote any
scalars in $\fld$. Then 
\begin{equation}
0 =\lbrack A,A^2A^*-\beta AA^*A + 
A^*A^2 -\gamma (AA^*+A^*A)-\varrho A^*\rbrack 
\label{eq:comeqzerozpreS99n}
\end{equation}
if and only if (i)--(iii) hold below.
\begin{enumerate}
\item The sequence $\theta_0,\theta_1,\ldots,\theta_d$ is
$(\beta,\gamma,\varrho)$-recurrent.
\item The sequence
$\theta^*_0,\theta^*_1,\ldots,\theta^*_d$ is $\beta$-recurrent.
\item The sequence
$\vartheta_0,\vartheta_1,\ldots,\vartheta_{d+1}$ from
Definition \ref{def:varthetadefS99n} 
is $\beta$-recurrent.
\end{enumerate}

%

\end{theorem}

\begin{proof} First assume 
(\ref{eq:comeqzerozpreS99n}).
Then 
(\ref{eq:entriesdolangradyS99n}) is zero, so all its entries
given in  
Lemma 
\ref{lem:theentriesofdolangradyS99n} are zero.
In particular, for $1 \leq i \leq d$, the expression 
in Lemma 
\ref{lem:theentriesofdolangradyS99n}(v) is zero. 
In that expression the two factors on the left   are nonzero,
so the remaining factor
\beast
\theta^2_{i-1}-\beta \theta_{i-1}\theta_i + \theta_i^2
-\gamma (\theta_{i-1} +\theta_i)  -\varrho 
\eeast
is zero.
 Now 
$\theta_0,\theta_1,\ldots,\theta_d$ is
$(\beta,\gamma,\varrho)$-recurrent by Definition
\ref{lem:beginthreetermS99}(iv). 
 For $2 \leq i\leq d-1$, the expresssion
 in 
Lemma \ref{lem:theentriesofdolangradyS99n}(i) is zero,
so 
$\theta^*_0,\theta^*_1,\ldots,\theta^*_d$ is $\beta$-recurrent
by
Definition
\ref{lem:beginthreetermS99}(ii). 
For $2 \leq i\leq d$, 
the expression in 
Lemma \ref{lem:theentriesofdolangradyS99n}(ii) is zero. 
Consider the four lines in that expression.
The sequence 
$\theta_0,\theta_1,\ldots,\theta_d$ is
$(\beta,\gamma)$-recurrent by Lemma 
\ref{lem:bgrecvsbgdrecS99}, so line 4 is zero.
The sequence 
$\theta_0,\theta_1,\ldots,\theta_d$ is
$\beta$-recurrent by Lemma 
\ref{lem:brecvsbgrecS99}, so line 2 is zero.
We mentioned 
$\theta^*_0,\theta^*_1,\ldots,\theta^*_d$ is $\beta$-recurrent,
so line 3 is zero.
Apparently line 1 is zero, so 
$\vartheta_0,\vartheta_1,\ldots,\vartheta_{d+1}$ is
$\beta$-recurrent in view of Definition 
\ref{lem:beginthreetermS99}(ii). 
We are now done in one direction.
To get the converse, suppose  (i)--(iii) hold.
Applying Lemma 
\ref{lem:brecvsbgrecS99} and Lemma 
\ref{lem:bgrecvsbgdrecS99},
 we find
 $\theta_0,\theta_1,\ldots,\theta_d$ is both
$\beta$-recurrent and 
$(\beta,\gamma)$-recurrent. From these facts
and the data in Lemma 
\ref{lem:theentriesofdolangradyS99n}, we find
all the entries of 
(\ref{eq:entriesdolangradyS99n})
 are zero, so
(\ref{eq:entriesdolangradyS99n})
 is zero. We now have
(\ref{eq:comeqzerozpreS99n}).

\end{proof}

\noindent
Modifying our point of view in 
Theorem
\ref{lem:forsplitshapedolangpreS99n}, we get the following result.

\begin{corollary}
\label{lem:forsplitshapedolangS99n}
With reference to 
Definition \ref{def:setupforsectionS99}, 
let $\beta $ denote any
scalar in $\fld$. Then there exists scalars
$\gamma, \varrho $ in $\fld$ such that
\begin{equation}
0 =\lbrack A,A^2A^*-\beta AA^*A + 
A^*A^2 -\gamma (AA^*+A^*A)-\varrho A^*\rbrack
\label{eq:comeqzerozS99n}
\end{equation}
if and only if (i)--(iii) hold below.
\begin{enumerate}
\item The sequence $\theta_0,\theta_1,\ldots,\theta_d$ is
$\beta$-recurrent.
\item The sequence
$\theta^*_0,\theta^*_1,\ldots,\theta^*_d$ is $\beta$-recurrent.
\item The sequence
$\vartheta_0,\vartheta_1,\ldots,\vartheta_{d+1}$ from
Definition \ref{def:varthetadefS99n} 
is $\beta$-recurrent.
\end{enumerate}
\end{corollary}

\begin{proof} Recall $\theta_0,\theta_1,\ldots, \theta_d$ are distinct
by
Definition \ref{def:setupforsectionS99}.
Applying Lemma
\ref{lem:brecvsbgrecS99} and Lemma 
\ref{lem:bgrecvsbgdrecS99}, we find
 $\theta_0,\theta_1,\ldots,\theta_d$ is
$\beta$-recurrent if and only if 
there exists 
 $\gamma,\varrho \in \fld$ such that
 $\theta_0,\theta_1,\ldots,\theta_d$ is
$(\beta,\gamma,\varrho)$-recurrent.
The result now follows in view of Theorem
\ref{lem:forsplitshapedolangpreS99n}.

\end{proof}

%
%

\noindent 
We now have enough information to obtain our 
classification theorem in one direction.

\begin{lemma}
\label{lem:mainthmdir1S99}
Let $\Phi$ denote a Leonard system   with 
eigenvalue sequence $\theta_0, \theta_1, \ldots, \theta_d$, 
dual eigenvalue sequence  
$\theta^*_0, \theta^*_1, \ldots, \theta^*_d$, $\varphi$-sequence
$ \varphi_1, \varphi_2, \ldots, \varphi_d $, and $\phi$-sequence
$\phi_1, \phi_2, \ldots, \phi_d$.
Then (i)--(iii) hold below.
\begin{enumerate}
\item $ {\displaystyle{ \varphi_i = \phi_1 \sum_{h=0}^{i-1}
{{\theta_h-\theta_{d-h}}\over{\theta_0-\theta_d}} 
\;+\;(\theta^*_i-\theta^*_0)(\theta_{i-1}-\theta_d) \qquad \;\;(1 \leq i \leq d)}}$,
\item $ {\displaystyle{ \phi_i = \varphi_1 \sum_{h=0}^{i-1}
{{\theta_h-\theta_{d-h}}\over{\theta_0-\theta_d}} 
\;+\;(\theta^*_i-\theta^*_0)(\theta_{d-i+1}-\theta_0) \qquad (1 \leq i \leq d)}}$,
\item The expressions
\begin{equation}
{{\theta_{i-2}-\theta_{i+1}}\over {\theta_{i-1}-\theta_i}},\qquad \qquad  
 {{\theta^*_{i-2}-\theta^*_{i+1}}\over {\theta^*_{i-1}-\theta^*_i}} 
 \qquad  \qquad 
\label{eq:defbetaplusonedir1S99}
\end{equation} 
 are equal and independent of $i$, for $\;2\leq i \leq d-1$.  
\end{enumerate}
\end{lemma}


\begin{proof} We write $\ls$, and begin by proving (iii).

\noindent
 (iii). Applying Lemma 
\ref{lem:dolangradyS99n} to the split canonical form of $\Phi$, 
we find 
there exists scalars 
 $\beta, 
\gamma, \varrho $ such that 
(\ref{eq:dolangrady1textS99n}) holds.
Applying Corollary
\ref{lem:forsplitshapedolangS99n},
we find both the eigenvalue  sequence and the dual eigenvalue
sequence of $\Phi$ are $\beta$-recurrent.
Using this, we find  
the expressions
(\ref{eq:defbetaplusonedir1S99})
equal
$\beta +1$ for $2 \leq i \leq d-1$. In particular, these 
expressions are equal and independent of $i$.

\noindent (i).
We first claim
\begin{equation}
\vartheta_i = \vartheta_1 \sum_{h=0}^{i-1}
{{\theta_h-\theta_{d-h}}\over {\theta_0-\theta_d}}
\qquad \qquad (1 \leq i \leq d),
\label{eq:mainthmAdir1S99}
\end{equation} 
\noindent where $\vartheta_1, \vartheta_2,\ldots, \vartheta_d$
are from 
Definition \ref{def:varthetadefS99n}.
To see  
(\ref{eq:mainthmAdir1S99}),
we combine Lemma 
\ref{lem:charofsumsS99}
 and Theorem 
\ref{thm:whenproductsvanishS99}.
Assume $d\geq 2$; otherwise 
(\ref{eq:mainthmAdir1S99})
is trivial.
Observe $E_dA^*E_i=0$ for $0 \leq i \leq d-2$; applying
Theorem \ref{thm:whenproductsvanishS99} to the split canonical form
of $\Phi$,
 we find
\begin{equation}
\vartheta_{i+1} = \vartheta_i\,{{\theta_i-\theta_{d-1}}\over {\theta_{i-1}-\theta_d}}+ 
\vartheta_1 \qquad \qquad (1\leq i \leq d-1).
\label{eq:almostthereinmainthmzS99n}
\end{equation}
By 
(\ref{eq:almostthereinmainthmzS99n}), 
and since $\theta_0,\theta_1,\ldots,\theta_d$
is recurrent,
 we obtain the assumptions of
Lemma 
\ref{lem:charofsumsS99} (with $r=d$). Applying  that lemma, 
we obtain 
(\ref{eq:mainthmAdir1S99}).   
In  
(\ref{eq:mainthmAdir1S99}),   
we eliminate $\vartheta_i$ on
the left
using 
Definition \ref{def:varthetadefS99n}, 
and set  $\vartheta_1=\phi_1$ on the right, to get the result.

\noindent (ii).
Apply (i) above to $\Phi^\Downarrow$, and use Theorem \ref{thm:newrepLSS99}.

\end{proof}

\noindent We finish this section by proving Theorem
\ref{eq:lastchancedolangradyS99} from the Introduction.

%

\medskip
\noindent {\it Proof of Theorem
\ref{eq:lastchancedolangradyS99}}:
Let $\Phi$ denote a Leonard system
associated with $(A,A^*)$, and
abbreviate $\alg = \alg(\Phi)$.
If the diameter $d\leq 2$, then
 any pair of multiplicity-free elements $A, A^*$ from $\alg$ satisfy
(\ref{eq:qdolangrady199n}), 
(\ref{eq:qdolangrady2S99n}), if we choose $\beta = -1$
and appropriate $\gamma$,
$\gamma^*$,
$\varrho$,  
$\varrho^*$.
For the rest of the proof,  assume $d\geq 3$.
Applying
Lemma
\ref{lem:dolangradyS99n} to the split canonical
form of $\Phi$, we find 
 there exist scalars 
$\beta, \gamma, \varrho$ in $\fld $ such that
(\ref{eq:qdolangrady199n}) holds.
Applying the above argument to $\Phi^*$,
we find there exist scalars 
$\beta^*, \gamma^*, \varrho^*$ in $\fld $ such that
\beast
0 = \lbrack A^*,A^{*2}A-\beta^* A^*AA^* + AA^{*2} -\gamma^* (A^*A+AA^*)-
\varrho^* A\rbrack.   
\eeast
We show $\beta^* = \beta $. 
Applying  
Corollary
\ref{lem:forsplitshapedolangS99n} to the split canonical
form of $\Phi$, 
we find both the eigenvalue sequence
and dual eigenvalue sequence of $\Phi$
are  $\beta$-recurrent.
From this we find 
 $\beta+1$ equals the common value  
of 
(\ref{eq:defbetaplusonedir1S99}).
Applying this argument to
 $\Phi^*$, we find $\beta^*+1$ also  equals
the common value of 
(\ref{eq:defbetaplusonedir1S99}), so 
 $\beta = \beta^*$.
We now have 
(\ref{eq:qdolangrady2S99n}).
 Concerning uniqueness, we showed   
$\beta +1$ equals the common value of 
(\ref{eq:defbetaplusonedir1S99}), so 
 $\beta $ 
is uniquely determined
by $(A,A^*)$. Applying Theorem 
\ref{lem:forsplitshapedolangpreS99n}
to the split
canonical form of $\Phi$, we find 
the eigenvalue sequence is $(\beta,\gamma,\varrho)$-recurrent,
so $\gamma$, $\varrho$ are determined by this sequence.
Applying this argument to $\Phi^*$,  
we find
the dual eigenvalue sequence is 
$(\beta,\gamma^*,\varrho^*)$-recurrent,
 so $\gamma^*$, $\varrho^*$ are
determined by this sequence. 
\par\noindent $\Box$\par

\section{Some vanishing products}

\noindent In this section, we establish  a few facts that we need
to complete the  proof of the classification theorem.

\begin{definition}
\label{def:setupforlastsec}
In this section, we assume we are in the situation of
Definition \ref{def:setupforsectionS99},  
and we further assume (i), (ii) below.
\begin{enumerate}
\item
$
{\displaystyle{
\varphi_i = \phi_1 \,\sum_{h=0}^{i-1} {{\theta_h-\theta_{d-h} }\over
{\theta_0-\theta_d}} \;+\;(\theta^*_i-\theta^*_0)(\theta_{i-1}-\theta_d)
\qquad \qquad (1 \leq i \leq d).}}$  
\item
The expressions 
\begin{equation}
{{\theta_{i-2}-\theta_{i+1}}\over {\theta_{i-1}-\theta_i}},
\qquad \qquad 
{{\theta^*_{i-2}-\theta^*_{i+1}}\over {\theta^*_{i-1}-\theta^*_i}}
\qquad \qquad 
\label{eq:equalandindS99n}
\end{equation}
are equal and independent of $i$, for $2 \leq i \leq d-1$.
\end{enumerate}
(the scalar $\phi_1$ is from 
(\ref{eq:extenddefphiS99})).

\end{definition}

\begin{lemma}
\label{lem:conclusionabttripleprodstrongS99n}
With reference to 
Definition 
\ref{def:setupforlastsec}, lines (i), (ii) hold below. 
\begin{enumerate}
\item 
$E_iA^*E_j = 0  
\quad \hbox{if} \quad i-j>1, \qquad \qquad 
(0 \leq i,j\leq d)$.
\item 
$E^*_iAE^*_j = 0 
\quad \hbox{if} \quad j-i>1, \qquad
\qquad (0 \leq i,j\leq d)$.
\end{enumerate}
\end{lemma}

\begin{proof}  
(i) We assume $d\geq 2$; otherwise there is nothing to prove.
We first claim 
\begin{equation}
\vartheta_i= 
\phi_1 \,\sum_{h=0}^{i-1} {{\theta_h-\theta_{d-h} }\over
{\theta_0-\theta_d}}
\label{eq:backtovarthetaS99}
\end{equation}
for $0 \leq i \leq d+1$,
where 
$\vartheta_0, \vartheta_1,\ldots, \vartheta_{d+1}$ are from 
Definition \ref{def:varthetadefS99n}.
To  obtain 
(\ref{eq:backtovarthetaS99}) for $1 \leq i \leq d$,  
 eliminate $\varphi_i$ in
Definition \ref{def:setupforlastsec}(i) using 
(\ref{eq:varthetadefS99n}). Line
(\ref{eq:backtovarthetaS99}) holds for   
$i=0$ and $i=d+1$, since in these cases
both sides of 
(\ref{eq:backtovarthetaS99}) are zero.
We proceed in two steps. We first  show 
\begin{eqnarray}
E_iA^*E_j &=& 0  
\quad \hbox{if} \quad 1<i-j <d, \qquad \qquad 
(0 \leq i,j\leq d). \qquad \quad 
\label{eq:mostprodzeroAsS99n}
\end{eqnarray}
To do this, we apply Lemma
\ref{lem:conclusionabttripleprodS99n}. 
By Definition
\ref{def:setupforlastsec}(ii),
there exists $\beta \in \fld$ such that both
$\theta_0,\theta_1,\ldots,\theta_d$ and
$\theta^*_0,\theta^*_1,\ldots,\theta^*_d$ are $\beta$-recurrent.
Now by 
(\ref{eq:backtovarthetaS99}) and 
Lemma \ref{lem:varthetacharacS99},
we find 
$\vartheta_0,\vartheta_1,\ldots,\vartheta_{d+1}$ is $\beta$-recurrent.
Now conditions (i)--(iii) hold in Corollary
\ref{lem:forsplitshapedolangS99n}. Applying that  corollary,
we find there exists scalars  $\gamma, \varrho $ in $\fld $ such that
(\ref{eq:comeqzerozS99n}) holds.
%
Applying
Lemma 
\ref{lem:conclusionabttripleprodS99n},
we obtain
(\ref{eq:mostprodzeroAsS99n}).
To remove the restriction $i-j<d$ in 
(\ref{eq:mostprodzeroAsS99n}), 
we show $E_dA^*E_0=0$.
By 
(\ref{eq:backtovarthetaS99}),
Lemma 
\ref{lem:shortcharofsumsAS99}(i), and
since 
$\theta_0,\theta_1,\ldots,\theta_d$ is recurrent, 
we find 
\begin{equation}
\vartheta_{i+1}= \vartheta_i
{{\theta_i-\theta_{d-1}}\over {\theta_{i-1}-\theta_d}} +\vartheta_1
\qquad \qquad (1 \leq i \leq d).
\label{eq:altrecforvarthetaS99n}
\end{equation}
In particular Theorem
\ref{thm:whenproductsvanishS99}(ii) holds.
Applying that theorem,
 we find 
$E_dA^*E_i=0$ for $0 \leq i \leq d-2$, and in particular
$E_dA^*E_0=0$.

\noindent (ii) Consider the matrices
 $A':=ZGA^*G^{-1}Z$ and $A^{*\prime}:=ZGAG^{-1}Z$ from 
Lemma \ref{lem:insteadofd4actionS99}(iii). From that lemma,
we  observe
 $A', A^{*\prime}$ satisfy the conditions  of 
Definition \ref{def:setupforsectionS99}. We show   
 $A', A^{*\prime}$ satisfy the conditions (i), (ii)
 of Definition
\ref{def:setupforlastsec}. To this end, define 
\begin{eqnarray}
&&\theta'_i := \theta^*_{d-i},\qquad \qquad \theta^{*\prime}_i := \theta_{d-i}
\qquad \qquad (0 \leq i \leq d),
\label{eq:thetaprimeS99}
\\
&&\qquad \qquad \varphi'_i := \varphi_{d-i+1}\qquad \qquad (1\leq i\leq d),
\label{eq:varphiprimeS99}
\end{eqnarray}
and put
\begin{equation}
\phi'_1 := \varphi'_1 -(\theta^{*\prime}_1-\theta^{*\prime}_0)(\theta'_0-\theta'_d)
\label{eq:phioneprimeS99}
\end{equation}
in view of (\ref{eq:extenddefphiS99}). 
We show
\begin{equation}
\varphi'_i = \phi'_1 \,\sum_{h=0}^{i-1} {{\theta'_h-\theta'_{d-h} }\over
{\theta'_0-\theta'_d}} \;+\;(\theta^{*\prime}_i-\theta^{*\prime}_0)(\theta'_{i-1}-\theta'_d)
\label{eq:varphiprimedesS99}
\end{equation}
for $1 \leq i \leq d$. Assume $d\geq 1$, and let  $i$ be given.  By Definition 
\ref{def:setupforlastsec}(i), 
\begin{equation}
\varphi_d = 
 \phi_1 +(\theta^*_d-\theta^*_0)(\theta_{d-1}-\theta_d).
\label{eq:whatisphi1okS99}
\end{equation}
Evaluating the right side of 
(\ref{eq:phioneprimeS99}) using 
(\ref{eq:thetaprimeS99}),
(\ref{eq:varphiprimeS99}),
(\ref{eq:whatisphi1okS99}), we obtain 
\begin{equation}
\phi'_1=\phi_1.
\label{eq:phiequalsphiprS99}
\end{equation}
By 
Lemma \ref{lem:dualitypreservedS99} and 
Definition \ref{def:setupforlastsec}(ii),
\begin{equation}
{{\theta_h-\theta_{d-h}}\over {\theta_0-\theta_d}}
=
{{\theta^*_h-\theta^*_{d-h}}\over {\theta^*_0-\theta^*_d}}
\qquad \qquad (0 \leq h\leq d).
\label{eq:dualitysettleS99}
\end{equation}
By 
(\ref{eq:varthetaprelims2S99}),  
(\ref{eq:thetaprimeS99}), and
(\ref{eq:dualitysettleS99}), 
\begin{equation}
\sum_{h=0}^{i-1} {{\theta'_h-\theta'_{d-h}}\over
{\theta'_0-\theta'_d}} 
\;=\;
\sum_{h=0}^{d-i} {{\theta_h-\theta_{d-h} }\over
{\theta_0-\theta_d}}. 
\label{eq:twosumsS99}
\end{equation}
The right side of 
 (\ref{eq:varphiprimedesS99}), upon simplification   
using (\ref{eq:thetaprimeS99}),
(\ref{eq:phiequalsphiprS99}),
and (\ref{eq:twosumsS99}), becomes 
\begin{equation}
\phi_1 \,\sum_{h=0}^{d-i} {{\theta_h-\theta_{d-h} }\over
{\theta_0-\theta_d}} \;+\;(\theta_{d-i}-\theta_d)(\theta^*_{d-i+1}-\theta^*_0).
\label{eq:shouldbevarphiS99}
\end{equation}
Replacing $i$ by $d-i+1$ in 
Definition \ref{def:setupforlastsec}(i), we find
(\ref{eq:shouldbevarphiS99}) equals $\varphi_{d-i+1}$. 
Recall 
$\varphi_{d-i+1}=\varphi'_i$ by
(\ref{eq:varphiprimeS99}), so
 (\ref{eq:varphiprimedesS99}) holds.
It is clear 
Definition \ref{def:setupforlastsec}(ii) holds after we replace
$\theta_j, \theta^*_j$ by 
$\theta'_j, \theta^{*\prime}_j$ for $0 \leq j\leq d$. We have now shown
$A', A^{*\prime}$ satisfy the conditions (i), (ii) of 
Definition \ref{def:setupforlastsec}, so we can apply part (i) of the 
present lemma to that pair. For $0 \leq i \leq d$, let $E'_i$
denote the primitive idempotent of $A'$ associated with $\theta'_i$,
and observe
\begin{equation}
E'_i = ZGE^*_{d-i}G^{-1}Z\qquad \qquad (0 \leq i \leq d).
\label{eq:primprimS99}
\end{equation}
By part (i) of the 
present lemma,
\begin{equation}
E'_iA^{*\prime}E'_j = 0 \quad \hbox{if} \quad i-j>1,\qquad \qquad (0 \leq i,j\leq d).
\label{eq:primeprodS99}
\end{equation}
Evaluating 
(\ref{eq:primeprodS99}) using 
(\ref{eq:primprimS99}) and the definition of $A^{*\prime}$, we obtain
\beast
E^*_{d-i}AE^*_{d-j}
= 0 \quad \hbox{if} \quad i-j>1,\qquad \qquad (0 \leq i,j\leq d).
\eeast
Replacing $i$ and $j$  in the above line by 
$d-i$ and $d-j$, respectively, we obtain
\beast
E^*_iAE^*_j
= 0 \quad \hbox{if} \quad j-i>1,\qquad \qquad (0 \leq i,j\leq d).
\eeast
\end{proof}

\begin{lemma}
\label{lem:phichangeS99}
With reference to 
Definition \ref{def:setupforlastsec}, the scalars 
$\phi_1, \phi_2,\ldots, \phi_d$ from
(\ref{eq:extenddefphiS99})
 are given by
\begin{equation}
\phi_i = \varphi_1 \,\sum_{h=0}^{i-1} {{\theta_h-\theta_{d-h} }\over
{\theta_0-\theta_d}} \;+\;(\theta^*_i-\theta^*_0)(\theta_{d-i+1}-\theta_0)
\qquad \qquad (1 \leq i \leq d).  
\label{eq:mergetwophiS99}
\end{equation}
\end{lemma}

\begin{proof} Assume $d\geq 1$; otherwise there is nothing to prove.
By
Lemma \ref{lem:dualitypreservedS99}, 
Lemma \ref{lem:factabtsumS99},
and Definition
\ref{def:setupforlastsec}(ii), 
\begin{equation}
{{\theta^*_0-\theta^*_1+\theta^*_{i-1}-\theta^*_i}\over {\theta^*_0-\theta^*_i}}
\sum_{h=0}^{i-1}\,{{\theta_h-\theta_{d-h}}\over {\theta_0-\theta_d}}
=  
{{\theta_0+\theta_{i-1}-\theta_{d-i+1}-\theta_d}\over {\theta_0-\theta_d}}.
 \label{eq:pieceofpuz4S99n}  
\end{equation}
By (\ref{eq:extenddefphiS99}),
\begin{eqnarray}
\phi_1 &=& \varphi_1 -   
 (\theta^*_1-\theta^*_0)(\theta_{0}-\theta_d).
\label{eq:pieceofpuz3S99n}
\end{eqnarray}
Adding 
(\ref{eq:extenddefphiS99}) to the equation in 
Definition \ref{def:setupforlastsec}(i), and simplifying the result
using 
(\ref{eq:pieceofpuz4S99n}),
(\ref{eq:pieceofpuz3S99n}), we routinely obtain
(\ref{eq:mergetwophiS99}).

\end{proof}


\begin{lemma}
\label{lem:whataboutphi3S99n}
With reference to 
Definition \ref{def:setupforlastsec},  suppose the scalars
$\phi_1, \phi_2,\ldots, \phi_d$ from 
(\ref{eq:extenddefphiS99})
are all nonzero. Then
the products 
 $E_{i}A^*E_{i-1} $ 
and $E^*_{i-1}AE^*_i$ 
are nonzero for 
$1 \leq i \leq d$.

\end{lemma}

\begin{proof}
 Let $V=\fld^{d+1}$ denote the irreducible left module
for $\Mdf$.
By Theorem \ref{thm:submodcondS99}, and since  each of
$\phi_1,\phi_2,\ldots, \phi_d$ is nonzero, we find
$V$ is irreducible as an $(A,A^*)$-module.
Suppose  there exists an integer $i$ $(1 \leq i \leq d)$ such 
that 
 $E_{i}A^*E_{i-1}=0 $, and consider the sum 
\begin{equation}
E_0V + E_1V + \cdots  + E_{i-1}V.
\label{eq:whatmightbemodS99}
\end{equation}
Applying Lemma 
\ref{lem:modcharS99}
to the set $S=\lbrace 0,1,\ldots, i-1\rbrace $, and using 
Lemma 
\ref{lem:conclusionabttripleprodstrongS99n}(i), we find
(\ref{eq:whatmightbemodS99}) is an $(A,A^*)$-module. The
module 
(\ref{eq:whatmightbemodS99}) is not 0 or $V$  by 
(\ref{eq:VdecompS99}), and since $1\leq i \leq d$. This
contradicts our above comment that $V$ is irreducible
as an $(A,A^*)$-module, so we conclude
 $E_{i}A^*E_{i-1}\not=0 $  for $1 \leq i \leq d$.
Next suppose 
 there exists an integer $i$ $(1 \leq i \leq d)$ such 
that 
 $E^*_{i-1}AE^*_i=0$,
  and consider the sum
\begin{equation}
E^*_iV + E^*_{i+1}V + \cdots  + E^*_dV.
\label{eq:whatmightbemodsS99}
\end{equation}
Applying Lemma 
\ref{lem:modcharstarS99}
to the set $S^*=\lbrace i,i+1,\ldots, d\rbrace $,
and using
Lemma
\ref{lem:conclusionabttripleprodstrongS99n}(ii), we find
(\ref{eq:whatmightbemodsS99})
 is an $(A,A^*)$-module. 
We observe (\ref{eq:whatmightbemodsS99}) is not $0$ or $V$,
contradicting the fact that $V$ is irreducible
as an $(A,A^*)$-module. We conclude
 $E^*_{i-1}AE^*_i\not=0$ for
 $1 \leq i \leq d$.

\end{proof}

\section{A classification of Leonard systems}

\noindent We are now ready to prove our classification
theorem for Leonard systems, which is Theorem
\ref{thm:newrepLScharagainS99} from the Introduction.

\medskip
\noindent {\it  Proof of Theorem \ref{thm:newrepLScharagainS99}}:
To prove the theorem in one direction, 
let
\beast
\ls 
\eeast
denote a Leonard system over $\fld $ with
eigenvalue sequence $\;\theta_0, \theta_1, \ldots, \theta_d$, 
dual eigenvalue sequence  
$\;\theta^*_0, \theta^*_1, \ldots, \theta^*_d$, $\varphi$-sequence
$\; \varphi_1, \varphi_2, \ldots, \varphi_d $, and $\phi$-sequence
$\;\phi_1, \phi_2, \ldots, \phi_d$.
We verify conditions (i)--(v) in the statement of the theorem.
Condition (i) holds by 
Definition
\ref{def:varphidefS99},
Definition
\ref{def:phidefS99}, and the last assertion of 
Theorem
\ref{thm:splitformexistS99}. 
Condition (ii) holds by
Definition \ref{def:lsdefcommentsS99},
and since $A$ and $A^*$ are multiplicity-free.
Conditions (iii)--(v) are immediate from
Lemma  \ref{lem:mainthmdir1S99}, and we are done in one direction.

\medskip
\noindent To obtain the converse, suppose (i)--(v) hold
in the present theorem, and put 
\beast
A = 
\left(
\begin{array}{c c c c c c}
\theta_0 & & & & & {\bf 0} \\
1 & \theta_1 &  & & & \\
& 1 & \theta_2 &  & & \\
& & \cdot & \cdot &  &  \\
& & & \cdot & \cdot &  \\
{\bf 0}& & & & 1 & \theta_d
\end{array}
\right),
&&\quad 
A^* = 
\left(
\begin{array}{c c c c c c}
\theta^*_0 &\varphi_1 & & & & {\bf 0} \\
 & \theta^*_1 & \varphi_2 & & & \\
&  & \theta^*_2 & \cdot & & \\
& &  & \cdot & \cdot &  \\
& & &  & \cdot & \varphi_d \\
{\bf 0}& & & &  & \theta^*_d
\end{array}
\right).
\eeast
We observe $A$ (resp. $A^*$) 
is multiplicity-free, with eigenvalues
$\theta_0,\theta_1,\ldots, \theta_d$,
(resp. $\theta^*_0,\theta^*_1,\ldots, \theta^*_d$).
For $0 \leq i \leq d$, let $E_i$ (resp. $E^*_i$) denote
the primitive idempotent of $A$ (resp. $A^*$) associated
with $\theta_i$ (resp. $\theta^*_i$).
We show
\begin{equation}
\Phi:=(A;E_0,E_1,\ldots, E_d;A^*;E^*_0,E^*_1,\ldots,E^*_d)
\label{eq:isthisleonardsysS99}
\end{equation}
is a Leonard system in $\Mdf$.
To do this, we show $\Phi$ satisfies the conditions (i)--(v)
of Definition
\ref{def:deflstalkS99}.
Conditions (i)--(iii) are clearly satisfied,
so consider conditions
 (iv), (v).
By Lemma 
\ref{lem:commentoneis1S99}, 
\begin{eqnarray}
 E_iA^*E_j = \left\{ \begin{array}{ll}
                   0  & \mbox{if $\;j-i > 1; $ } \\
				  \not= 0 & \mbox{if $\;j-i=1 $ }
			   \end{array}
			\right. \qquad \qquad  (0 \leq i,j\leq d)
\label{eq:thirdoffourpartsAS99}
\end{eqnarray}
and 
\begin{eqnarray}
 E^*_iAE^*_j = \left\{ \begin{array}{ll}
                   0  & \mbox{if $\;i-j > 1; $ } \\
				  \not= 0 & \mbox{if $\;i-j=1 $ }
			   \end{array}
			\right. \qquad \qquad  (0 \leq i,j\leq d).
\label{eq:fourthoffourpartsAS99}
\end{eqnarray}
By assumption (iv), the scalar $\phi_1$ in the present theorem
equals
\beast
 \varphi_1 - (\theta^*_1-\theta^*_0)(\theta_0-\theta_d),
\eeast
and  therefore equals the scalar denoted $\phi_1$ in
Definition \ref{def:extenddefphiS99}.
Combining this with assumptions (iii), (v) in
the present theorem, 
we find $A$ and $A^*$ satisfy conditions (i), (ii) of
Definition
\ref{def:setupforlastsec}. Applying 
Lemma \ref{lem:conclusionabttripleprodstrongS99n},
we find
\begin{eqnarray}
&&E_iA^*E_j = 0 \quad \hbox{if}\quad i-j>1, \qquad \quad (0 \leq i,j\leq d),
\qquad \qquad 
\label{eq:firstoffourpartsS99}
\\
&&E^*_iAE^*_j = 0 \quad \hbox{if}\quad j-i>1, \qquad \quad 
(0 \leq i,j\leq d).
\label{eq:secondoffourpartsS99}
\end{eqnarray}
By assumption (iv) and Lemma 
\ref{lem:phichangeS99}, 
the sequence
$\phi_1,\phi_2,\ldots, \phi_d$ in the present theorem
equals the corresponding 
 sequence from 
Definition \ref{def:extenddefphiS99}.
The elements of this 
sequence are nonzero
by assumption (i), 
 so by Lemma 
\ref{lem:whataboutphi3S99n},
\begin{eqnarray}
&&E_iA^*E_j \not= 0 
\quad \hbox{if}\quad i-j=1 \qquad \qquad (0 \leq i,j\leq d),
\label{eq:fifthoffourpartsS99}
\\
&&E^*_iAE^*_j \not= 0 
\quad \hbox{if}\quad j - i=1 \qquad \qquad (0 \leq i,j\leq d).
\label{eq:sixthoffourpartsS99}
\end{eqnarray}
Combining
(\ref{eq:thirdoffourpartsAS99})--(\ref{eq:sixthoffourpartsS99}),
we obtain conditions (iv), (v) of 
Definition
\ref{def:deflstalkS99}, so 
$\Phi$ is a Leonard system in $\Mdf$.
By Lemma
\ref{lem:def41andlsS99}, we find
$\Phi$ has eigenvalue sequence
$\theta_0,\theta_1,\ldots,\theta_d$,  dual eigenvalue  
sequence 
$\theta^*_0,\theta^*_1,\ldots,\theta^*_d$,
and $\varphi$-sequence  
$\varphi_1, \varphi_2,\ldots, \varphi_d$.
We mentioned the sequence 
$\phi_1, \phi_2,\ldots, \phi_d$ from the present theorem
is the same as the corresponding sequence from
Definition \ref{def:extenddefphiS99},
so
this  
is the $\phi$-sequence of $\Phi$ in view of 
Lemma
\ref{lem:varphiminusphiS99}(i).
The Leonard system
$\Phi$ is unique up to isomorphism by Lemma 
\ref{lem:paramsdetisoS99}.
\par\noindent $\Box$\par



\begin{corollary}
\label{thm:wegetallfrom9S99}
Let $\Phi$ denote a Leonard system over $\fld $ with diameter $d\geq 3$, 
eigenvalue sequence $\theta_0, \theta_1, \ldots, \theta_d$, 
dual eigenvalue sequence  
$\theta^*_0, \theta^*_1, \ldots, \theta^*_d$, $\varphi$-sequence
$ \varphi_1, \varphi_2, \ldots, \varphi_d $, and $\phi$-sequence
$\phi_1, \phi_2, \ldots, \phi_d$.
Consider a sequence $\cal S$ of 9 parameters consisting of the
sequence in (i) below, followed by either parameter in (ii)
below, followed by any one of the parameters in (iii) below:
\begin{enumerate}
\item $d,\theta_0,\theta_1, \theta_2,\theta^*_0,\theta^*_1,\theta^*_2$,
\item $\theta_3, \theta^*_3$,
\item $\varphi_1, \phi_1, \varphi_d,\phi_d $.
\end{enumerate}
Then
the isomorphism class of $\Phi$ as a Leonard system over $\fld$
is determined by $\cal S$.

\end{corollary}

\begin{proof} 
From  Theorem
\ref{thm:newrepLScharagainS99}(v), 
we recursively
obtain $\theta_i$, $\theta^*_i$ for $0 \leq i \leq d$.
Using Theorem \ref{thm:newrepLScharagainS99}(iii),(iv) (with $i=1$, $i=d$),
we obtain $\phi_1$. 
Using  
Theorem \ref{thm:newrepLScharagainS99}(iii), we obtain
$\varphi_i$  for $1\leq i \leq d$.
Applying 
Lemma 
\ref{lem:paramsdetisoS99},
 we find 
the isomorphism class of $\Phi$ is determined by $\cal S$.

\end{proof}

\begin{corollary}
\label{cor:lastcommentonsplitS99}
Let 
$d$ denote a nonnegative integer, let $\fld $ denote a field,
and let 
$A$ and $A^*$ denote  matrices in $\Mdf$ of the form 
\beast
A = 
\left(
\begin{array}{c c c c c c}
\theta_0 & & & & & {\bf 0} \\
1 & \theta_1 &  & & & \\
& 1 & \theta_2 &  & & \\
& & \cdot & \cdot &  &  \\
& & & \cdot & \cdot &  \\
{\bf 0}& & & & 1 & \theta_d
\end{array}
\right),
&&\quad 
A^* = 
\left(
\begin{array}{c c c c c c}
\theta^*_0 &\varphi_1 & & & & {\bf 0} \\
 & \theta^*_1 & \varphi_2 & & & \\
&  & \theta^*_2 & \cdot & & \\
& &  & \cdot & \cdot &  \\
& & &  & \cdot & \varphi_d \\
{\bf 0}& & & &  & \theta^*_d
\end{array}
\right).
\eeast

Then the following are equivalent.
\begin{enumerate}
\item $(A,A^*)$ is a Leonard pair  in $\Mdf$.
\item There exist a sequence of scalars $\phi_1, \phi_2, \ldots, \phi_d$ taken
from $\fld$
such that (i)--(v) hold in Theorem
\ref{thm:newrepLScharagainS99}.
\end{enumerate}
Suppose (i),(ii) hold above. Then
\begin{equation}
(A;E_0,E_1,\ldots, E_d;A^*;E^*_0,E^*_1,\ldots, E^*_d)
\label{eq:whatislslastqS99}
\end{equation}
is a Leonard system in $\Mdf$, where 
$E_i$ (resp. $E^*_i$) denotes the primitive idempotent
of $A$ (resp. $A^*$) associated with $\theta_i$ (resp. $\theta^*_i$),
for $0 \leq i \leq d$. The Leonard system
(\ref{eq:whatislslastqS99})
has eigenvalue sequence 
    $\theta_0, \theta_1,\ldots,
\theta_d$,
dual eigenvalue sequence 
 $\theta^*_0, \theta^*_1,\ldots,
\theta^*_d$,
$\varphi$-sequence 
 $\varphi_1, \varphi_2,\ldots, \varphi_d$, 
and $\phi$-sequence  
 $\phi_1, \phi_2,\ldots, \phi_d$. 

\end{corollary}

\begin{proof} $(i)\rightarrow (ii)$.
We first show
(\ref{eq:whatislslastqS99}) is a Leonard system.
To do this, we apply Lemma
\ref{lem:def41andlsS99}. In order to do that, we 
verify $A$ and $A^*$ satisfy the conditions of
Definition
\ref{def:setupforsectionS99}.
Certainly $\theta_0, \theta_1,\ldots, \theta_d$ are distinct,
since $A$ is multiplicity-free.
Similarily
$\theta^*_0, \theta^*_1,\ldots, \theta^*_d$ are distinct.
Observe $\varphi_1,\varphi_2,\ldots, \varphi_d$ are nonzero;
otherwise the left module $V=\fld^{d+1}$ of 
$\Mdf$ 
is reducible as an $(A,A^*)$-module,
contradicting Lemma
\ref{lem:VirredaastarmoduleS99}.
We have now shown 
 $A$ and $A^*$ satisfy  the conditions of 
Definition 
\ref{def:setupforsectionS99}, so we can apply
Lemma 
\ref{lem:def41andlsS99}.  
By that lemma,
we find
(\ref{eq:whatislslastqS99}) is a Leonard system,
with 
   eigenvalue sequence $\theta_0, \theta_1,\ldots,
\theta_d$,  
dual eigenvalue sequence 
 $\theta^*_0, \theta^*_1,\ldots,
\theta^*_d$, and 
$\varphi$-sequence 
 $\varphi_1, \varphi_2,\ldots,
\varphi_d$. 
Let $\phi_1, \phi_2,\ldots \phi_d$
denote the $\phi$-sequence of the Leonard system
(\ref{eq:whatislslastqS99}).
Applying Theorem
\ref{thm:newrepLScharagainS99} to this system, we find 
(i)--(v) hold in that theorem.

\noindent 
$(ii)\rightarrow (i)$.
By 
Theorem
\ref{thm:newrepLScharagainS99},  
there exists a Leonard system 
$\Phi$   
over $\fld$ with eigenvalue
sequence
 $\theta_0, \theta_1,\ldots,
\theta_d$,
dual eigenvalue
sequence
 $\theta^*_0, \theta^*_1,\ldots,
\theta^*_d$,
$\varphi$-sequence 
 $\varphi_1, \varphi_2,\ldots,
\varphi_d$, and $\phi$-sequence 
 $\phi_1, \phi_2,\ldots,
\phi_d$. 
The Leonard system $\Phi$ has split canonical form
(\ref{eq:whatislslastqS99})  by  
Definition
\ref{def:varphidefS99}, so 
(\ref{eq:whatislslastqS99})  is a Leonard system  in $\Mdf$. 
In particular
$(A,A^*)$ is a Leonard pair in $\Mdf$, as desired.

\noindent Suppose (i),(ii). From the proof of
 $(ii)\rightarrow (i)$ we find  
(\ref{eq:whatislslastqS99}) is a Leonard system in $\Mdf$,
with the required
eigenvalue, dual eigenvalue, $\varphi$- and $\phi$-sequences.

\end{proof}

\section{Appendix: Leonard systems and polynomials}

\noindent There is a theorem due to Doug Leonard 
 \cite{Leodual},
 \cite[p260]{BanIto} 
that gives a characterization of the $q$-Racah polynomials
and some related polynomials in the Askey scheme
\cite{Ask}, \cite{AskAW}, \cite{ARSClass},
\cite{KoeSwa},
\cite{KooGroup}.
The situation  considered in that theorem
is closely connected 
to the subject of the present paper,
and it is this connection that motivates our terminology. 
We sketch the connection here without proof; details will be
provided in a future paper.


\medskip
\noindent
Let $\fld$ denote any field, and let
\begin{equation}
\ls 
\label{eq:deflsinlastchS99}
\end{equation}
denote a Leonard system over $\fld$.
Then there
exists a unique sequence of monic polynomials 
\beast
 p_0, p_1, \ldots, p_{d+1};\qquad \qquad 
p^*_0, p^*_1, \ldots, p^*_{d+1}  
\eeast
 in $\fld\lbrack \lambda \rbrack $ such that  
\beast
   deg\, p_i = i,\qquad &&\qquad  
\qquad deg \, p^*_i = i \qquad \qquad  \quad (0 \leq i \leq d+1), 
\\
p_i(A)E^*_0 =E^*_iA^iE^*_0,&& \qquad   
p^*_i(A^*)E_0 =E_iA^{*i}E_0 \qquad \qquad  (0\leq i \leq d), \qquad \qquad  
\\
p_{d+1}(A)=0,\quad \;\; &&\qquad  
\quad \;\;p^*_{d+1}(A^*)=0. 
\eeast
These polynomials satisfy  
\begin{eqnarray}
&& \qquad p_0=1, \qquad \qquad p^*_0=1, 
\label{eq:rec1introS99}
\\
\lambda p_i &=& p_{i+1} + a_ip_i +x_ip_{i-1} \qquad \qquad (0 \leq i \leq d),
\label{eq:rec2introS99}
\\
\lambda p^*_i &=& p^*_{i+1} + a^*_ip^*_i +x^*_ip^*_{i-1} 
\qquad \qquad (0 \leq i \leq d),
\label{eq:rec3introS99}
\end{eqnarray}
where $\,x_0,\; x^*_0, \; 
p_{-1},\; 
p^*_{-1}\,$ are all zero,  and where 
\beast 
 a_i =tr \,E^*_iA,\quad && \qquad    \quad   
a^*_i =tr \,E_iA^* \qquad \qquad \quad  (0 \leq i \leq d), \qquad 
\\
x_i =tr \,E^*_iAE^*_{i-1}A,&& \qquad     
x^*_i =tr\, E_iA^*E_{i-1}A^* \qquad  \quad   (1 \leq i \leq d). 
\eeast
In fact
\begin{equation}
x_i\not=0, \qquad \qquad 
x^*_i\not=0 \qquad \qquad (1 \leq i \leq d).
\label{eq:rec4introS99}
\end{equation}
We call $p_0, p_1,\ldots, p_{d+1}$  the 
{\it monic polynomial sequence} (or {\it MPS}) of $\Phi$.
We call 
$p^*_0, p^*_1,\ldots, p^*_{d+1}$  the  {\it dual MPS} of  $\Phi$.
Let $\theta_0,\theta_1,\ldots, \theta_d$ 
(resp. $\theta^*_0,\theta^*_1,\ldots,\theta^*_d$)
denote the eigenvalue sequence (resp. dual eigenvalue sequence)
of $\Phi$, so that
\begin{eqnarray}
&&\theta_i\not=\theta_j, \qquad \quad 
\theta^*_i\not=\theta^*_j \qquad \quad \hbox{if} \quad i\not=j,  
\qquad \;\;  (0 \leq i,j\leq d), \qquad  \qquad 
\label{eq:rec5introS99}
\\
&&\qquad p_{d+1}(\theta_i) = 0, \qquad \qquad 
p^*_{d+1}(\theta^*_i) = 0, \qquad  (0 \leq i \leq d). \qquad 
\label{eq:rec6introS99}
\end{eqnarray}
Then
\begin{equation} 
p_i(\theta_0)\not=0, \qquad \qquad  
p^*_i(\theta^*_0)\not=0 \qquad \qquad  (0 \leq i \leq d),
\label{eq:rec7introS99}
\end{equation}
and  
\begin{equation}
{{p_i(\theta_j)}\over {
p_i(\theta_0)}}
= 
{{p^*_j(\theta^*_i)}\over {
p^*_j(\theta^*_0)}}\qquad \qquad (0 \leq i,j\leq d). 
\label{eq:rec8introS99}
\end{equation}

\medskip
\noindent
Conversely, given
 polynomials
\beast
\qquad \qquad \qquad \qquad 
  p_0, p_1, \ldots, p_{d+1}; \qquad \qquad
p^*_0, p^*_1, \ldots, p^*_{d+1}
\qquad \qquad \qquad \qquad \qquad  \; \;
(174L),\; (174R)
\eeast
in 
 $\fld\lbrack \lambda \rbrack $ satisfying  
(\ref{eq:rec1introS99})--(\ref{eq:rec4introS99}), 
and given   
 scalars 
\beast
\qquad \qquad \qquad \qquad 
\theta_0, \theta_1, \ldots, \theta_d; \qquad \qquad  \qquad
\theta^*_0, \theta^*_1, \ldots, \theta^*_d 
\qquad \qquad \qquad \qquad \qquad \; \;(175L),\;(175R)
\eeast
in 
 $\fld $ 
satisfying 
(\ref{eq:rec5introS99})--(\ref{eq:rec8introS99}),
 there exists  
 a Leonard system  $\Phi$ over $\fld$
with MPS (174L), dual MPS (174R), eigenvalue sequence (175L), 
and dual eigenvalue sequence (175R).
The system $\Phi $ is unique up to isomorphism of Leonard systems.   

\medskip
\noindent In the above paragraph, we described a bijection 
between the Leonard systems and the systems 
(174L)--(175R)  satisfying 
(\ref{eq:rec1introS99})--(\ref{eq:rec8introS99}).
For the case
 $\fld = \R$,  
the systems
(174L)--(175R)  satisfying 
(\ref{eq:rec1introS99})--(\ref{eq:rec8introS99}) were previously
classified by Leonard
\cite{Leodual} and Bannai and Ito
 \cite[p260]{BanIto}. 
They found the polynomials involved 
are $q$-Racah polynomials or
related polynomials from the Askey scheme.
Their classification has come to be known as Leonard's theorem.
Given the above bijection, we may view 
Theorem \ref{thm:newrepLScharagainS99} in the present paper
as a ``linear algebraic version'' of Leonard's theorem. To see one 
advantage of our version, compare it with the previous version
 \cite[p260]{BanIto}. In that version, the statement of the theorem
 takes 11 pages. 
 We believe the main   value of our version lies in the conceptual
 simplicity and alternative point of view it provides for the study
 of orthogonal polynomials.

%

\medskip
\noindent For the benefit of researchers in special functions 
we now give  more detail on the polynomials 
that come from  Leonard systems.  In what follows, we freely
use the notation of 
Definition
\ref{def:d4actionlsS99}.
Let $\Phi $ denote the Leonard system 
from
(\ref{eq:deflsinlastchS99}).
In view of 
(\ref{eq:rec8introS99}) we define  the polynomials
\beast
u_i = {{p_i}\over {p_i(\theta_0)}} \qquad \qquad (0 \leq i \leq d),
\eeast
so that
\beast
u_i(\theta_j) = u^*_j(\theta^*_i) \qquad \qquad (0 \leq i,j\leq d).
\eeast
The $u_0, u_1, \ldots, u_d$ satisfy a  recurrence similar
to 
(\ref{eq:rec2introS99}), as we now explain.
There exists
a unique sequence of scalars $c_1, c_2, \ldots, c_d;$
 $b_0, b_1,\ldots, b_{d-1}$ taken from $\fld$  such that 
\beast
x_i &=& b_{i-1}c_i \qquad \qquad \qquad (1\leq i \leq d),
\\
\theta_0 &=& c_i+a_i+b_i \qquad \qquad (0 \leq i\leq d),
\eeast
where $c_0=0, b_d=0$. 
Then
\beast
\lambda u_i &=& c_iu_{i-1}+a_iu_i+b_iu_{i+1} \qquad \qquad (0 \leq i \leq d-1), 
\eeast
and  $\; 
\lambda u_d - c_du_{d-1}-a_du_d \;$ vanishes on  each of 
$\theta_0, \theta_1, \ldots, \theta_d$. 

\medskip
\noindent The polynomials $p_i$ and $u_i$ both  satisfy 
 orthogonality relations.
Set  
\beast
m_i:= tr\,E_iE^*_0 \qquad \qquad (0 \leq i \leq d).
\eeast
Then each of $m_0, m_1, \ldots, m_d$ is nonzero, and  
the orthogonality for the  $p_i$ is
\beast
 \sum_{r=0}^d p_i(\theta_r)p_j(\theta_r)m_r &=& \delta_{ij}
x_1x_2\cdots x_i \qquad \qquad (0 \leq i,j\leq d),
\\
\sum_{i=0}^d {{ p_i(\theta_r) p_i(\theta_s)}\over {x_1x_2\cdots x_i}}
&=& \delta_{rs}m_r^{-1} \qquad \qquad (0 \leq r,s\leq d).
\eeast
We remark that since the ground field  $\fld$ is arbitrary, the 
question of whether the $m_i$ are positive or not does not arise.

\medskip
\noindent 
Turning to the  $u_i$, observe
that 
 $m_0=m^*_0$; let $n$ denote the
multiplicative inverse of this common value, and set
\beast
k_i:=m^*_in 
\qquad \qquad  (0 \leq i \leq d).
\eeast
The orthogonality for the  $u_i$ is
\beast
\sum_{r=0}^d u_i(\theta_r)u_j(\theta_r)m_r &=&
\delta_{ij}k^{-1}_i \qquad \qquad (0 \leq i,j\leq d),
\\
\sum_{i=0}^d u_i(\theta_r)u_i(\theta_s)k_i &=&
\delta_{rs} m_r^{-1}   \qquad \qquad (0 \leq r,s\leq d).
\eeast
We remark 
\beast
k_i={{b_0b_1\cdots b_{i-1}}\over {c_1c_2\cdots c_i}}
\qquad \qquad 
(0 \leq i \leq d)
\eeast
and 
\beast
n=k_0+k_1+\cdots+k_d.
\eeast
All the polynomials and scalars we have described   in this section
are given by 
fairly simple
rational expressions involving  the eigenvalues, dual eigenvalues, 
$\varphi$-sequence and $\phi$-sequence of $\Phi$. 
For example the $a_i, a^*_i$ are 
given in 
Lemma 
\ref{thm:aiintermsofphivarphiS99}.
Using the notation
of Definition \ref{def:thetauS99}, we have
\beast
b_i = \varphi_{i+1} {{\tau^*_i(\theta^*_i)}\over {
\tau^*_{i+1}(\theta^*_{i+1})}}, \qquad \qquad 
b^*_i = \varphi_{i+1} {{\tau_i(\theta_i)}\over {
\tau_{i+1}(\theta_{i+1})}} 
\eeast
for $0 \leq i \leq d-1$,   
\beast
c_i = \phi_i {{\eta^*_{d-i}(\theta^*_i)}\over{\eta^*_{d-i+1}(\theta^*_{i-1})}},
\qquad \qquad 
c^*_i = \phi_{d-i+1} {{\eta_{d-i}(\theta_i)}\over{\eta_{d-i+1}(\theta_{i-1})}}
\eeast
for $1\leq i\leq d$, and 
\beast
n = {{\eta_d(\theta_0) \eta^*_d(\theta^*_0)
}\over {\phi_1 \phi_2 \cdots \phi_d 
}}.
\eeast
Moreover
\beast
p_i=\sum_{h=0}^i {{\varphi_1\varphi_2\cdots \varphi_i}\over{
\varphi_1\varphi_2\cdots \varphi_h}}
{{\tau^*_h(\theta^*_i)}\over 
{\tau^*_i(\theta^*_i)}} \tau_h,
\qquad \qquad 
p^*_i=\sum_{h=0}^i {{\varphi_1\varphi_2\cdots \varphi_i}\over{
\varphi_1\varphi_2\cdots \varphi_h}}
{{\tau_h(\theta_i)}\over 
{\tau_i(\theta_i)}} \tau^*_h,
\\
\eeast
\begin{eqnarray}
\label{eq:uifinal}
u_i=\sum_{h=0}^i {{\tau^*_h(\theta^*_i)}\over 
{\varphi_1\varphi_2\cdots \varphi_h}} \tau_h 
,\qquad \qquad 
u^*_i=\sum_{h=0}^i {{\tau_h(\theta_i)}\over 
{\varphi_1\varphi_2\cdots \varphi_h}} \tau^*_h 
\end{eqnarray}
for  $0 \leq i \leq d$. Immediately after Theorem
\ref{thm:newrepLScharagainS99} in the Introduction, we displayed 
a parametric solution. We now assume the parameters of $\Phi$
are given by this solution, and consider the effect on the
associated polynomials.
Using this solution and either equation in 
(\ref{eq:uifinal}), we find that for $0 \leq i,j\leq d$, 
the common value of 
$u_i(\theta_j), u^*_j(\theta^*_i)$ is given by
\begin{eqnarray}
\sum_{n=0}^d {{(q^{-i};q)_n (s^*q^{i+1};q)_n 
(q^{-j};q)_n (sq^{j+1};q)_n q^n
}\over
{(r_1q;q)_n(r_2q;q)_n (q^{-d};q)_n(q;q)_n }}, 
\label{eq:uihyper}
\end{eqnarray}
where 
\beast
(a;q)_n := (1-a)(1-aq)(1-aq^2)\cdots (1-aq^{n-1})\qquad \qquad n=0,1,2\ldots 
\eeast
Observe   
(\ref{eq:uihyper}) is the basic hypergeometric series
\beast
 {}_4\phi_3 \Biggl({{q^{-i}, \;s^*q^{i+1},\;q^{-j},\;sq^{j+1}}\atop
{r_1q,\;\;r_2q,\;\;q^{-d}}};\; q,\;q\Biggr).
\eeast
The  $q$-Racah polynomials are defined in 
\cite{KoeSwa}. Comparing that definition with the above data, and
recalling $r_1r_2=s s^*q^{d+1}$, 
we find the $u_i$ and the $u^*_i$ are $q$-Racah polynomials. 

\medskip



\end{document}